\documentclass[11pt]{amsart}
\usepackage{graphicx,rlepsf,latexsym,amssymb}

\usepackage[all]{xy}
\CompileMatrices

\begin{document}
%
%
%
%
%
%

%
% Definition of theorem environments
%

\theoremstyle{plain}
\newtheorem{theorem}{Theorem}[section]
\newtheorem{lemma}[theorem]{Lemma}
\newtheorem{proposition}[theorem]{Proposition}
\newtheorem{corollary}[theorem]{Corollary}
\newtheorem{definition}[theorem]{Definition}

\theoremstyle{definition}
\newtheorem{example}[theorem]{Example}
\newtheorem{remark}[theorem]{Remark}
\newtheorem{summary}[theorem]{Summary}
\newtheorem{notation}[theorem]{Notation}
\newtheorem{problem}[theorem]{Problem}

\theoremstyle{remark}
\newtheorem{claim}[theorem]{Claim}
\newtheorem{sublemma}[theorem]{Sub-lemma}
\newtheorem{innerremark}[theorem]{Remark}

\numberwithin{equation}{section}
\numberwithin{figure}{section}

%
% Preferences
%

\renewcommand{\labelitemi}{$\centerdot$}

%
% Definition of new commands
%

\def \N{\mathbb{N}}
\def \Z{\mathbb{Z}}
\def \Q{\mathbb{Q}}
\def \R{\mathbb{R}}

\def \A{\mathcal{A}}
\def \As{\mathcal{A}^s}
\def \AY{\mathcal{A}^Y}
\def \AYc{\mathcal{A}^{Y,c}}
\def \btT{{}_{b}^t\mathcal{T}}
\def \Cob{\mathcal{C}ob}
\def \Cyl{\mathcal{C}yl}
\def \F{\mathcal{F}}
\def \I{\mathcal{I}}
\def \LCob{\mathcal{LC}ob}
\def \dLCob{{}^d\!\mathcal{LC}ob}
\def \LqCob{\mathcal{LC}ob_q}
\def \dLqCob{{}^d\!\mathcal{LC}ob_q}
\def \M{\mathcal{M}}
\def \P{\mathcal{P}}
\def \qTCub{\mathcal{T}_q\mathcal{C}ub}
\def \QLCob{\Q\mathcal{LC}ob}
\def \QLqCob{\Q\mathcal{LC}ob_q}
\def \sLCob{{}^s\!\mathcal{LC}ob}
\def \sLqCob{{}^s\!\mathcal{LC}ob_q}
\def \T{\mathcal{T}}
\def \tsA{{}^{ts}\!\!\mathcal{A}}
\def \Ztilde{\widetilde{Z}}
\def \ZtildeY{\widetilde{Z}^{Y}}
\def \ZZ{\widetilde{\mathsf{Z}}}

\def \osqcup{\hphantom{}^<_\sqcup}

\def \aug{{\rm aug}}
\def \Aut{{\rm Aut}}
\def \cl{{\rm cl}}
\def \Coker{{\rm Coker}}
\def \deg{{\rm deg}}
\def \ideg{{\rm i\hbox{-}deg}}
\def \edeg{{\rm e\hbox{-}deg}}
\def \End{{\rm End}}
\def \Gr{{\rm Gr}}
\def \Hom{{\rm Hom}}
\def \incl{{\rm incl}}
\def \Id{{\rm Id}}
\def \Img{{\rm Im}}
\def \Ker{{\rm Ker}}
\def \Lk{{\rm Lk}}
\def \mod{{\rm mod}}
\def \ord{{\rm ord}}
\def \pr{{\rm pr}}
\def \rk{{\rm rk}}
\def \sgn{{\rm sgn}}
\def \Tors{{\rm Tors}}
\def \ud{{\rm d}}

\newcommand{\set}[1]{\lfloor #1\rceil}
\newcommand{\Star}[3]{\stackrel{{\scriptsize #1},{\scriptsize #2}}{\star}}

\newcommand{\figtotext}[3]{\begin{array}{c}\includegraphics[width=#1pt,height=#2pt]{#3}\end{array}}

\newcommand{\capleft}{\figtotext{10}{10}{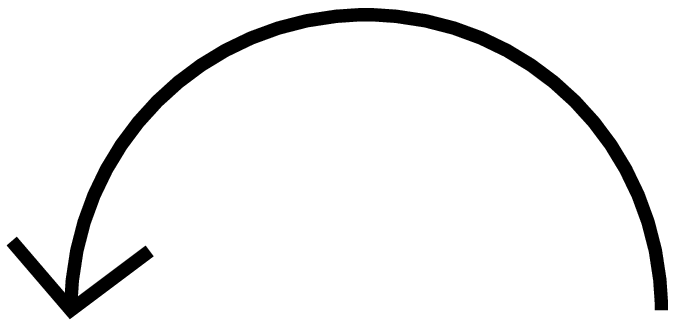}}
\newcommand{\cupright}{\figtotext{10}{10}{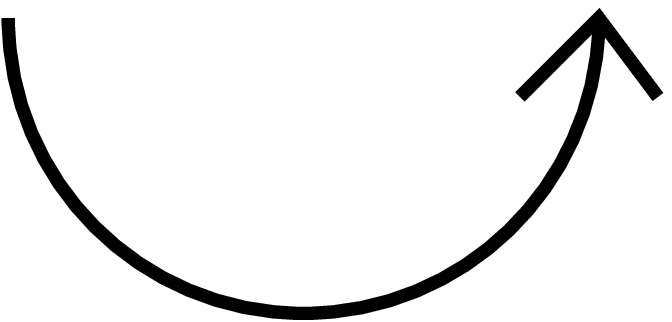}}
\newcommand{\crossing}{\figtotext{12}{12}{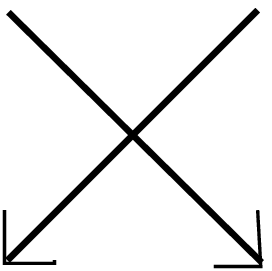}}
\newcommand{\pcrossing}{\figtotext{12}{12}{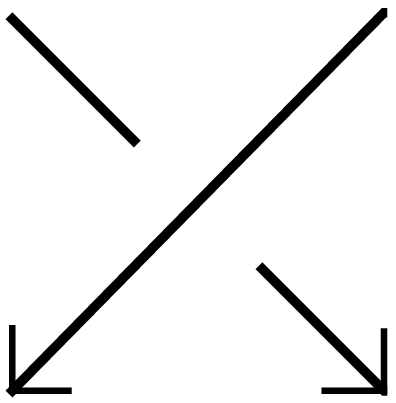}}
\newcommand{\ncrossing}{\figtotext{12}{12}{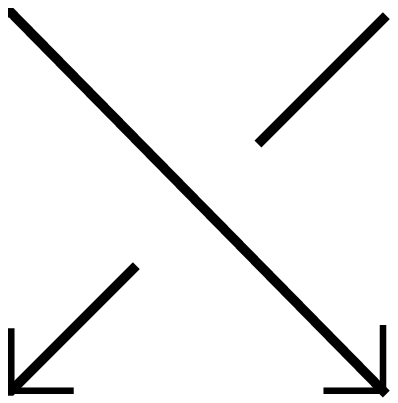}}

\newcommand{\thetagraph}
{\hspace{-0.2cm} \figtotext{12}{12}{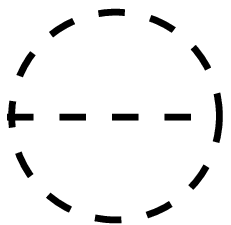} \hspace{-0.2cm}}

\newcommand{\strutgraph}[2]
{\begin{array}{c} \\[-0.2cm] \!\! {\relabelbox \small
\epsfxsize 0.07truein \epsfbox{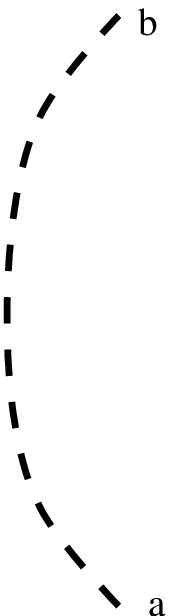}
\adjustrelabel <0cm,-0cm> {a}{\scriptsize $#1$}
\adjustrelabel <0cm,-0cm> {b}{\scriptsize $#2$}
\endrelabelbox} \end{array}}

\newcommand{\strutgraphbot}[2]
{\begin{array}{c} \\[-0.2cm] \!\! {\relabelbox \small
\epsfxsize 0.3truein \epsfbox{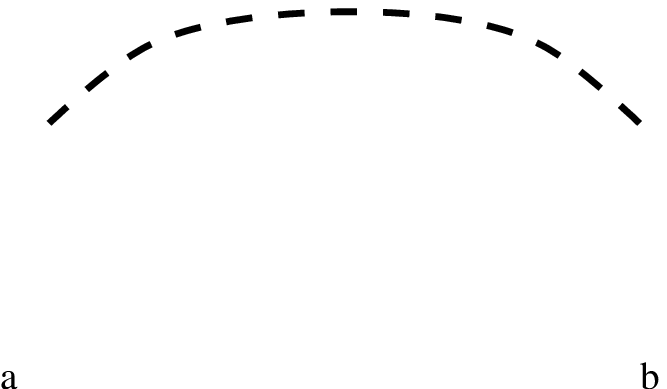}
\adjustrelabel <0cm,-0cm> {a}{\scriptsize $#1$}
\adjustrelabel <0cm,-0cm> {b}{\scriptsize $#2$}
\endrelabelbox} \end{array}}

\newcommand{\Hgraph}[4]
{\begin{array}{c} \\[-0.2cm] \!\! {\relabelbox \small
\epsfxsize 0.2truein \epsfbox{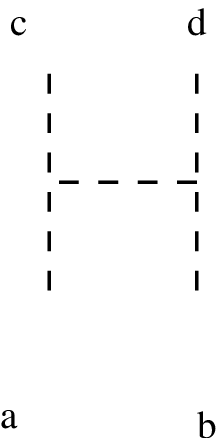}
\adjustrelabel <0cm,-0cm> {a}{\scriptsize $#1$}
\adjustrelabel <0cm,-0cm> {b}{\scriptsize $#2$}
\adjustrelabel <0cm,-0cm> {c}{\scriptsize $#3$}
\adjustrelabel <0cm,-0cm> {d}{\scriptsize $#4$}
\endrelabelbox} \end{array}}

\newcommand{\Hgraphtopright}[4]
{\begin{array}{c} \\[-0.2cm] \!\! {\relabelbox \small
\epsfxsize 0.3truein \epsfbox{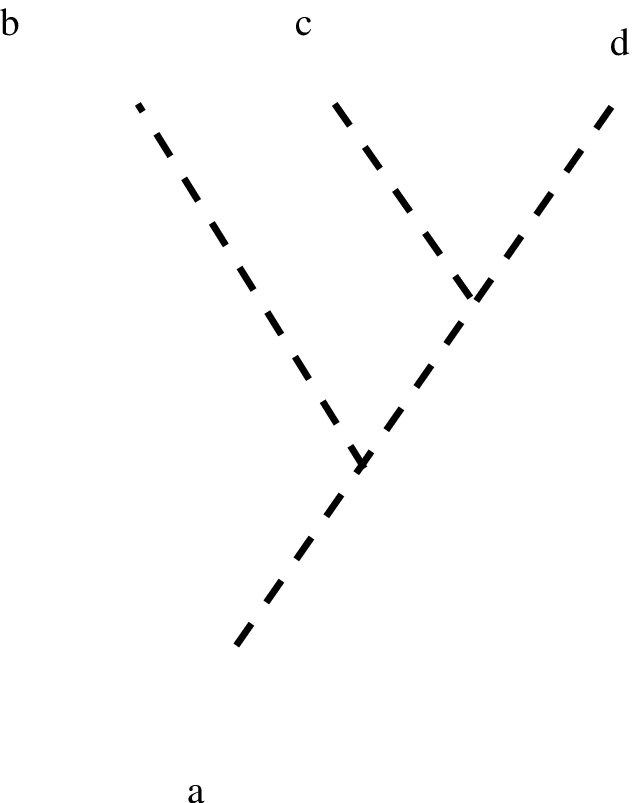}
\adjustrelabel <0cm,-0cm> {a}{\scriptsize $#1$}
\adjustrelabel <0cm,-0cm> {b}{\scriptsize $#2$}
\adjustrelabel <0cm,-0cm> {c}{\scriptsize $#3$}
\adjustrelabel <0cm,-0cm> {d}{\scriptsize $#4$}
\endrelabelbox} \end{array}}

\newcommand{\Hgraphtopleft}[4]
{\begin{array}{c} \\[-0.2cm] \!\! {\relabelbox \small
\epsfxsize 0.3truein \epsfbox{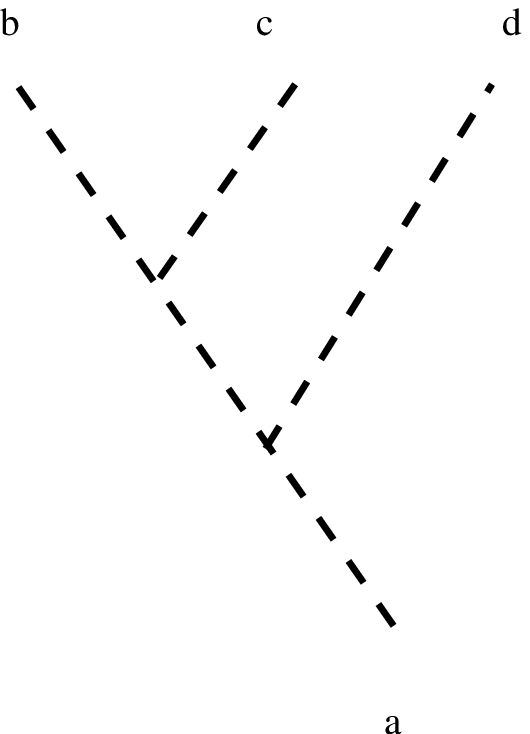}
\adjustrelabel <0cm,-0cm> {a}{\scriptsize $#1$}
\adjustrelabel <0cm,-0cm> {b}{\scriptsize $#2$}
\adjustrelabel <0cm,-0cm> {c}{\scriptsize $#3$}
\adjustrelabel <0cm,-0cm> {d}{\scriptsize $#4$}
\endrelabelbox} \end{array}}

\newcommand{\Hgraphbotright}[4]
{\begin{array}{c} \\[-0.2cm] \!\! {\relabelbox \small
\epsfxsize 0.3truein \epsfbox{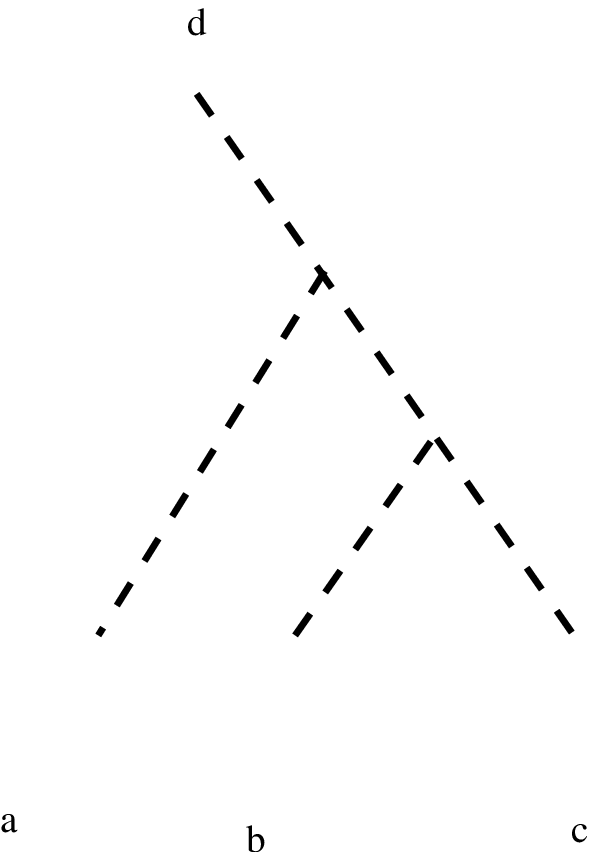}
\adjustrelabel <0cm,-0cm> {a}{\scriptsize $#1$}
\adjustrelabel <0cm,-0cm> {b}{\scriptsize $#2$}
\adjustrelabel <0cm,-0cm> {c}{\scriptsize $#3$}
\adjustrelabel <0cm,-0cm> {d}{\scriptsize $#4$}
\endrelabelbox} \end{array}}

\newcommand{\Hgraphbotleft}[4]
{\begin{array}{c} \\[-0.2cm] \!\! {\relabelbox \small
\epsfxsize 0.3truein \epsfbox{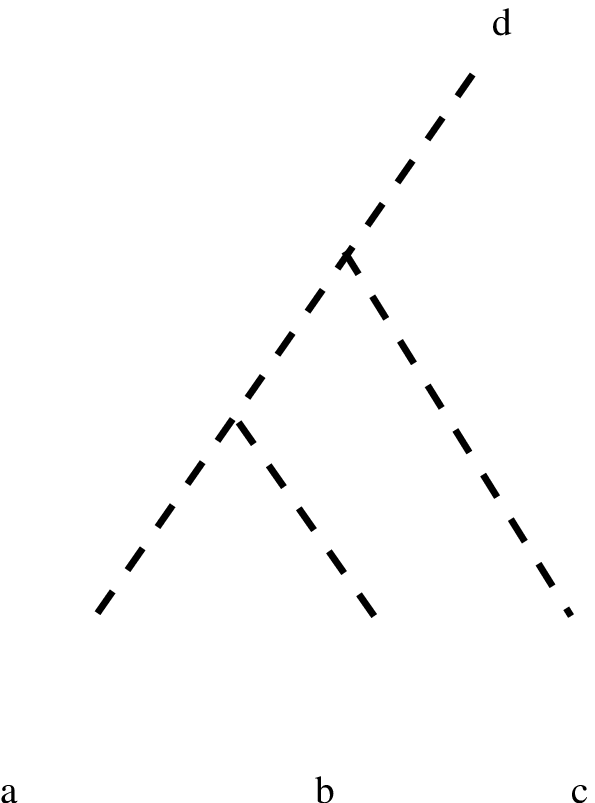}
\adjustrelabel <0cm,-0cm> {a}{\scriptsize $#1$}
\adjustrelabel <0cm,-0cm> {b}{\scriptsize $#2$}
\adjustrelabel <0cm,-0cm> {c}{\scriptsize $#3$}
\adjustrelabel <0cm,-0cm> {d}{\scriptsize $#4$}
\endrelabelbox} \end{array}}

\newcommand{\Hgraphcup}[4]
{\begin{array}{c} \\[-0.2cm] \!\! {\relabelbox \small
\epsfxsize 0.5truein \epsfbox{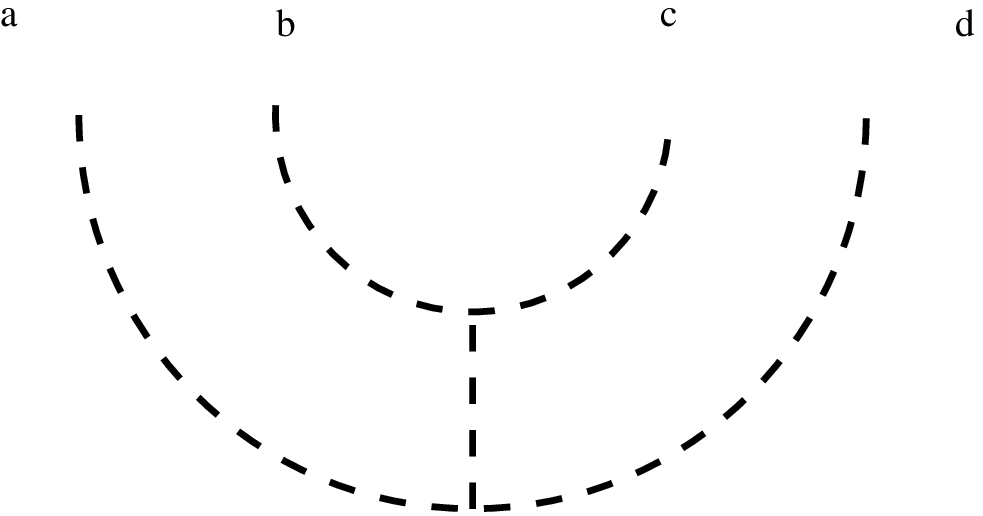}
\adjustrelabel <0cm,-0cm> {a}{\scriptsize $#1$}
\adjustrelabel <0cm,-0cm> {b}{\scriptsize $#2$}
\adjustrelabel <0cm,-0cm> {c}{\scriptsize $#3$}
\adjustrelabel <0cm,-0cm> {d}{\scriptsize $#4$}
\endrelabelbox} \end{array}}

\newcommand{\Hgraphcap}[4]
{\begin{array}{c} \\[-0.2cm] \!\! {\relabelbox \small
\epsfxsize 0.5truein \epsfbox{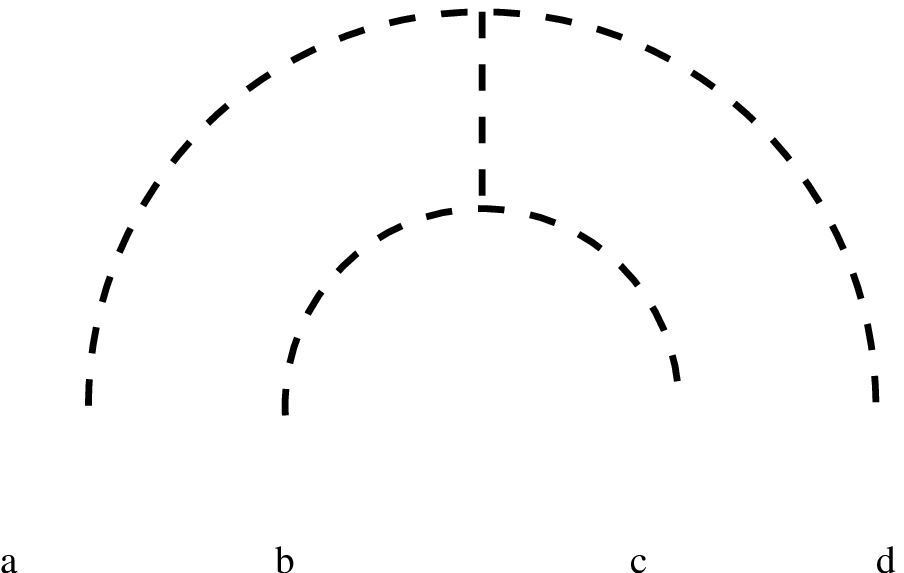}
\adjustrelabel <0cm,-0cm> {a}{\scriptsize $#1$}
\adjustrelabel <0cm,-0cm> {b}{\scriptsize $#2$}
\adjustrelabel <0cm,-0cm> {c}{\scriptsize $#3$}
\adjustrelabel <0cm,-0cm> {d}{\scriptsize $#4$}
\endrelabelbox} \end{array}}

\newcommand{\phigraph}[2]
{\begin{array}{c} \\[-0.2cm] \!\! {\relabelbox \small
\epsfxsize 0.1truein \epsfbox{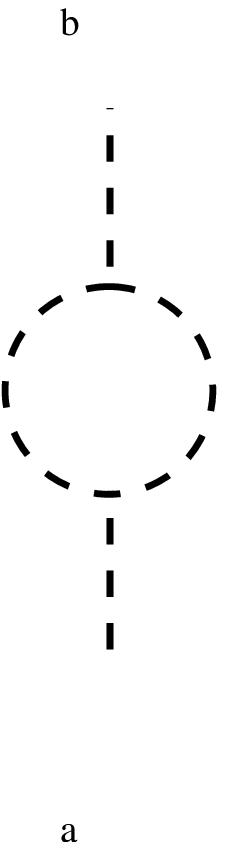}
\adjustrelabel <0cm,-0cm> {a}{\scriptsize $#1$}
\adjustrelabel <0cm,-0cm> {b}{\scriptsize $#2$}
\endrelabelbox} \end{array}}

\newcommand{\phigraphtop}[2]
{\begin{array}{c} \\[-0.2cm] \!\! {\relabelbox \small
\epsfxsize 0.25truein \epsfbox{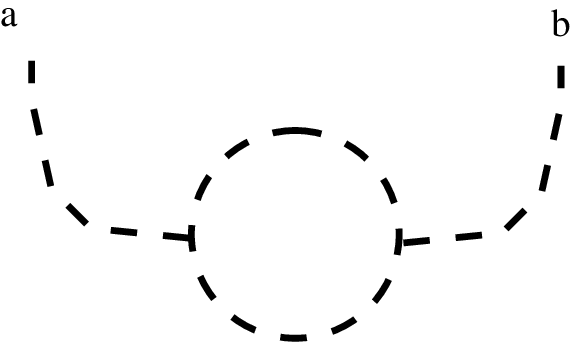}
\adjustrelabel <0cm,-0cm> {a}{\scriptsize $#1$}
\adjustrelabel <0cm,-0cm> {b}{\scriptsize $#2$}
\endrelabelbox} \end{array}}

\newcommand{\phigraphbot}[2]
{\begin{array}{c} \\[-0.2cm] \!\! {\relabelbox \small
\epsfxsize 0.35truein \epsfbox{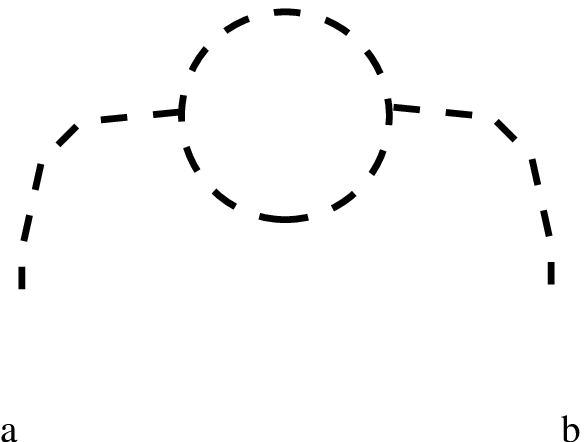}
\adjustrelabel <0cm,-0cm> {a}{\scriptsize $#1$}
\adjustrelabel <0cm,-0cm> {b}{\scriptsize $#2$}
\endrelabelbox} \end{array}}

\newcommand{\Ygraphtop}[3]
{\begin{array}{c} \\[-0.2cm] \!\! {\relabelbox \small
\epsfxsize 0.3truein \epsfbox{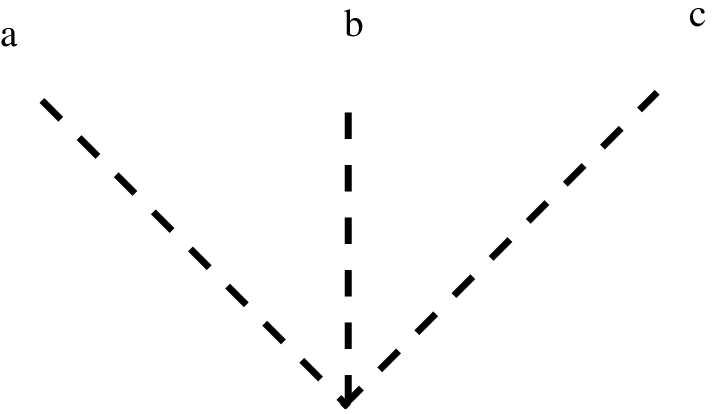}
\adjustrelabel <0cm,-0cm> {a}{\scriptsize $#1$}
\adjustrelabel <0cm,-0cm> {b}{\scriptsize $#2$}
\adjustrelabel <0cm,-0cm> {c}{\scriptsize $#3$}
\endrelabelbox} \end{array}}

\newcommand{\Ygraphbottoptop}[3]
{\begin{array}{c} \\[-0.2cm] \!\! {\relabelbox \small
\epsfxsize 0.2truein \epsfbox{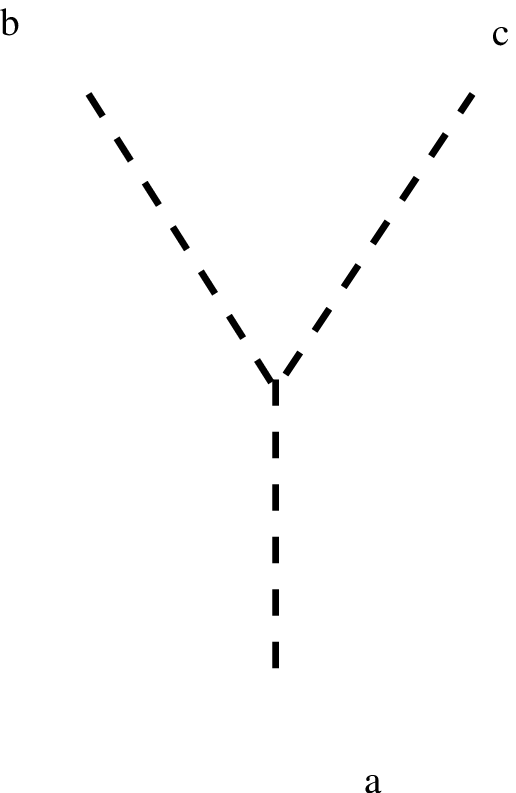}
\adjustrelabel <0cm,-0cm> {a}{\scriptsize $#1$}
\adjustrelabel <0cm,-0cm> {b}{\scriptsize $#2$}
\adjustrelabel <0cm,-0cm> {c}{\scriptsize $#3$}
\endrelabelbox} \end{array}}

\newcommand{\Ygraphbotbottop}[3]
{\begin{array}{c} \\[-0.2cm] \!\! {\relabelbox \small
\epsfxsize 0.3truein \epsfbox{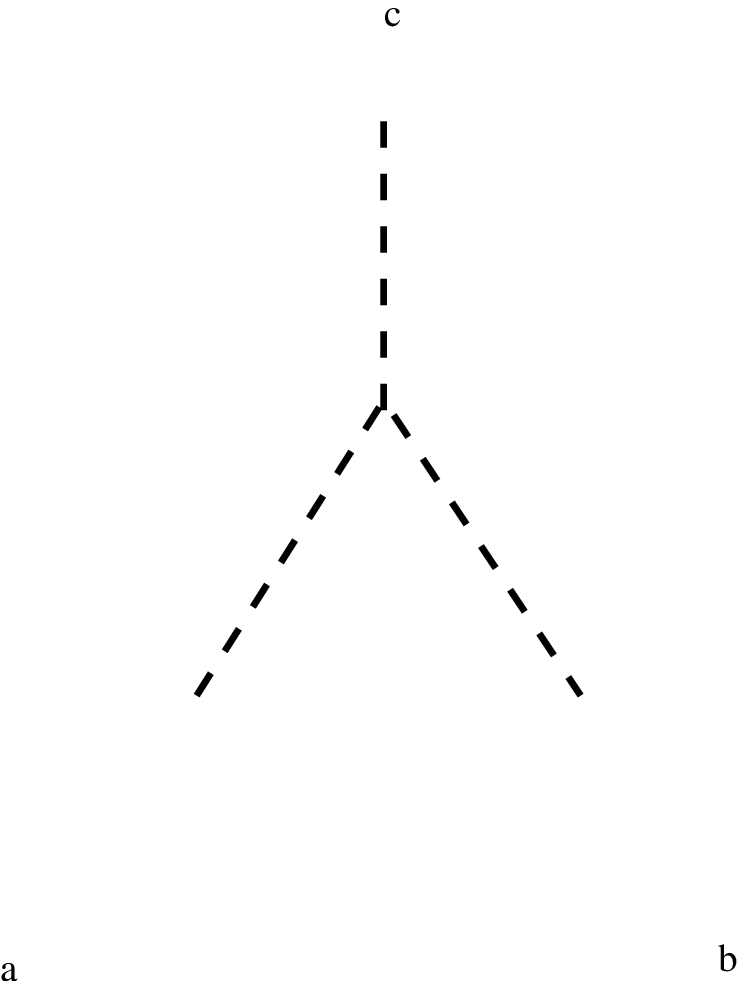}
\adjustrelabel <0cm,-0cm> {a}{\scriptsize $#1$}
\adjustrelabel <0cm,-0cm> {b}{\scriptsize $#2$}
\adjustrelabel <0cm,-0cm> {c}{\scriptsize $#3$}
\endrelabelbox} \end{array}}

\title[]{\large A functorial LMO invariant\\ for Lagrangian cobordisms}

\date{March 22, 2007}

\author[]{Dorin Cheptea}
\address{Center for the Topology and Quantization of Moduli Spaces,
University of Aarhus, Bygning 1530, Ny Munkegade, 8000 Aarhus C, Denmark}
\email{cheptea@imf.au.dk}

\author[]{Kazuo Habiro}
\address{Research Institute for Mathematical Sciences, Kyoto University, Kyoto 606-8502, Japan}
\email{habiro@kurims.kyoto-u.ac.jp}

\author[]{Gw\'ena\"el Massuyeau}
\address{Institut de Recherche Math\'ematique Avanc\'ee, CNRS -- Universit\'e Louis Pasteur,
7 rue Ren\'e Descartes, 67084 Strasbourg, France}
\email{massuyeau@math.u-strasbg.fr}

%\thanks{}

\keywords{3-manifold, finite-type invariant, LMO invariant, Kontsevich integral, cobordism category,
Lagrangian cobordism, homology cylinder, bottom-top tangle, Jacobi diagram, clasper}
\subjclass[2000]{57M27, 57M25}

\begin{abstract}
Lagrangian cobordisms are three-dimensional compact oriented cobordisms 
between once-punctured surfaces, subject to some homological conditions.
We extend the Le--Murakami--Ohtsuki invariant of homology three-spheres to a functor
from the category of Lagrangian cobordisms to a certain category of Jacobi diagrams. 
We prove some properties of this functorial LMO invariant,
including its universality among rational finite-type invariants of Lagrangian cobordisms. 
Finally, we apply the LMO functor to the study of homology cylinders 
from the point of view of their finite-type invariants. 
\end{abstract}

\maketitle

\vspace{0.5cm}

\tableofcontents 

\section{Introduction}

The Kontsevich integral is an invariant of links in $S^3$, the standard $3$-sphere.
In their papers \cite{LM1,LM2}, Le and Murakami extended this invariant 
to a functor from the category of tangles in the standard cube $[-1,1]^3$
to the category of Jacobi diagrams based on $1$-manifolds.
One of the main features of the Kontsevich integral is its 
universality among rational-valued finite-type invariants of tangles (in the Goussarov--Vassiliev sense).

Le, Murakami and Ohtsuki constructed in \cite{LMO} an invariant of closed oriented $3$-manifolds,
which is called the Le--Murakami--Ohtsuki invariant.
The LMO invariant is defined from the Kontsevich integral via surgery presentations of $3$-manifolds in $S^3$. 
For rational homology $3$-spheres, the LMO invariant is universal
among rational-valued finite-type invariants (in the Ohtsuki sense).
Later, Murakami and Ohtsuki  extended in \cite{MO}
the LMO invariant to an invariant of $3$-manifolds with boundary, 
which satisfies modified axioms of TQFT's. 
More recently, Le and the first author constructed from the LMO invariant
a functor from a certain category of $3$-dimensional cobordisms to a certain category of modules \cite{CL1}.
Let us recall that each of those two constructions \cite{MO,CL1} starts with the following two steps: 
\emph{(i)} Extend the Kontsevich integral to framed trivalent
graphs in $S^3$; \emph{(ii)} Unify the extended Kontsevich integral 
and the LMO invariant into a single invariant $Z(M,G)$ of couples $(M,G)$, 
where $M$ is a closed oriented $3$-manifold and $G\subset M$ is an embedded framed trivalent graph.
Then, a compact oriented $3$-manifold with boundary is obtained from each couple $(M,G)$ 
by cutting a regular neighborhood $\hbox{N}(G)$ of $G$ in $M$. 
If the connected components of $G$ were split into two parts,
say the ``top'' part $G^+$ and the ``bottom'' part $G^-$, 
then $M \setminus \hbox{N}(G)$ can be regarded as a cobordism between \emph{closed} surfaces, 
namely from $\partial \hbox{N}(G^+)$ to $- \partial \hbox{N}(G^-)$. 
Finally, the LMO invariant of the cobordism $M \setminus \hbox{N}(G)$ 
is defined in \cite{MO,CL1} to be the Kontsevich--LMO invariant $Z(M,G)$ of the couple $(M,G)$.\\

In this paper, we propose an alternative solution to the problem of extending the LMO invariant
to $3$-manifolds with boundary. In contrast with the previous two constructions \cite{MO,CL1}, 
we prefer to work with cobordisms between \emph{once-punctured} surfaces.
This technical choice has two advantages: On the one hand, it avoids us to extend the Kontsevich integral 
to trivalent graphs in $S^3$; On the other hand, 
it allows us to work with monoidal categories, and to construct tensor-preserving functors. 

Moreover, we \emph{normalize} the Kontsevich--LMO invariant $Z$ to obtain
an invariant $\Ztilde$ of $3$-manifolds with boundary. 
This normalization is done in such a way that the gluing formula satisfied by $\Ztilde$
can be described by a simple combinatorial formula.

Thus, our main result is the extension of the LMO invariant to a functor $\Ztilde$, 
which is defined on a certain category of cobordisms between once-punctured surfaces,
and is valued in a certain category of Jacobi diagrams with a facile composition law.
In order to present this LMO functor in more details, 
we need to specify first the kind of cobordisms to which it applies:

\subsection{Lagrangian cobordisms}

Let $\Cob$ denote the category of cobordisms between once-punctured surfaces,
as introduced by Crane and Yetter \cite{CY} and independently by Kerler \cite{Kerler}.
The objects of $\Cob$ are nonnegative integers $g$, 
to each of which is assigned a compact oriented connected surface $F_g$ 
of genus $g$ with one boundary component.
The morphisms from $g_+$ to $g_-$ are the homeomorphism classes (relative to boundary
parameterizations) of cobordisms between the surfaces $F_{g_+}$ and $F_{g_-}$. 
Observe that such cobordisms can be glued ``side-by-side'',
which gives $\Cob$ a monoidal structure.

The subcategory $\LCob$ of $\Cob$ will consist of ``Lagrangian cobordisms.'' 
Here, let us give a rough description of this notion. 
Let $F_+$ and $F_-$ be two compact connected oriented surfaces with one boundary component.
Let $A_+$ and $A_-$ be Lagrangian subgroups of $H_1(F_+;\Z )$ and
$H_1(F_-;\Z )$, respectively.  A {\em Lagrangian cobordism} between
$(F_+,A_+)$ and $(F_-,A_-)$ is a cobordism $M$ between $F_+$ and $F_-$ which
satisfies
\begin{enumerate}
\item $H_1(M;\Z ) = m_{-,*}(A_-) + m_{+,*}(H_1(F_+;\Z ))$,
\item $m_{+,*}(A_+) \subset m_{-,*}(A_-)$ in $H_1(M;\Z )$.
\end{enumerate}
Here $m_{\pm}\colon F_\pm \hookrightarrow M$ is the
inclusion, and $m_{\pm,*}$ is the induced map on $H_1 (\ \centerdot\ ;\Z )$.
Using the Mayer--Vietoris theorem, one easily checks that the
composite of two Lagrangian cobordisms is again Lagrangian.  Thus
there is a category $\mathbf{LCob}$ whose objects are pairs $(F,A)$ of
a punctured surface and a Lagrangian subgroup $A\subset H_1(F;\Z )$, 
and whose morphisms are homeomorphism classes 
(relative to boundary parameterizations) of Lagrangian cobordisms.

The subcategory $\LCob$ of $\Cob$, which will be defined later with more care, 
is essentially a skeleton of $\mathbf{LCob}$. We choose a  Lagrangian subgroup $A_g$ 
for the standard surface $F_g$; The objects in $\LCob$ are then 
nonnegative integers $g$ and the morphisms from $g_+$ to $g_-$ in $\LCob$ are morphisms from
$(F_{g_+},A_{g_+})$ to $(F_{g_-},A_{g_-})$ in $\mathbf{LCob}$. 
The category $\LCob$ contains as a subcategory a ``punctured version'' of the
category of ``semi-Lagrangian cobordisms'' used by Le and the first author in \cite{CL1}.
Actually, $\LCob$ can be identified with the category of ``bottom tangles in homology handlebodies'' 
as defined by the second author in \cite{Habiro_BT}.

Lagrangian cobordisms may be considered 
as a {\em natural generalization} of integral homology cubes, since
the latter are the morphisms from $0$ to $0$ in the category $\LCob$.  
Another reason to be interested in the class of Lagrangian cobordisms is that it contains 
homology cylinders, which have been introduced by Goussarov and the second author in \cite{Goussarov_note,Habiro_claspers}
and play an important role in the study of finite-type invariants.
In particular, let us observe that the Torelli group of the surface $F_g$ 
embeds into the monoid $\LCob(g,g)$ via the mapping cylinder construction.

To define the LMO functor on  Lagrangian cobordisms, it is
convenient to enhance the category $\LCob$ to a category $\LqCob$
of \emph{Lagrangian $q$-cobordisms}. The sets of morphisms of $\LqCob$ are the same as in $\LCob$,
but the objects are now parenthesized words $w$ in a single letter $\bullet$. 
Thus, there is a natural functor $\LqCob \rightarrow \LCob$ which simply forgets the parenthesization
(e.g. $((\bullet\bullet)\bullet)\mapsto 3$) and, as a category, $\LqCob$ is equivalent to $\LCob$.

\subsection{The category of top-substantial Jacobi diagrams}

Let us now roughly describe the category $\tsA$ in which our LMO functor
takes values, assuming that the reader has a certain familiarity with Jacobi diagrams.
The objects of $\tsA$ are nonnegative integers.  The
set of morphisms $\tsA(g,f)$ from $g$ to $f$ in $\tsA$ is the
$\Q $-vector space of ``top-substantial'' Jacobi diagrams with univalent vertices labeled by the set
$$
\{1^+,\ldots ,g^+\}\cup\{1^-,\ldots ,f^-\}.
$$
Here, \emph{top-substantiality} means
that no component of the graph is a strut whose two univalent
vertices are colored by elements of $\{1^+,\ldots ,g^+\}$.
For example, here is a Jacobi diagram defining a morphism from $4$ to $5$ in $\tsA$:\\[0.3cm]

\centerline{\relabelbox \small
\epsfxsize 2truein \epsfbox{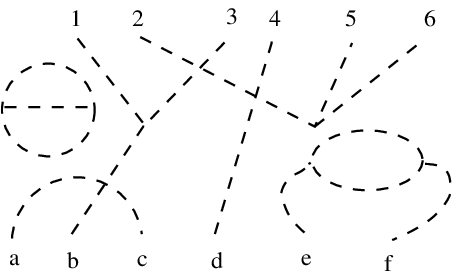}
\adjustrelabel <0cm,-0cm> {1}{$1^+$}
\adjustrelabel <0cm,-0cm> {2}{$1^+$}
\adjustrelabel <0cm,-0cm> {3}{$2^+$}
\adjustrelabel <0cm,-0cm> {4}{$2^+$}
\adjustrelabel <0cm,-0cm> {5}{$3^+$}
\adjustrelabel <0cm,-0cm> {6}{$4^+$}
\adjustrelabel <0cm,-0cm> {a}{$1^-$}
\adjustrelabel <0cm,-0cm> {b}{$1^-$}
\adjustrelabel <0cm,-0cm> {c}{$2^-$}
\adjustrelabel <0cm,-0cm> {d}{$3^-$}
\adjustrelabel <0cm,-0cm> {e}{$5^-$}
\adjustrelabel <0cm,-0cm> {f}{$5^-$}
\endrelabelbox}
\vspace{0.4cm}

\noindent As usual, the space $\tsA$ is completed with respect to the degree of diagrams, 
so that we consider formal series of Jacobi diagrams.  
The composition map
\begin{gather*}
  \circ : \tsA(g,f)\times \tsA(h,g) \longrightarrow  \tsA(h,f)
\end{gather*}
is simply defined as follows: Given $x\in \tsA(g,f)$ and $y\in \tsA(h,g)$,
$x\circ y\in \tsA(h,f)$ is obtained from $x \sqcup y$ by ``contracting''
the $i^+$-colored vertices in $x$ and the $i^-$-colored vertices in
$y$ for all $i=1,\ldots ,g$.
The identity morphism of the object $g$ in $\tsA$ is then given by
$$
\Id_g = \exp_\sqcup \left( \sum_{i=1}^g \strutgraph{i^-}{i^+} \right).
$$
There is also a natural monoidal structure on $\tsA$.

\subsection{The LMO functor}

Thus, our main construction is a tensor-preserving functor
$$
\Ztilde : \LqCob \longrightarrow \tsA.
$$
At the level of objects, $\Ztilde$ just sends a non-associative word $w$ to its length $|w|$.
At the level of morphisms, the series of Jacobi diagrams $\Ztilde(M)$
assigned to a Lagrangian $q$-cobordism $M\in \LqCob(w,v)$ is defined as follows.
First, we present the cobordism $M$ by a couple $(B,\gamma)$ of
an integral homology cube $B$ and a framed tangle $\gamma \subset B$  of a certain type, which we call a ``bottom-top tangle.''
This is inspired from the way cobordisms between closed surfaces are presented in \cite{CL1}:
Since our surfaces have one boundary component, we work with tangles in homology cubes 
rather than with trivalent graphs in homology spheres.
Next, we normalize the Kontsevich--LMO invariant $Z(B,\gamma)$ to
an invariant $\Ztilde(M)\in \tsA(|w|,|v|)$ in such a way that $\Ztilde$ is functorial.
To define the Kontsevich--LMO invariant $Z$ and to carry out this normalization,
we use the Aarhus integral developed by Bar-Natan, Garoufalidis, Rozansky and D. Thurston \cite{BGRT1,BGRT2}.

By construction, $\Ztilde$ sends a homology cube $B\in\LCob(0,0)$ 
to the LMO invariant of the homology sphere $\hat{B}$ obtained by ``recapping'' $B$.  
Thus, $\Ztilde$ should be considered as a functorial extension of the LMO invariant.
In particular, the reduction of $\Ztilde$ to Jacobi diagrams with no more than two trivalent vertices
defines a functorial extension of the Casson invariant. 

As announced in \cite{Habiro_claspers}, $\LCob$ is finitely generated as a monoidal category, 
so that any Lagrangian $q$-cobordism can be decomposed 
(with respect to the composition law of cobordisms and their tensor product) into ``building blocks''
of a finite number of types. Therefore, the functor $\Ztilde$ is determined 
by its values on those ``building blocks.'' 
As an illustration, we have computed those values at the Casson invariant level, 
i.e$.$ modulo Jacobi diagrams with more than two trivalent vertices.

Furthermore, after a suitable reduction, our functor $\Ztilde$ factors 
through the category of Lagrangian $q$-cobordisms between \emph{closed} surfaces.
We show that the TQFT constructed in \cite{CL1} can be recovered from this reduction of $\Ztilde$.

\subsection{Properties and applications of the LMO functor}

Just as is the case for the Kontsevich--LMO invariant, the functor $\Ztilde$ takes group-like values.
More precisely, the series of Jacobi diagram $\Ztilde(M)$ assigned to each cobordism $M\in \LqCob(w,v)$
splits into two group-like elements:
The ``$s$-part'' $\Ztilde^s(M)$ which only contains struts, 
and the ``$Y$-part'' $\ZtildeY(M)$ which does not contain strut at all. 
Whereas the former only contains homological information about $M$,
the latter is very rich in the sense that it contains all rational-valued finite-type invariants.
This universality among finite-type invariants is deduced from the functoriality of $\Ztilde$ 
using the ``clasper calculus'' of \cite{Goussarov_clovers, Habiro_claspers, GGP}.

A quite illustrative application of our results is offered by homology cylinders. 
In this case, the ``$Y$-part'' of the LMO functor restricts to a homomorphism 
$$
\ZtildeY : \Cyl(F_g) \longrightarrow \AY\left(\{1^+,\ldots ,g^+\}\cup\{1^-,\ldots ,g^-\}\right)
$$
from the monoid of homology cylinders over $F_g$ to the space of Jacobi diagrams with no strut, 
which is equipped with a certain multiplication $\star$. 
In contrast with the LMO-type invariant of homology cylinders 
introduced by Habegger in \cite{Habegger}, our universal invariant $\ZtildeY$ is \emph{multiplicative}.
This property allows us to compute diagrammatically the algebra dual 
to rational finite-type invariants of homology cylinders, 
as well as the ``Lie algebra of homology cylinders'' introduced by the second author in \cite{Habiro_claspers}.
Moreover, by adapting Habegger's method \cite{Habegger}, we explain 
how the first non-vanishing Johnson homomorphism of a  homology cylinder $M$ 
can be extracted from its LMO invariant $\ZtildeY(M)$.
We deduce that the LMO  homomorphism $\ZtildeY$ is injective on the Torelli group of $F_g$.\\[0.5cm]

Before going into the details of the constructions and proofs, we would like to fix a few conventions.
In this paper, we agree that
\begin{itemize}
\item unless otherwise specified, all homology groups are computed with integer coefficients;
\item one-dimensional objects that are drawn on diagrams, such as graphs or links,
are given the ``blackboard framing'', i.e$.$ are thickened along the plan;
\item  a \emph{tensor-preserving} functor $F:\mathcal{C} \to \mathcal{C'}$,
between two monoidal categories $\mathcal{C}$ and $\mathcal{C'}$, 
is a functor which strictly respects the tensor products
at the level of objects and morphisms, and which strictly preserves the unit objects. 
(However, $F$ is not required to preserve the associativity and unitality constraints.)
\end{itemize}

\vspace{0.5cm}

\section{Cobordisms and tangles}

We start by introducing the various categories of tangles and cobordisms 
that are used throughout the paper. Our goal is to define Lagrangian cobordisms,
and to explain how they can be presented as tangles in homology cubes of a certain type.

\subsection{The category $\Cob$ of cobordisms}

\label{subsec:cob}

First of all, we recall the category $\Cob$ of cobordisms between surfaces with one boundary component.
This category has been introduced by Crane and Yetter \cite{CY} as well as Kerler \cite{Kerler}.

\begin{figure}[h]
\centerline{\relabelbox \small
\epsfxsize 5truein \epsfbox{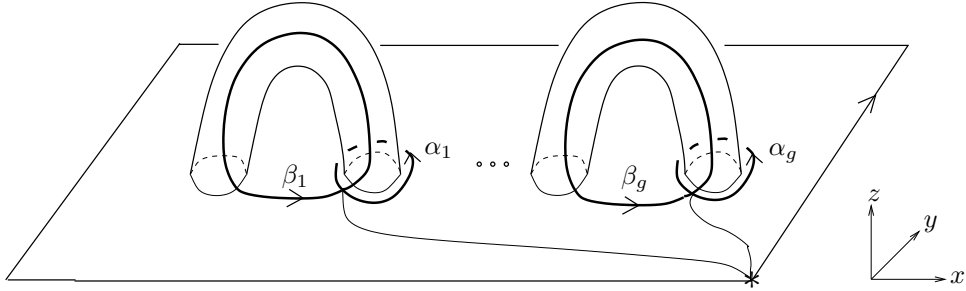}
\adjustrelabel <0cm,-0cm> {x}{$x$}
\adjustrelabel <0cm,-0cm> {y}{$y$}
\adjustrelabel <0cm,-0cm> {z}{$z$}
\adjustrelabel <0cm,-0cm> {alpha1}{$\alpha_1$}
\adjustrelabel <0cm,-0cm> {alphag}{$\alpha_g$}
\adjustrelabel <0cm,-0cm> {beta1}{$\beta_1$}
\adjustrelabel <0cm,-0cm> {betag}{$\beta_g$}
\adjustrelabel <-0.05cm,-0cm> {*}{{\Large $*$}}
\endrelabelbox}
\caption{The standard surface $F_g$ and its system of meridians and parallels $(\alpha,\beta)$.}
\label{fig:surface}
\end{figure}

For each integer $g\geq 0$, let $F_g$ be a compact connected oriented surface of genus $g$ with
one boundary component, which is \emph{fixed} once and for all. 
We think of it as embedded in the ambient space $\mathbb{R}^3$ (with coordinates $x$, $y$, $z$) 
and obtained from the square $[-1,1] \times [-1,1] \times 0$ by adding
$g$ handles uniformly in the $x$ direction. See Figure \ref{fig:surface} where 
the orientation on $F_g$ is materialized through the orientation it induces on $\partial F_g$.
We also \emph{fix} a base point $*$ and a basis for
$\pi_1\left(F_{g},*\right)$ by choosing
a system of meridians and parallels
$(\alpha_1,\beta_1,\dots, \alpha_g,\beta_g)$ as shown on the same picture.

\begin{remark}
\label{rem:homeomorphism}
In order to identify (up to isotopy) a  surface $S$ of genus $g$ with $F_g$,
it is enough to specify which are the images of $*$,
$\alpha_1,\beta_1,\dots, \alpha_g,\beta_g$ on $S$
and the induced orientation on $\partial S$.
\end{remark}

\begin{figure}[h]
\centerline{\relabelbox \small
\epsfxsize 3truein \epsfbox{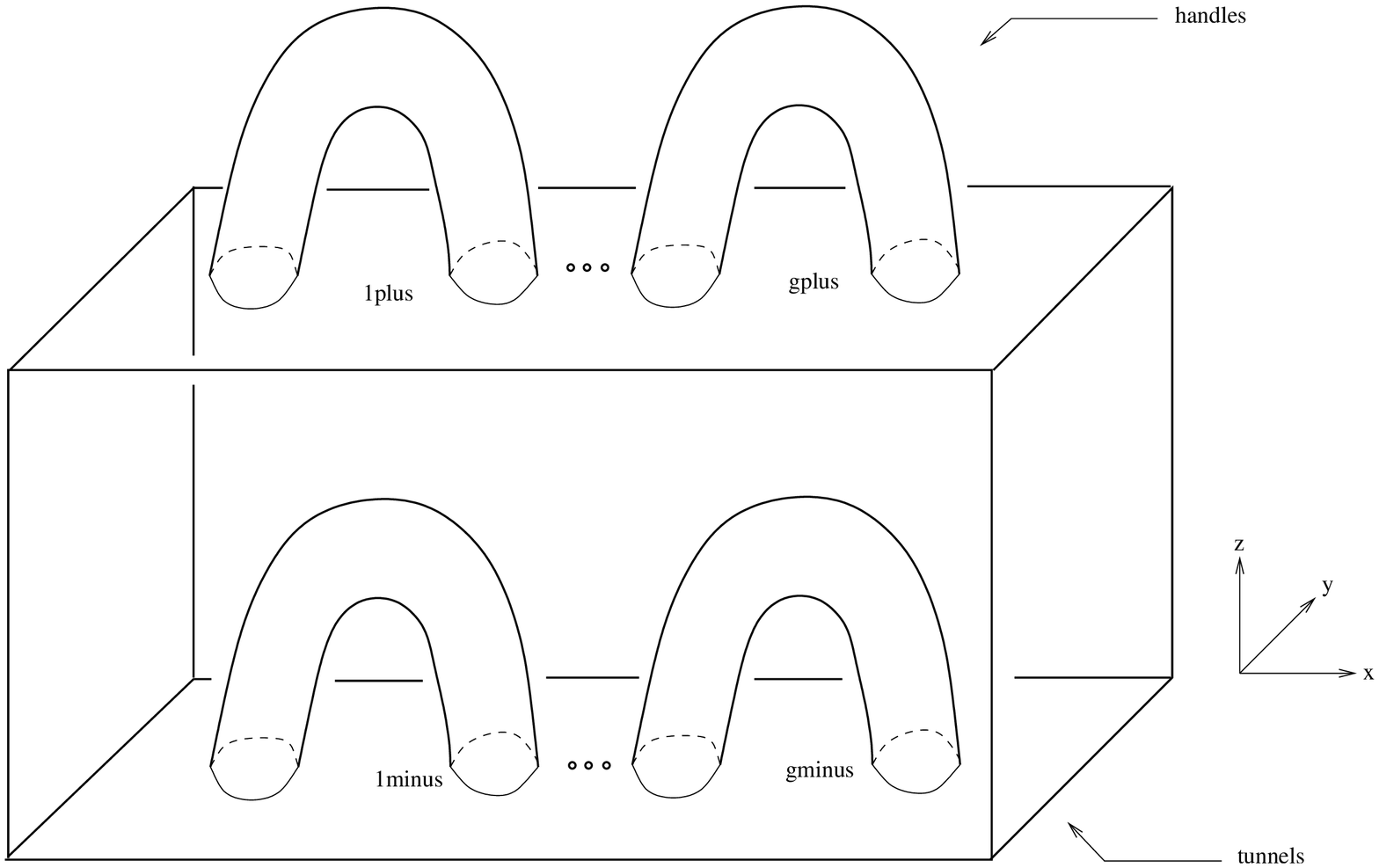}
\adjustrelabel <0cm,-0cm> {x}{$x$}
\adjustrelabel <0cm,-0cm> {y}{$y$}
\adjustrelabel <0cm,-0cm> {z}{$z$}
\adjustrelabel <0cm,-0.05cm> {1minus}{$1$}
\adjustrelabel <0cm,-0.05cm> {1plus}{$1$}
\adjustrelabel <0cm,-0cm> {gminus}{$g_-$}
\adjustrelabel <0cm,-0cm> {gplus}{$g_+$}
\adjustrelabel <0cm,-0cm> {handles}{handles}
\adjustrelabel <0cm,-0cm> {tunnels}{tunnels}
\endrelabelbox}
\caption{The cube $C_{g_-}^{g_+}$ with $g_-$ tunnels and $g_+$ handles.}
\label{fig:cube}
\end{figure}

We denote by $C_{g_-}^{g_+}$ the \emph{cube with $g_-$ tunnels and $g_+$ handles}:
This is the compact oriented $3$-manifold obtained
from the cube $[-1,1]^3$ by adding $g_\pm$ $1$-handles along
$[-1,1]\times [-1,1] \times (\pm 1)$,
uniformly in the $x$ direction, as shown in Figure \ref{fig:cube}.
We note the two canonical embeddings
\begin{equation}
\label{eq:bottom_top_surfaces}
F_{g_-} \hookrightarrow -\partial C_{g_-}^{g_+}
\quad \quad \textrm{and} \quad \quad
F_{g_+} \hookrightarrow \partial C_{g_-}^{g_+}
\end{equation}
obtained by appropriate translations in the $z$ direction.

\begin{definition}
Let $g_-\geq 0$ and $g_+ \geq 0$ be integers.
A \emph{cobordism} from $F_{g_+}$ to $F_{g_-}$ is an equivalence class of couples $(M,m)$ where
\begin{itemize}
\item $M$ is a compact connected oriented $3$-manifold,
\item $m: \partial C_{g_-}^{g_+} \to M$ is an orientation-preserving homeomorphism onto $\partial M$,
\end{itemize}
two such couples $(M,m)$ and $(M',m')$ being considered as \emph{equivalent} if there exists 
an orientation-preserving homeomorphism $f:M \to M'$ such that $f \circ m=m'$.

By the inclusions (\ref{eq:bottom_top_surfaces}), $m$ restricts to two embeddings
$$
m_-:F_{g_-} \hookrightarrow M
\quad \textrm{and} \quad
m_+:F_{g_+} \hookrightarrow M
$$
whose images are called \emph{bottom surface} and \emph{top surface} of the cobordism $M$ respectively.
\end{definition}

\begin{figure}[h]
\centerline{\relabelbox \small
\epsfxsize 4truein \epsfbox{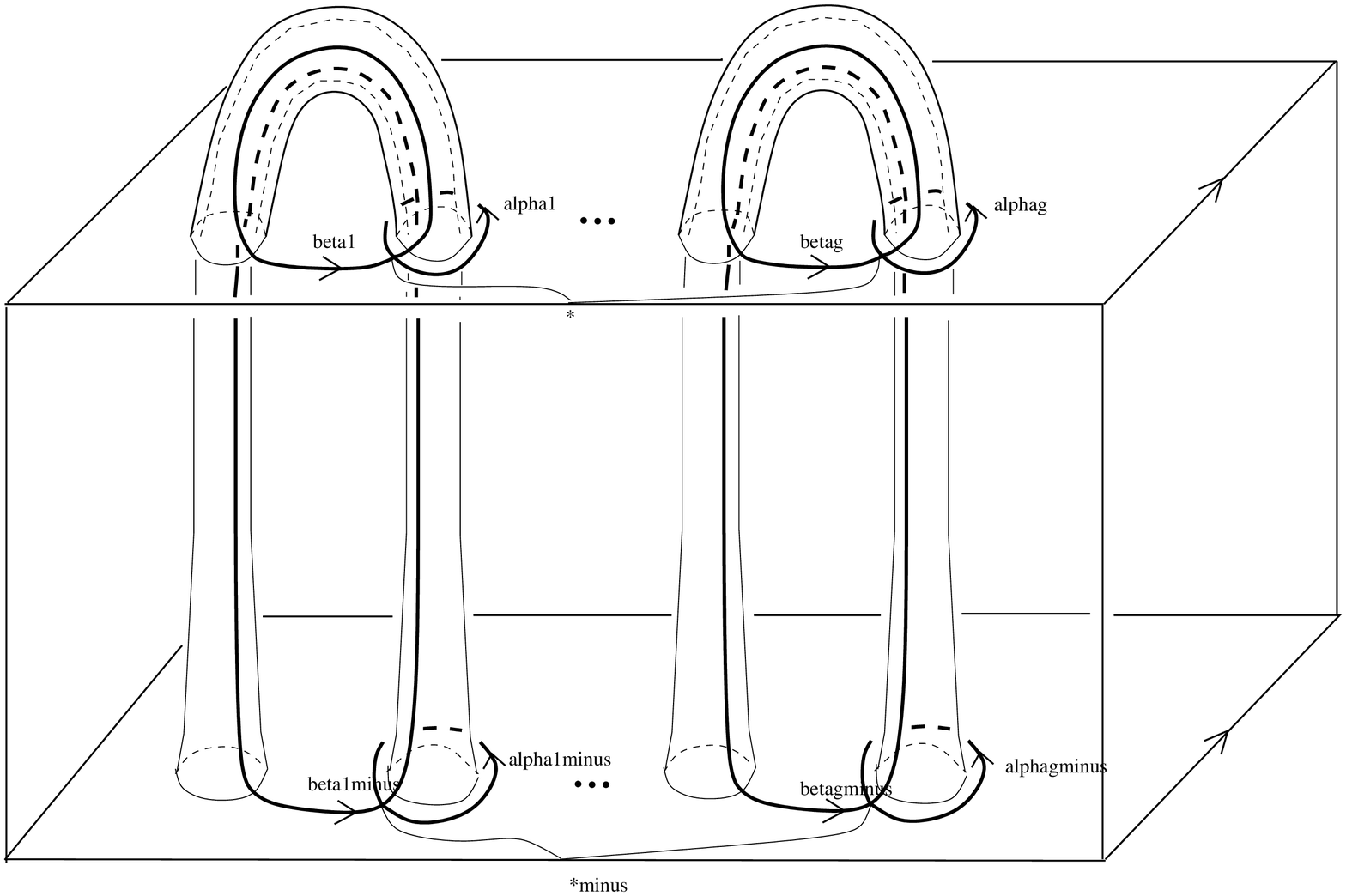}
\adjustrelabel <0cm,-0cm> {alpha1}{$\alpha_1$}
\adjustrelabel <0cm,-0cm> {beta1}{$\beta_1$}
\adjustrelabel <0cm,-0cm> {alphag}{$\alpha_g$}
\adjustrelabel <0cm,-0cm> {betag}{$\beta_g$}
\adjustrelabel <0cm,-0cm> {alpha1minus}{$\alpha_1$}
\adjustrelabel <0cm,-0cm> {beta1minus}{$\beta_1$}
\adjustrelabel <0cm,-0cm> {alphagminus}{$\alpha_g$}
\adjustrelabel <0cm,0.1cm> {betagminus}{$\beta_g$}
\adjustrelabel <0cm,0.05cm> {*}{$*$}
\adjustrelabel <-0.1cm,0.1cm> {*minus}{$*$}
\endrelabelbox}
\caption{The cobordism $\Id_g$ from $F_{g}$ to $F_{g}$.}
\label{fig:cob_id}
\end{figure}

Given cobordisms $(M,m)$ from $F_{g_+}$ to $F_{g_-}$ 
and $(N,n)$ from $F_{h_+}$ to $F_{h_-}$  such that $g_+=h_-$,
one obtains a new cobordism $(M,m) \circ (N,n)$ from $F_{h_+}$ to $F_{g_-}$ 
by ``stacking'' $N$ on the top of $M$ and parametrizing the boundary
of the new manifold in the obvious way. 
Thus, one obtains a category  $\Cob$ 
whose objects are non-negative integers $g$
and whose sets of morphisms $\Cob(g_+,g_-)$ are cobordisms from $F_{g_+}$ to $F_{g_-}$.
The identity of $g\geq 0$ in the category $\Cob$ is $F_g \times [-1,1]$ 
with the obvious parameterization shown on Figure \ref{fig:cob_id}.

The category $\Cob$ is monoidal (in the strict sense), with tensor product $\otimes$ given
by horizontal juxtaposition of cobordisms in the $x$ direction. 
So, to sum up, we have two operations on cobordisms:
$$
M \circ N  := \begin{array}{|c|}\hline N \\ \hline M\\ \hline \end{array}
\quad \quad \hbox{and} \quad \quad
M \otimes N  :=  \begin{array}{|c|c|}  \hline M & N \\ \hline \end{array}\ .
$$

\begin{example}
\label{ex:mapping_cylinder}
Let $\M(F_g)$ denote the \emph{mapping class group} of the surface $F_g$, i.e$.$
the group of isotopy classes of homeomorphisms $F_g \to F_g$ that fix $\partial F_g$ pointwise. 
The \emph{mapping cylinder} construction
$$
\M(F_g) \longrightarrow \Cob(g,g), \ h \longmapsto \left(F_g\times [-1,1], (\Id \times (-1)) \cup (h \times 1)\right)
$$
defines a monoid homomorphism, which is injective.
\end{example}

\subsection{The category $\LCob$ of Lagrangian cobordisms}

\label{subsec:LCob}

We now introduce the subcategory of $\Cob$ in which we are interested.
For this, we distinguish the following Lagrangian subgroup of $H_1(F_g)$:
$$
A_g:= \Ker\left( \incl_*: H_1(F_g) \to H_1(C_0^g)\right) = \langle \alpha_1, \dots, \alpha_g \rangle.
$$

\begin{definition}
\label{def:LCob}
A cobordism $(M,m)$ from $F_{g_+}$ to $F_{g_-}$ is  \emph{Lagrangian-preserving} 
(or, for short, \emph{Lagrangian}) if the following two conditions are satisfied:
\begin{enumerate}
\item $H_1(M) = m_{-,*}(A_{g_-}) + m_{+,*}(H_1(F_{g_+}))$,
\item $m_{+,*}(A_{g_+}) \subset m_{-,*}(A_{g_-})$ as subgroups of $H_1(M)$.
\end{enumerate}
\end{definition}

\noindent
If we consider the following supplement to $A_g$: 
$$
B_g := \Ker\left( \incl_*: H_1(F_g) \to H_1(C_g^0)\right) = \langle \beta_1, \dots, \beta_g \rangle,
$$
then, it is easily seen that condition (1) can be replaced in presence of (2) by
\begin{itemize}
\item[(1')] $m_{-,*} \oplus m_{+,*}: A_{g_-} \oplus B_{g_+} \to H_1(M)$ is an isomorphism.
\end{itemize}
With a Mayer--Vietoris argument, one checks that the composition of two Lagrangian cobordisms 
is Lagrangian as well. 
We denote by $\LCob$ the monoidal subcategory of $\Cob$ consisting of Lagrangian cobordisms.

\begin{example}
The mapping cylinder of an $h\in \M(F_g)$ is a Lagrangian cobordism if, and only if, 
$h_*:H_1(F_g) \to H_1(F_g)$ sends $A_g$ to itself.
\end{example}

Among Lagrangian cobordisms, some have a more specific property:

\begin{definition}
\label{def:sLCob}
A cobordism $(M,m)$ from $F_{g_+}$ to $F_{g_-}$ is  \emph{special Lagrangian} 
if it satisfies $C_0^{g_-} \circ M = C_0^{g_+}$.
\end{definition}

\noindent
A Mayer--Vietoris argument shows that ``special Lagrangian'' implies ``Lagrangian.''
The composition of two special Lagrangian cobordisms is special Lagrangian as well.
We denote by $\sLCob$ the monoidal subcategory of $\LCob$ 
consisting of special Lagrangian cobordisms.

\begin{remark}
The category $\LCob$ is isomorphic to the category of
``bottom tangles in homology handlebodies'' introduced by the second author
in \cite[\S 14.5]{Habiro_BT}, while the subcategory $\sLCob$ corresponds
to the category of ``bottom tangles in handlebodies''.
\end{remark}

\begin{remark}
If one takes homology with coefficients in $\Q$ instead of $\Z$, 
one defines in the same way the category $\QLCob$
of \emph{rational Lagrangian cobordisms}.
\end{remark}

\subsection{The category $\btT$ of bottom-top tangles}

\label{subsec:btT}

We now describe a way of presenting cobordisms between surfaces with one boundary component.

\begin{figure}[h]
\centerline{\relabelbox \small
\epsfxsize 3truein \epsfbox{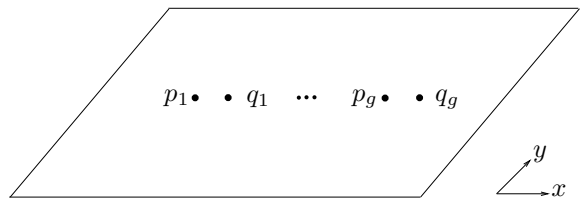}
\adjustrelabel <-0.05cm,-0cm> {p1}{$p_1$}
\adjustrelabel <-0.1cm,-0cm> {pg}{$p_g$}
\adjustrelabel <0cm,-0cm> {q1}{$q_1$}
\adjustrelabel <0cm,-0cm> {qg}{$q_g$}
\adjustrelabel <0cm,-0cm> {x}{$x$}
\adjustrelabel <0cm,-0cm> {y}{$y$}
\endrelabelbox}
\caption{The standard pairs of points
$(p_1,q_1), \dots,$ $(p_g,q_g)$ on $[-1,1]^2$.}
\label{fig:points}
\end{figure}

For all integer $g\geq 1$, denote by $(p_1,q_1), \dots,$ $(p_g,q_g)$
the $g$ pairs of points on $[-1,1]^2$ taken uniformly in the $x$ direction
as shown in Figure \ref{fig:points}.

\begin{definition}
\label{def:bt_tangle}
A \emph{bottom-top tangle} of \emph{type} $(g_+,g_-)$ is an equivalence class of couples $(B,\gamma)$ where $B=(B,b)$ 
is a cobordism from $F_0$ to $F_0$  and $\gamma=(\gamma^+,\gamma^-)$ is a framed oriented tangle with
$g_-$ \emph{bottom} components $\gamma_1^-,\dots,\gamma_{g_-}^-$
and $g_+$ \emph{top} components $\gamma_1^+,\dots,\gamma_{g_+}^+$ such that
\begin{itemize}
\item each $\gamma_j^-$ runs from $q_j\times (-1) $ to $p_j\times (-1)$,
\item each $\gamma_j^+$ runs from $p_j\times 1$ to $q_j\times 1$,
\end{itemize}
two such couples $(B,\gamma)$ and $(B',\gamma')$ being considered as \emph{equivalent} if there exists
an equivalence $(B,b) \to (B',b')$ sending $\gamma$ to $\gamma'$.
\end{definition}

\begin{figure}[h]
\centerline{\relabelbox \small
\epsfxsize 3truein \epsfbox{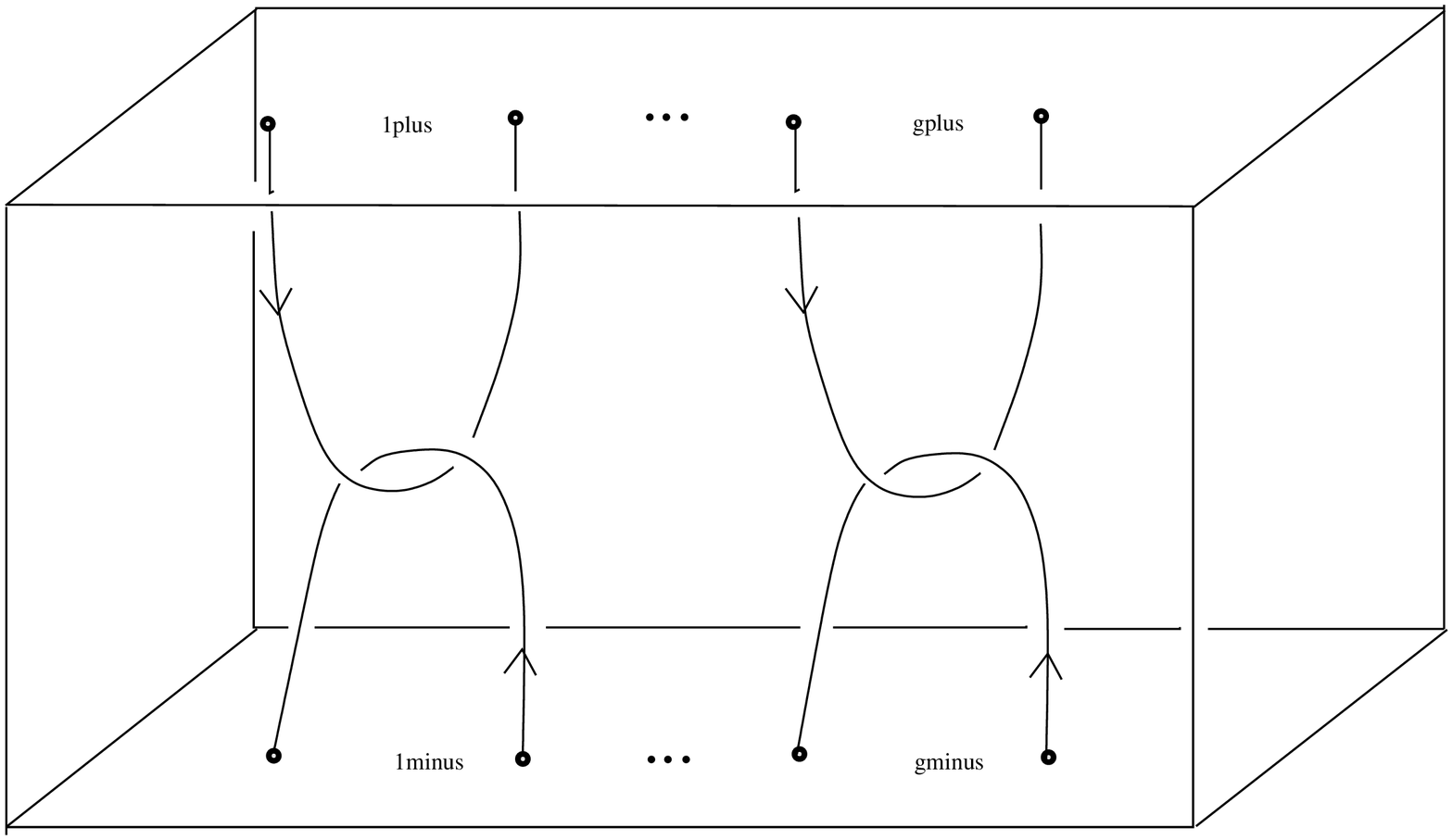}
\adjustrelabel <0cm,-0cm> {1minus}{$1$}
\adjustrelabel <0cm,-0cm> {1plus}{$1$}
\adjustrelabel <0cm,-0cm> {gminus}{$g$}
\adjustrelabel <0cm,-0cm> {gplus}{$g$}
\endrelabelbox}
\caption{The bottom-top tangle $\left([-1,1]^3,T_g\right)$ of type $(g,g)$.}
\label{fig:bt_tangle_T}
\end{figure}

\begin{figure}[h]
\centerline{\relabelbox \small
\epsfxsize 3truein \epsfbox{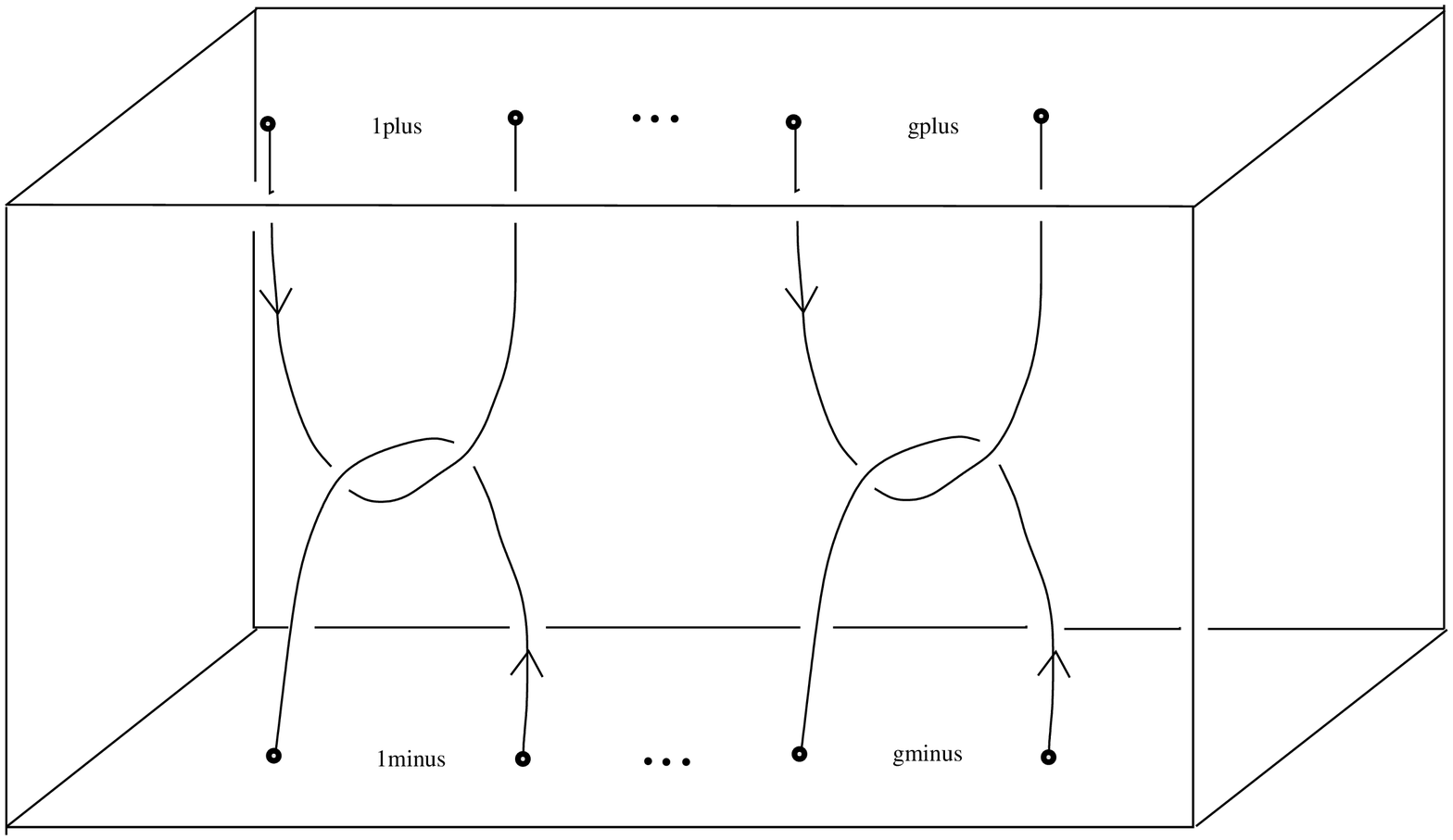}
\adjustrelabel <0cm,-0cm> {1minus}{$1$}
\adjustrelabel <0cm,-0cm> {1plus}{$1$}
\adjustrelabel <0cm,-0cm> {gminus}{$g$}
\adjustrelabel <0cm,-0cm> {gplus}{$g$}
\endrelabelbox}
\caption{The bottom-top tangle $\Id_g$ of type $(g,g)$.}
\label{fig:bt_tangle_id}
\end{figure}

Given bottom-top tangles $(B,\gamma)$ of type $(g_+,g_-)$ and $(C,\upsilon)$ of  type $(h_+,h_-)$
such that $g_+=h_-$, one obtains a new bottom-top tangle $(B,\gamma) \circ (C,\upsilon)$ of type $(h_+,g_-)$
as follows: The new manifold is
$$
(B \circ C)_{\gamma^+\ \cup\ T_{g_+}\ \cup\ \upsilon^-}
$$ 
i.e$.$ the composition $B\circ C$ in the category $\Cob$,
followed by surgery along the $(2g_+)$-component framed link  obtained by inserting the 
tangle $T_{g_+} \subset [-1,1]^3$ shown on Figure \ref{fig:bt_tangle_T} ``between'' $\upsilon^{-}$ and $\gamma^+$;
the new tangle is $(\upsilon^+,\gamma^-)$.
Thus, one obtains a category $\btT$ whose objects are non-negative integers $g$
and whose sets of morphisms $\btT(g_+,g_-)$ are bottom-top tangles of type $(g_+,g_-)$.
The identity of $g\geq 0$ in the category $\btT$ is drawn on Figure \ref{fig:bt_tangle_id}.

The category $\btT$ is monoidal (in the strict sense), with tensor product $\otimes$ given
by horizontal juxtaposition of bottom-top tangles in the $x$ direction.
So, to sum up, we have two operations on bottom-top tangles:
\begin{eqnarray*}
(B,\gamma) \circ (C,\upsilon) &:=& \begin{array}{|c|}
\hline \upsilon \subset C \\ \hline T_{g_+} \subset [-1,1]^3 \\ \hline \gamma \subset B \\ \hline \end{array}_{\begin{array}{r}
\hbox{\scriptsize \& do surgery along} \\ {}^{\gamma^+\ \cup\ T_{g_+}\ \cup\ \upsilon^-} \end{array}} \\
\hbox{and} \quad \quad
(B,\gamma) \otimes (C,\upsilon) & := & \begin{array}{|c|c|}  \hline \gamma \subset B  & \upsilon \subset C \\ \hline \end{array}\ .
\end{eqnarray*}

The study of bottom-top tangles is equivalent to the study of three-dimensional cobordisms.
More precisely, we have the following

\begin{theorem}
\label{th:iso}
There exists an isomorphism of monoidal categories $\hbox{D}: \btT \to \Cob$.
\end{theorem}

\noindent
This is very close to Kerler's presentation of cobordisms \cite{Kerler}, 
as well as the presentation of cobordisms between \emph{closed} surfaces described in \cite{CL2}.

\begin{proof}
Given a bottom-top tangle $(B,\gamma)$ of type $(g_+,g_-)$,
one obtains a cobordism from $F_{g_+}$ to $F_{g_-}$ by ``digging'' along the components of $\gamma$.
One gets a compact oriented connected $3$-manifold $M$ whose boundary is
identified with $\partial C_{g_-}^{g_+}$ via a map $m$
that is defined by means of the given identification $b:\partial C_0^0 \to \partial B$
and the framing of $\gamma$.
This construction is shown on Figure \ref{fig:bt_to_cob}, where
Remark \ref{rem:homeomorphism} applies to describe the
parameterizations $m_-:F_{g_-} \to M$ and $m_+: F_{g_+} \to M$ of the bottom and top surfaces.

\begin{figure}[h]
\centerline{\relabelbox \small
\epsfxsize 5.5truein \epsfbox{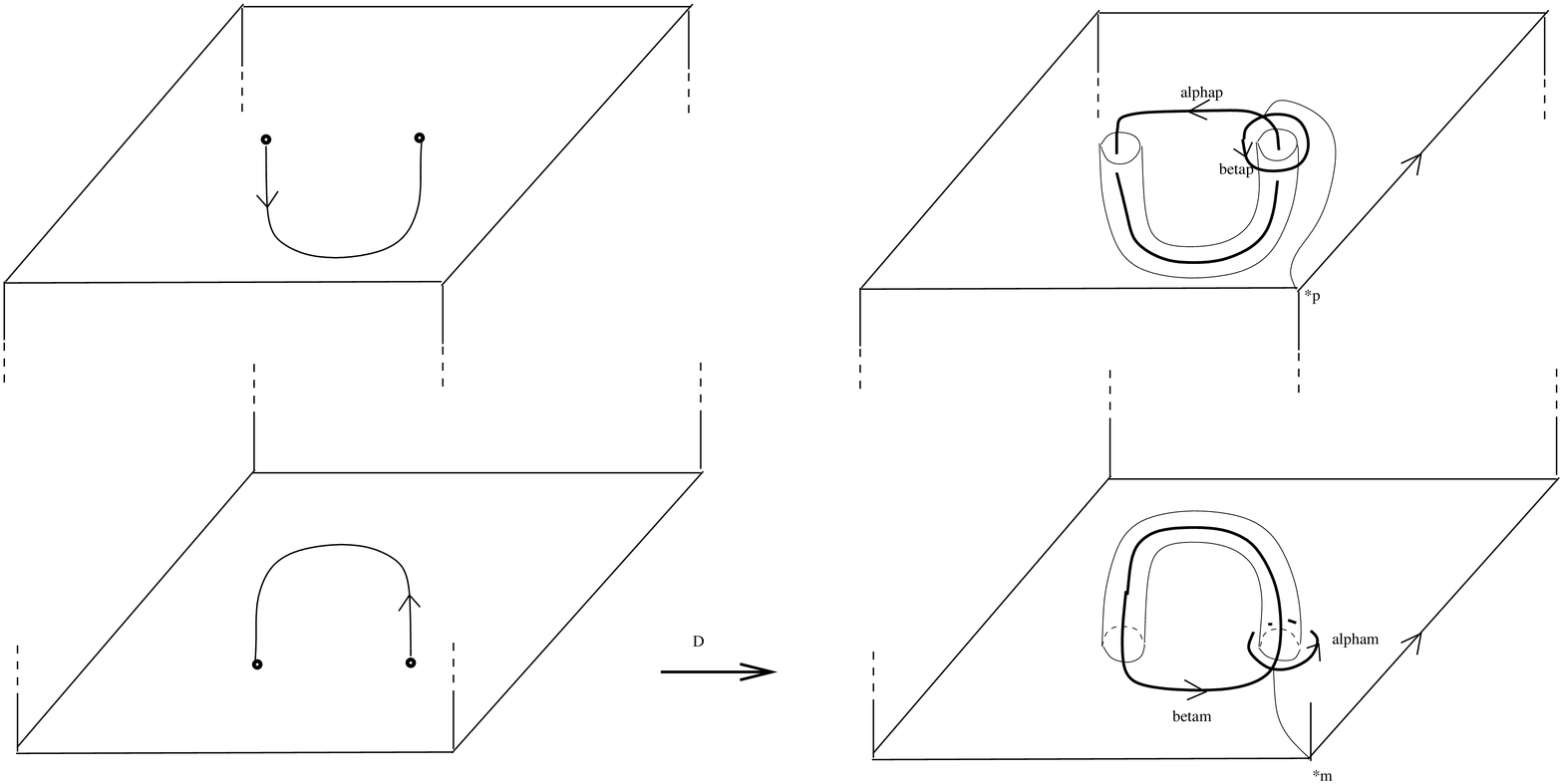}
\adjustrelabel <0cm,-0cm> {alpham}{$\alpha$}
\adjustrelabel <0cm,-0cm> {alphap}{$\alpha$}
\adjustrelabel <0cm,-0.1cm> {betam}{$\beta$}
\adjustrelabel <-0.1cm,-0.1cm> {betap}{$\beta$}
\adjustrelabel <0cm,-0cm> {*m}{$*$}
\adjustrelabel <0cm,-0cm> {*p}{$*$}
\adjustrelabel <0cm,-0cm> {D}{$\hbox{D}$}
\endrelabelbox}
\caption{From a bottom-top tangle to a cobordism (here $g_-=g_+=1$).}
\label{fig:bt_to_cob}
\end{figure}

The above construction is denoted by $\hbox{D}$, and we have to check its functoriality.
First, one easily sees that $\hbox{D}$ sends $\Id_g$ in $\btT$
to $\Id_g$ in $\Cob$, i.e$.$ Figure \ref{fig:bt_tangle_id} to Figure \ref{fig:cob_id}.
Next, let $(B,\gamma)$ and $(C,\upsilon)$ be bottom-top tangles 
of type $(g_+,g_-)$ and $(h_+,h_-)$ respectively, such that $g_+=h_-$,
and let $(M,m)$ and $(N,n)$ be the corresponding cobordisms by $\hbox{D}$.
A top component $\gamma_j^+$ of $\gamma$ may not bound a disk in $B$
but, after introducing a surgery link in $B$,
we can always assume that this is the case.
Then, $\gamma_j^+$ bounds a disk which is crossed
by a parallel family $X$ of strands, some
belonging to bottom components of $\gamma$ and some others
belonging to the added surgery link.
The rest of the argument is shown on Figure \ref{fig:functoriality}.

\begin{figure}[h]
\centerline{\relabelbox \small
\epsfxsize 5.5truein \epsfbox{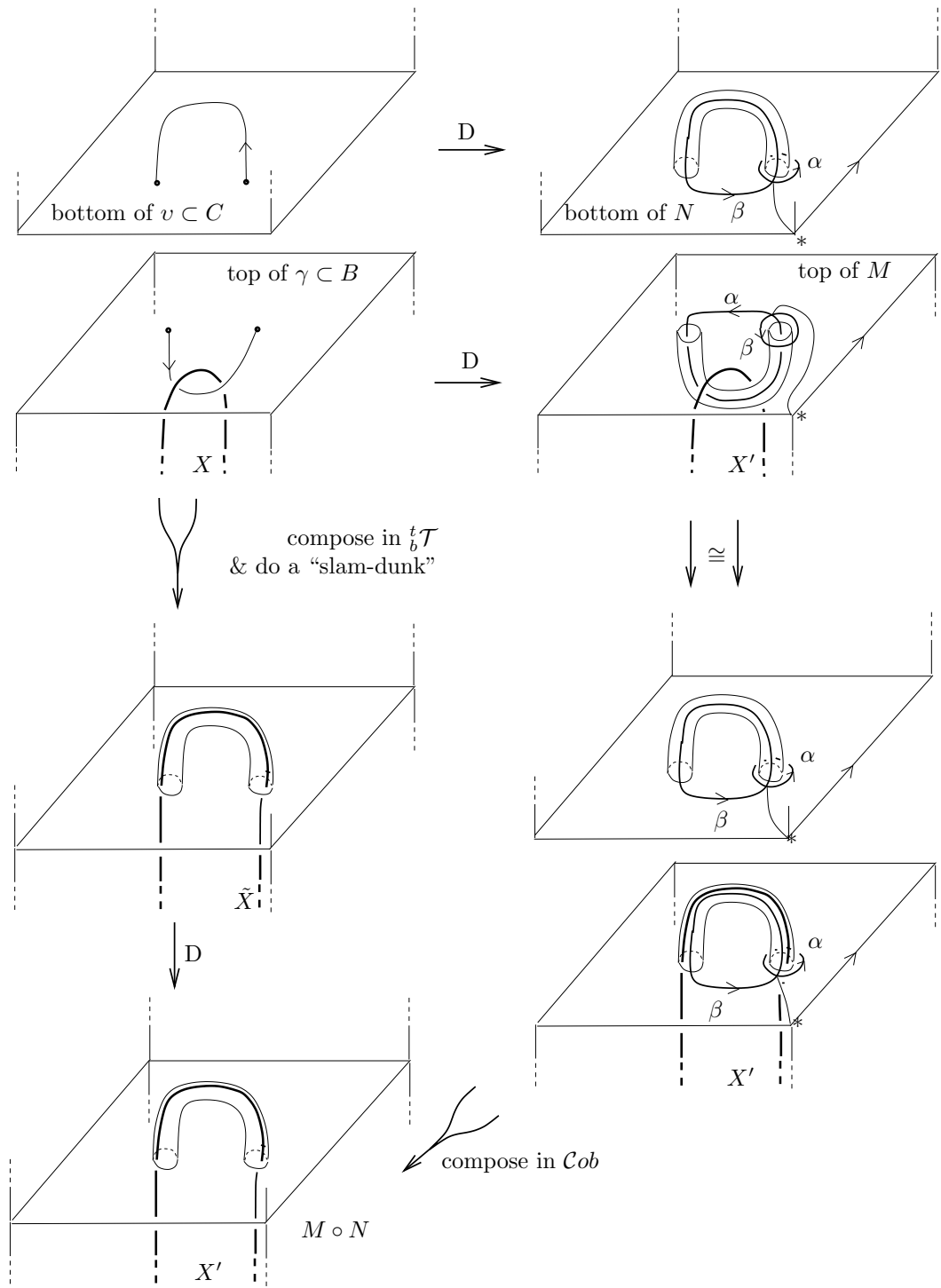}
\adjustrelabel <0cm,-0cm> {bottomofupsilon}{bottom of $\upsilon\subset C$}
\adjustrelabel <0cm,-0cm> {topofgamma}{top of $\gamma \subset B$}
\adjustrelabel <0cm,-0cm> {X}{$X$}
\adjustrelabel <0cm,-0cm> {i}{$\hbox{D}$}
\adjustrelabel <0cm,-0cm> {i2}{$\hbox{D}$}
\adjustrelabel <0cm,-0cm> {bottomofN'}{bottom of $N$}
\adjustrelabel <0cm,-0cm> {topofM'}{top of $M$}
\adjustrelabel <0cm,-0cm> {X2}{$X'$}
\adjustrelabel <0cm,-0cm> {alpham}{$\alpha$}
\adjustrelabel <0cm,-0cm> {alphap}{$\alpha$}
\adjustrelabel <0cm,-0.1cm> {betam}{$\beta$}
\adjustrelabel <-0.1cm,-0.1cm> {betap}{$\beta$}
\adjustrelabel <0cm,-0cm> {*m}{$*$}
\adjustrelabel <0cm,-0cm> {*p}{$*$}
\adjustrelabel <0cm,-0cm> {homeo}{$\cong$}
\adjustrelabel <-0.1cm,0.1cm> {*m2}{$*$}
\adjustrelabel <-0.1cm,0.1cm> {*p2}{$*$}
\adjustrelabel <0cm,-0cm> {alpham2}{$\alpha$}
\adjustrelabel <0cm,-0cm> {alphap2}{$\alpha$}
\adjustrelabel <0cm,-0.1cm> {betam2}{$\beta$}
\adjustrelabel <-0.1cm,-0.1cm> {betap2}{$\beta$}
\adjustrelabel <0cm,-0cm> {X3}{$X'$}
\adjustrelabel <0cm,-0cm> {compositioninCob}{compose in $\Cob$}
\adjustrelabel <0cm,-0cm> {X4}{$X'$}
\adjustrelabel <0cm,-0cm> {N'M'}{$M\circ N$}
\adjustrelabel <0cm,-0cm> {i3}{$\hbox{D}$}
\adjustrelabel <0cm,-0cm> 
{compositioninBT}{\begin{tabular}{r}compose in $\btT$\\ \& do a ``slam-dunk'' \end{tabular}}
\adjustrelabel <0cm,-0cm> {X5}{$\tilde{X}$}
\endrelabelbox}
\caption{Functoriality of the map $\hbox{D}$ (here $g_+=h_-=1$).
The ``slam-dunk'' move is recalled on Figure \ref{fig:slam-dunk}.}
\label{fig:functoriality}
\end{figure}

\begin{figure}[h]
\centerline{\relabelbox \small
\epsfxsize 2truein \epsfbox{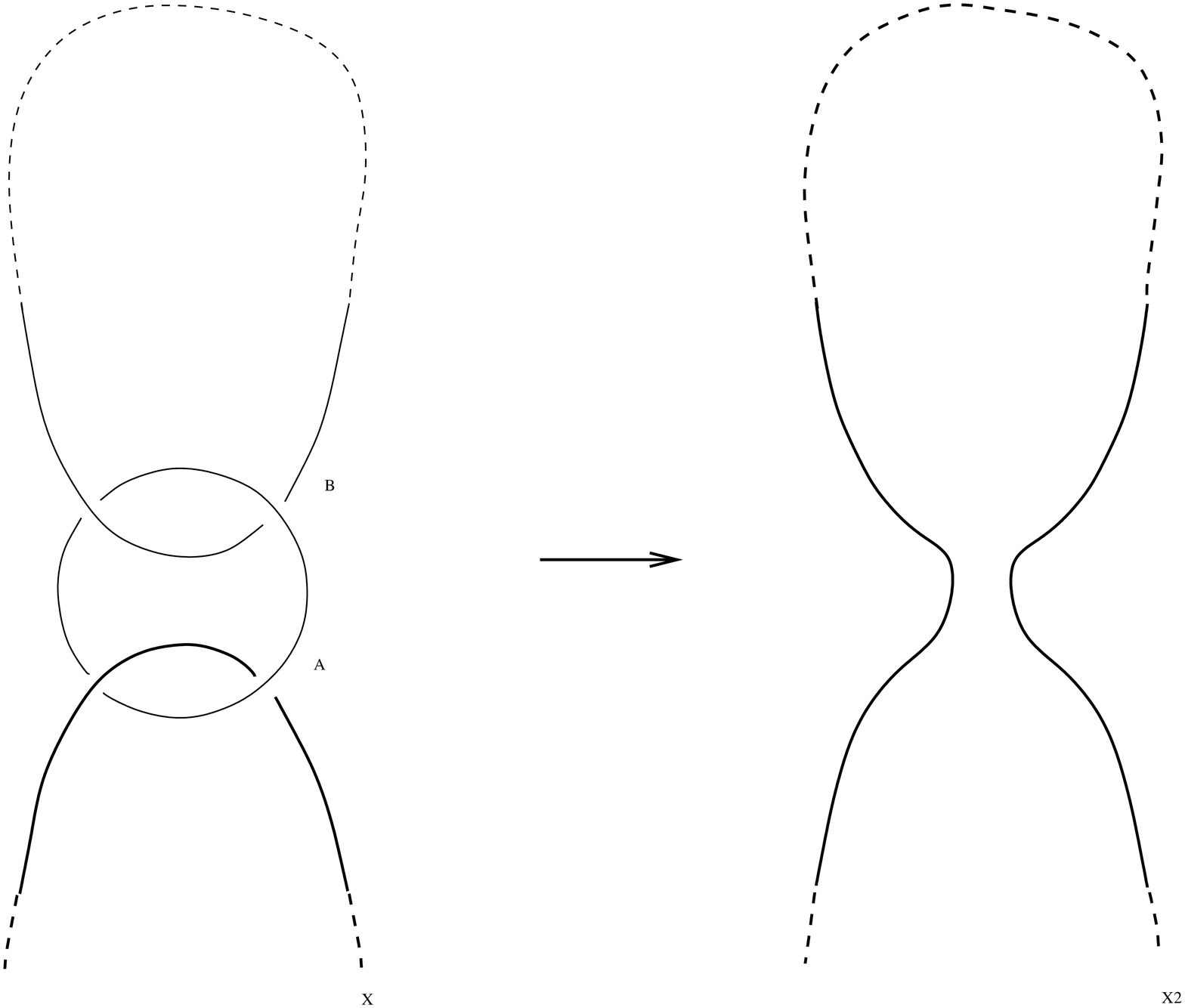}
\adjustrelabel <0cm,-0cm> {A}{$R$}
\adjustrelabel <0cm,-0cm> {B}{$V$}
\adjustrelabel <0cm,-0cm> {X}{$X$}
\adjustrelabel <0cm,-0cm> {X2}{$\tilde{X}$}
\endrelabelbox}
\caption{The \emph{slam-dunk} move: Surgery is performed along the two-component framed link $(R,V)$,
it produces a homeomorphic manifold and the corresponding homeomorphism changes
a parallel family of strands $X$ to $\tilde{X}$.}
\label{fig:slam-dunk}
\end{figure}

Thus, one gets a functor $\hbox{D}: \btT \to \Cob$, which obviously preserves the tensor product.
Also, $\hbox{D}$ has an inverse functor defined by gluing $2$-handles as follows:
Given a cobordism $(M,m)$ from $F_{g_+}$ to $F_{g_-}$,
one obtains a manifold $B$ with $\partial B \cong \partial C_0^0$
by attaching one $2$-handle along each curve $m_-(\alpha_i)$ of the bottom surface,
and along each curve $m_+(\beta_i)$ of the top surface; the co-cores of those $2$-handles
define a bottom-top tangle $\gamma$ in $B$ of type $(g_+,g_-)$.
\end{proof}

The isomorphism $\Cob \simeq \btT$ allows one to regard $\LCob$,
and a fortiori $\sLCob$, as subcategories of $\btT$.
It follows from the definitions that
a bottom-top tangle $(B,\gamma)$ of type $(g_+,g_-)$ belongs to $\sLCob(g_+,g_-)$ 
if, and only if, $B$ is the standard cube $C_0^0=[-1,1]^3$ 
and $\gamma^+$ is the trivial $g_+$-component top tangle. 
In order to characterize $\LCob$ in $\btT$, we need the following 

\begin{definition}
Let $(B,\gamma)$ be a bottom-top tangle in a homology cube.\footnote{
A \emph{homology cube} $B$ is a cobordism $(B,b)$ 
from $F_0$ to $F_0$ such that $H_*(B) \simeq H_*([-1,1]^3)$.}
The \emph{linking matrix} of $\gamma$ in $B$ is the matrix,
whose rows and columns are indexed by the set of connected components of $\gamma$, defined by
$$
\Lk_B(\gamma) := \Lk_{\hat{B}}(\hat{\gamma}).
$$ 
Here, $\hat{B}:= B \cup_b \left(S^3 \setminus [-1,1]^3\right)$ 
is the homology sphere obtained by ``recapping'' $B$, 
$\hat{\gamma}$ is the framed oriented link in $\hat{B}$ 
whose component $\hat{\gamma}_i^\pm$ is $\gamma_i^\pm$ union with a small arc
connecting $p_i\times (\pm 1)$ to $q_i \times (\pm 1)$ in the $x$ direction
and $\Lk_{\hat{B}}(\hat{\gamma})$ denotes the usual linking matrix of $\hat{\gamma}$ in $\hat{B}$.
\end{definition}

\begin{lemma}
\label{lem:bh_bt}
A bottom-top tangle $(B,\gamma)$ of type $(g_+,g_-)$ belongs to $\LCob(g_+,g_-)$ if, and only if,
$B$ is a homology cube and the linking matrix of $\gamma^+$ in $B$ is trivial.
\end{lemma}

\begin{proof}
Let $(M,m)$ be the cobordism from $F_{g_+}$ to $F_{g_-}$ corresponding
to the bottom-top tangle $(B,\gamma)$ by Theorem \ref{th:iso}.
Recall that $M$ is the complement of a tubular neighborhood of $\gamma$ in $B$. 
Observe that $\alpha_k^-:=m_-(\alpha_k)$,
$\beta_k^-:=m_-(\beta_k)$, $\alpha_k^+:=m_+(\alpha_k)$
and $\beta_k^+:=m_+(\beta_k)$ are respectively oriented meridian of $\gamma_k^-$,
oriented longitude of $\gamma_k^-$, oriented longitude of $\gamma_k^+$
and oriented meridian of $\gamma_k^+$. The condition
$$
H_1(M)=m_{-,*}(A_{g_-}) \oplus m_{+,*}(B_{g_+})
$$
is equivalent to the condition that $B$ is a homology cube. Assuming this condition,
we have that $H_1(M)$ is free Abelian of rank $g_-+g_+$ with
basis given by the oriented meridians, namely $(\alpha^-_1,\dots,\alpha^-_{g_-},\beta_1^+,\dots,\beta^+_{g_+})$.
Since the columns of the linking matrix of $\gamma$ in $B$ express
how the oriented longitudes $\beta^-_1,\dots,\beta^-_{g_-},\alpha^+_1,\dots,\alpha^+_{g_+}$ expand in that basis,
the linking matrix of $\gamma^+$ is trivial if and only if $m_{-,*}(A_{g_-}) \supset m_{+,*}(A_{g_+})$. 
\end{proof}

\begin{remark}
The same proof shows that a bottom-top tangle $(B,\gamma)$ of type $(g_+,g_-)$
 belongs to $\QLCob(g_+,g_-)$ if, and only if,
$B$ is a $\Q$-homology cube and the linking matrix of $\gamma^+$ in $B$ is trivial.
\end{remark}

According to the previous lemma, 
the following definition makes sense and it will be used later:

\begin{definition}
The \emph{linking matrix} of a Lagrangian cobordism $M=(M,m)$ is 
$$
\Lk(M) := \Lk_B(\gamma)
$$
where $(B,\gamma)$ is the corresponding bottom-top tangle in a homology cube.
\end{definition}

\vspace{0.5cm}

\section{The Kontsevich--LMO invariant of tangles in homology cubes}

In this section, we review the Kontsevich--LMO invariant of tangles in homology cubes,
which will play the lead role in the next sections.

\subsection{Spaces of Jacobi diagrams}

\label{subsec:diagrams}

First of all, we need to recall some definitions and notations  about Jacobi diagrams, 
which come mainly from \cite{BN,BGRT1,BGRT2}. The reader is refered to those papers for details.

A \emph{uni-trivalent graph} $D$ is a finite graph whose vertices have valence $1$
(\emph{external} vertices) or $3$ (\emph{internal} vertices).
It is \emph{vertex-oriented} if each internal vertex comes with a cylic
order of its incident edges.
One defines the \emph{internal degree}, the \emph{external degree} and the \emph{degree} to be
$$
\left\{\begin{array}{rcl}
\ideg(D) & := & \hbox{number of internal vertices of $D$}\\
\edeg(D) & := & \hbox{number of external vertices of $D$}\\
\deg(D) & := &\left(\ideg(D)+\edeg(D)\right)/2.
\end{array} \right.
$$

In the sequel, let $X$ be a compact oriented $1$-manifold 
and let $C$ be a finite set.

\begin{example}
If $S$ is a finite set, $X$ can be the disjoint union indexed by $S$ of oriented circles (respectively intervals),
which is denoted by $\circlearrowleft^S$ (respectively $\uparrow^S$). 
Conversely, if $L$ is a compact oriented $1$-manifold, $C$ can be the set of its connected components,
which is denoted by $\pi_0(L)$.
\end{example}

\begin{example}
If $n$ is a positive integer and $*$ is an extra symbol (such as $+,-,$ etc.), 
$C$ can be the finite set $\{1^*,\dots,n^*\}$, which is denoted by $\set{n}^*$. 
\end{example}

A \emph{Jacobi diagram} $D$ based on $(X,C)$
is a vertex-oriented uni-trivalent graph
whose external vertices are either embedded into $X$ or are colored with elements from $C$.
Let $X'$ be another compact oriented $1$-manifold whose $\pi_0(X')$ is identified with
$\pi_0(X)$ and let $C'$ be another finite set identified with $C$.
Then, two Jacobi diagram $D$ and $D'$  based on $(X,C)$ and $(X',C')$ respectively
are \emph{equivalent} if there exists a homeomorphism $f:X\cup D\to X'\cup D'$
sending $X$ to $X'$ in such a way that orientations
and connected components are preserved, and sending $D$ to $D'$
in such a way that vertex-orientations and colors are respected. 
In the sequel, Jacobi diagrams $(X,C)$ are considered up to equivalence.
In pictures, the $1$-manifold part $X$ is drawn with bold lines
while the graph part $D$ is drawn with dashed lines, and the vertex-orientation is given 
by the trigonometric orientation of the blackboard.

\begin{example}
A \emph{strut} is a Jacobi diagram reduced to a single edge and whose vertices are colored with $C$. 
It is pictured as $\strutgraph{a}{b}$ where $a,b\in C$. 
\end{example}

The spaces of Jacobi diagrams that are needed in this paper, are always of the form  
$$
\A\left(X,C\right):=
\frac{\Q\cdot \left\{\hbox{Jacobi diagrams based on $(X,C)$}\right\}}
{\hbox{AS, IHX, STU}}
$$
where the AS, IHX and STU relations are as usual \cite{BN}:\\

\begin{center}
{\relabelbox \small
\epsfxsize 5.8truein \epsfbox{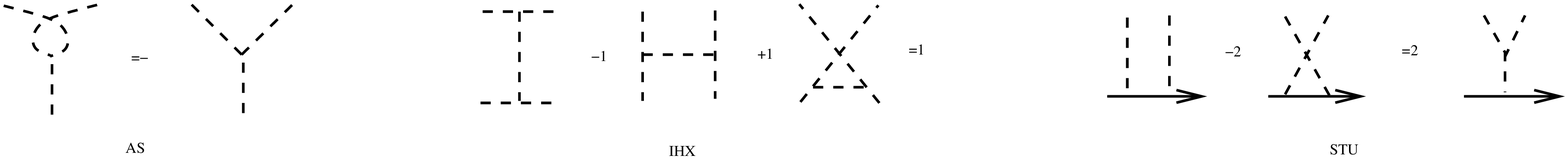}
\adjustrelabel <0cm,-0.2cm> {AS}{AS}
\adjustrelabel <0cm,-0.2cm> {IHX}{IHX}
\adjustrelabel <0cm,-0.2cm> {STU}{STU}
\adjustrelabel <0cm,-0cm> {=-}{$= \ -$}
\adjustrelabel <0cm,-0cm> {-1}{$-$}
\adjustrelabel <0cm,-0.05cm> {+1}{$+$}
\adjustrelabel <0cm,-0.1cm> {=1}{$=0$}
\adjustrelabel <0cm,0.cm> {-2}{$-$}
\adjustrelabel <0cm,-0.05cm> {=2}{$=$}
\endrelabelbox}
\end{center}

\vspace{-0.2cm}

\begin{example}
Any rational  matrix $M=\left(m_{ij}\right)_{i,j\in C}$,
whose rows and columns are indexed by $C$, defines a linear combination of Jacobi diagrams:
$$
M := \sum_{i, j \in C} m_{ij} \cdot \strutgraph{i}{j}  \ \in \A(X,C).
$$
\end{example}

The relations AS, IHX and STU being homogeneous with respect to the degree,
$\A\left(X,C\right)$ is graded by the degree of Jacobi diagrams:
The degree completion of $\A\left(X,C\right)$ is denoted the same way.
The STU relation is not homogeneous with respect to the internal degree:
Nevertheless, an element $x \in \A(X,C)$ is said to have 
\emph{i-filter} at least $n$ if it can be written as a linear sum of Jacobi diagrams
with at least $n$ internal vertices. 

Assume now that $X$ is empty, so that the STU relation becomes trivial. 
The disjoint union operation $\sqcup$ of Jacobi diagrams makes $\A\left(C\right)$ a commutative algebra,
whose identity element is the empty diagram $\varnothing$.
The exponential $\exp_\sqcup (x)$ of an $x\in \A(C)$, with respect to the multiplication $\sqcup$,
will often be denoted by
$$
[x] := \sum_{n\geq 0} \frac{1}{n!} \cdot \underbrace{x \sqcup \cdots \sqcup x}_{n\ {\rm \scriptsize{times}}}.
$$
The sub-space of $\A\left(C\right)$ spanned by Jacobi diagrams 
without strut (respectively, with only struts) is denoted by $\AY\left(C\right)$ (respectively, $\As(C)$), 
and is identified with the quotient of $\A\left(C\right)$ by the ideal generated by struts
(respectively, by Jacobi diagrams with at least one internal vertex). So, one has two projections
$$
\xymatrix@R=0pt{
{\As(C)} & {\A(C)} \ar@{->>}[l] \ar@{->>}[r] & {\AY(C)}\\
{x^s} & {x} \ar@{|->}[r] \ar@{|->}[l] & {x^Y}
}
$$
called the \emph{$s$-reduction} and the \emph{$Y$-reduction} respectively.
Observe that the degree completion of $\AY(C)$ (still denoted by $\AY(C)$) 
is canonically isomorphic to its i-degree completion.

The usual coproduct $\Delta$, defined by 
$$
\Delta(D) := \sum_{D = D' \sqcup D''} D' \otimes D'',
$$
enhances $\A(C)$ to a co-commutative Hopf algebra, 
whose counit is the linear map $\varepsilon: \A(C) \to \Q$ defined by 
$\varepsilon(D) := \delta_{D,\varnothing}$ for all Jacobi diagram $D$.
The space of primitives elements of $\A(C)$ is the sub-space $\A^c\left(C\right)$ 
of non-empty $c$onnected diagrams.
The following lemma is well-known,
and is deduced from the fact that group-like elements are exponentials of primitive elements:

\begin{lemma}
\label{lem:group-like}
An $x\in\A(C)$ is group-like if, and only if,
the $s$-reduction and the $Y$-reduction of $x$ are group-like and such that $x = x^s \sqcup x^Y$.
\end{lemma}

\noindent
Observe that a group-like element of $\As(C)$ is necessarily of the form 
$[M]$ where $M$ is a $C \times C$ matrix with rational entries.\\

We now recall some operations on Jacobi diagrams. Let $S$ be another finite set, disjoint from $C$.
There is defined in \cite{BN} a diagrammatic analogue of the Poincar\'e--Birkhoff--Witt isomorphism
$$
\chi_S: \A(X, C\cup S) \longrightarrow \A(X \uparrow^S,C).
$$
For $D$ a Jacobi diagram, $\chi_S(D)$ is the average of all possible ways of attaching, for all color $s\in S$,
the $s$-colored external vertices of $D$ to the $s$-indexed interval in $\uparrow^S$. 

A Jacobi diagram $D \in \A(X,C \cup S)$ is said to be \emph{$S$-substantial} if it contains no strut both of whose vertices
are colored by $S$. For all Jacobi diagrams $D,E \in \A(X,C \cup S)$ such that $D$ or $E$ is $S$-substantial, 
we define as in \cite{BGRT1}
$$
\langle E,D \rangle_{S} := 
\left(\begin{array}{c}
\hbox{sum of all ways of gluing the $s$-colored vertices of $D$}\\
\hbox{to the $s$-colored vertices of $E$, for all color $s\in S$}
\end{array}\right) \ \in \A(X,C).
$$
Next combinatorial result will be used at several places:
\begin{theorem}[Jackson--Moffatt--Morales \cite{JMM}]
\label{th:JMM}
For all group-like elements $D,E \in \A(C \cup S)$ such that $D$ or $E$ is $S$-substantial, 
$\langle E,D \rangle_{S}$ is group-like. In other words,
$$
\langle E,D \rangle_{S} = \exp_\sqcup \left(\hbox{connected part of }\langle E,D \rangle_{S}\right).
$$
\end{theorem}

A linear combination of Jacobi diagrams $G \in \A(X,C \cup S)$ is \emph{Gaussian} in the \emph{variable} $S$
if it can be written in the form
$$
G = \left[L/2 \right] \sqcup P
$$
where $P$ is $S$-substantial and $L$ is a rational symmetric $S\times S$ matrix.
When $\det(L)\neq 0$, the Gaussian $G$ is \emph{non-degenerate}.
In this case, the \emph{formal Gaussian integral} of $G$ \emph{along} $S$ is defined in \cite{BGRT1} by
$$
\int_S G := \left\langle \left[-L^{-1}/2\right], P \right\rangle_S \ \in \A(X,C).
$$

\begin{remark}
\label{rem:inv_link}
Some \emph{$S$-link} relations in $\A(X, C\cup S)$ 
are defined in \cite[\S 5.2]{BGRT2} so that
$$
\xymatrix{
{\A(X, C\cup S)} \ar@{->>}[d] \ar[r]_-{\simeq}^-{\chi_S} & {\A(X \uparrow^S,C)} \ar@{->>}[d]^-{\hbox{\footnotesize closure}}\\
{\frac{\A(X, C\cup S)}{\hbox{\footnotesize $S$-link}}} \ar@{-->}[r]^-{\simeq}_-{\chi_S} & {\A(X \circlearrowleft^S,C).}
}
$$
Let $G = \left[L/2 \right] \sqcup P$ and 
$G' =  \left[L'/2 \right] \sqcup P'$ be non-degenerate Gaussians 
in the variable $S$. According to \cite[Prop. 2.2]{BL}, if $G$ varies from $G'$ by some $S$-link relations, then 
$$
\int_S G  = \int_S G'.
$$
This fact will allow us to forget about link relations in computations.
\end{remark}

\subsection{The category $\qTCub$ of $q$-tangles in homology cubes}

We now define the domain of the Kontsevich--LMO invariant.

\begin{definition}
\label{def:tangle}
By a \emph{tangle}, we mean an equivalence class of couples $(B,\gamma)$ 
where $B$ is a cobordism from $F_0$ to $F_0$
and where $\gamma$ is a framed oriented tangle in $B$ whose boundary points 
(if any) are either in the bottom surface, or in the top surface.
We also assume that those points (at their respective levels) are uniformly distributed along
the segment $[-1,1] \times 0 \times 0$ in $[-1,1]^2 \times 0=F_0$. 
\end{definition}

\begin{example}
A bottom-top tangle is a tangle.
\end{example}

\noindent
If one associates to each boundary point of $\gamma$ the sign
\begin{equation}
\label{eq:words}
\left\{\begin{array}{ll}
+  & \hbox{if the orientation of $\gamma$ at that point goes ``downwards''}\\
-  & \hbox{if the orientation of $\gamma$ at that point goes ``upwards''}\\
\end{array} \right.
\end{equation}
one gets two associative words in the letters $(+,-)$,
one for the bottom and another one for the top.

\begin{definition}
\label{def:q-tangle}
A \emph{q-tangle} is a tangle $(B,\gamma)$
together with some lifts $w_t(\gamma)$ and $w_b(\gamma)$
to the free non-associative magma generated by $(+,-)$
of the top and bottom words defined by $\gamma$ in the free monoid generated by $(+,-)$.
\end{definition}

\begin{remark}
This definition slightly extends the notion of ``$q$-tangle'' given in \cite{LM1,LM2},
where the cobordism $B$ is required to be the cube $[-1,1]^3$.
\end{remark}

Given two $q$-tangles $(B,\gamma)$ and $(C,\upsilon)$ such that $w_t(\gamma)=w_b(\upsilon)$,
one can form the new tangle $(B \circ C, \gamma \cup \upsilon)$ and equip it
with the non-associative words $w_t(\upsilon\cup \gamma) := w_t(\upsilon)$
and $w_b(\upsilon \cup \gamma) := w_b(\gamma)$. Thus, one obtains
a category whose objects are non-associative words in the letters $(+,-)$
and whose  morphisms are $q$-tangles. There is a tensor product $\otimes$
given by horizontal juxtaposition of $q$-tangles in the $x$ direction.
So, we get a monoidal category (in the non-strict sense).

In the sequel, we will only need two subcategories of this: 
The monoidal category of $q$-tangles in \emph{homology cubes}, which we denote by $\qTCub$,
and the monoidal category of $q$-tangles in the standard cube $[-1,1]^3$, 
which we denote by $\mathcal{T}_q$.

\subsection{The category $\A$ of Jacobi diagrams on $1$-manifolds}

We now define the codomain of the Kontsevich--LMO invariant.

For all associative words $u$ and $v$ in the letters $(+,-)$, we define
$\A(v,u)$ to be the union of all the spaces $\A(X)$,
where $X$ runs over homeomorphism classes of compact oriented $1$-manifolds 
whose boundary is  identified with the set of letters of $u$ and $v$ as follows:
A positive point of $\partial X$ should be assigned either to a $-$ letter in $v$ or a 
$+$ letter in $u$, and vice-versa for a negative point of $\partial X$. 

\begin{example}
For all associative word $w$ in the letters $(+,-)$, 
we denote by $\downarrow^w$ the $1$-manifold obtained 
by taking one copy of $\downarrow$ or $\uparrow$
for each letter $+$  or $-$ respectively read in the word $w$.
Thus, $\A(w,w)$ contains $\A(\downarrow^w)$.
\end{example}

Given $a\in \A(X) \subset \A(v,u)$ and $b\in \A(Y) \subset \A(w,v)$, 
one obtains a new element $a \circ b$ of $\A(X \cup Y) \subset \A(w,u)$ by gluing
$b$ on the ``top'' of $a$.
Thus, one gets a category $\A$ whose objects are associative words in the letters $(+,-)$ 
and whose morphisms are linear combinations of Jacobi diagrams based on compact oriented $1$-manifolds. 
The identity of $w$ in the category $\A$ is the empty Jacobi diagram on $\downarrow^w$.

There is a tensor product $\otimes$ given by juxtaposition of Jacobi diagrams:
So, $\A$ is a monoidal category (in the strict sense).

\begin{notation}
\label{not:doubling}
Recall from \cite{LM1,LM2} that there are an ``orientation-reversal'' map $S: \A(X \downarrow) \to \A(X \uparrow)$
and a ``doubling'' map $\Delta : \A(X \downarrow) \to \A(X \downarrow \downarrow)$.
If $w$ is a word of length $g:=|w|$ in the letters $(+,-)$ 
and if $w_1,\dots,w_g$ are extra words, we denote by 
$$
\Delta_{w_1,\dots,w_g}^w: \A\left(X \downarrow^w\right) 
\longrightarrow \A\left(X \downarrow^{w_1} \cdots \downarrow^{w_g}\right)
$$
the map obtained by applying, for each $i=1,\dots,g$, $(|w_i|-1)$ times the $\Delta$ map to
the $i$-th component of $\downarrow^w$ and by applying the $S$ map to each new interval
whose corresponding letter in $w_i$ does not agree with the $i$-th letter of $w$. 
For example, we have $\Delta_{++}^+=\Delta$ and $\Delta_{-}^+=S$.
\end{notation}

\subsection{The Kontsevich integral $Z$}

\label{subsec:Kontsevich}

Le and Murakami have extended in \cite{LM1,LM2} the Kontsevich integral 
of links in $S^3$ to a tensor-preserving functor
$$
\widehat{Z}_f: \mathcal{T}_q \longrightarrow \A, \ \gamma \longmapsto \widehat{Z}_f(\gamma).
$$
At the level of objects, $\widehat{Z}_f$ just forgets the parenthesizings.
At the level of morphisms, $\widehat{Z}_f$ is determined 
by its values on the  ``elementary'' $q$-tangles, namely
$$
\begin{array}{llcll}
\widehat{Z}_f\left( \begin{array}{c} {}_{(++)} \\ \pcrossing \\ {}^{(++)} \end{array} \right) 
:= \left[\frac{1}{2} \figtotext{20}{20}{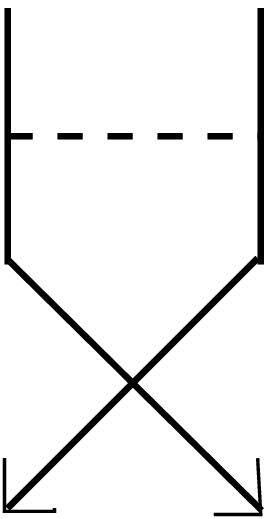}\right] \! & \! \in \A\left( \crossing \right) 
& \ &
\widehat{Z}_f\left( \begin{array}{c} {}_{(++)} \\ \ncrossing \\ {}^{(++)} \end{array} \right) 
:= \left[- \frac{1}{2} \figtotext{20}{20}{chord_on_crossing.eps}\right] \! & \! \in \A\left( \crossing \right) \\
\widehat{Z}_f\left(\begin{array}{c}   \capleft \\ {}^{(+-)}  \end{array}\right) := 
\begin{array}{c}{\relabelbox \small
\epsfxsize 0.5truein \epsfbox{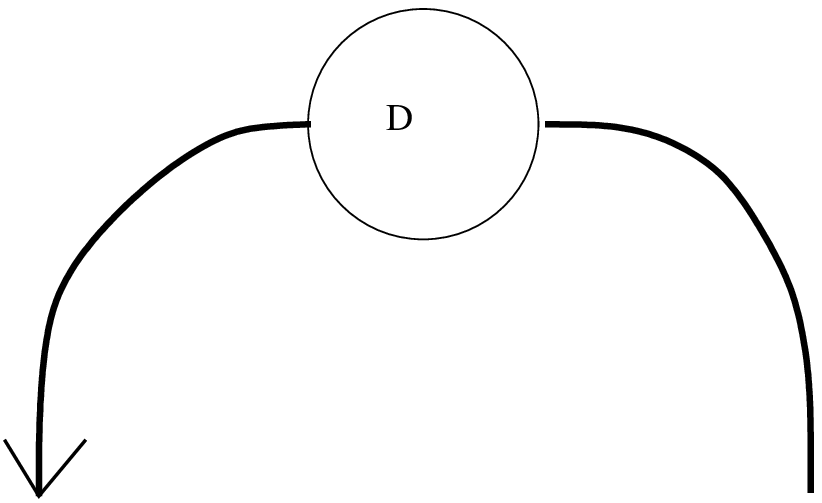}
\adjustrelabel <-0.05cm,-0cm> {D}{${}_{\nu}$}
\endrelabelbox} \end{array} \! & \! \in \A\left(\capleft\right)
& \ &
\widehat{Z}_f\left(\begin{array}{c} {}_{(+-)} \\ \cupright \end{array}\right) := 
\figtotext{20}{15}{cupright.eps} \! & \! \in \A\left(\cupright\right)
\end{array}
$$
and, for all non-associative words $u,v,w$, by 
$$
\widehat{Z}_f \left(
\begin{array}{c} {}_{(u\ (vw))} \\ \downarrow \swarrow \downarrow \\ {}^{((uv)\ w)}\end{array}
\right) := \Delta_{u,v,w}^{+++}(\Phi) \ \in \A\left(\downarrow^{u v w}\right).
$$
Here, $\Phi \ \in \A(\downarrow \downarrow \downarrow)$ is a Drinfeld associator which has to be choosen,
while $\nu = \widehat{Z}_f(\circlearrowleft_0) \in \A(\circlearrowleft) \simeq \A(\uparrow)$ 
is given as the value of the Kontsevich integral on the $0$-framed unknot.

\begin{remark}
We agree to fix a Drinfeld associator with \emph{rational} coefficients.
Nonetheless, if we had  defined Jacobi diagrams with complex coefficients, 
then we could have worked with the KZ associator as well.
\end{remark}

In this paper, we prefer for technical convenience to modify $\widehat{Z}_f$ as follows: 
For all $q$-tangle $\gamma$ in $[-1,1]^3$
with connected components $\gamma_1,\dots,\gamma_l$, we set
$$
Z(\gamma) := \widehat{Z}_f(\gamma)\ 
\sharp_1\ \nu^{d(\gamma_1)}\ \sharp_2\ \cdots\ \sharp_l\ \nu^{d(\gamma_l)}  
$$
where $\sharp_i$ means that a connected sum
of an element of $\A(\circlearrowleft)$ is taken with the $i$-th component of $\gamma$,
and where $d(\gamma_i)$ is $-1$, $0$ or $1$ if the component $\gamma_i$
is of type ``bottom-bottom'', ``bottom-top'' or ``top-top'' respectively.
In other words, $Z$ only differs from $\widehat{Z}_f$ 
by the values it takes on the ``cap'' and the ``cup'':
$$
Z\left(\begin{array}{c}\capleft\\ {}^{(+-)} \end{array}\right) 
= \figtotext{20}{15}{capleft.eps} \ \in  \A\left(\capleft\right)
\quad \quad  \quad \quad \quad 
Z\left(\begin{array}{c} {}_{(+-)} \\ \cupright  \end{array}\right) = 
\begin{array}{c}{\relabelbox \small
\epsfxsize 0.5truein \epsfbox{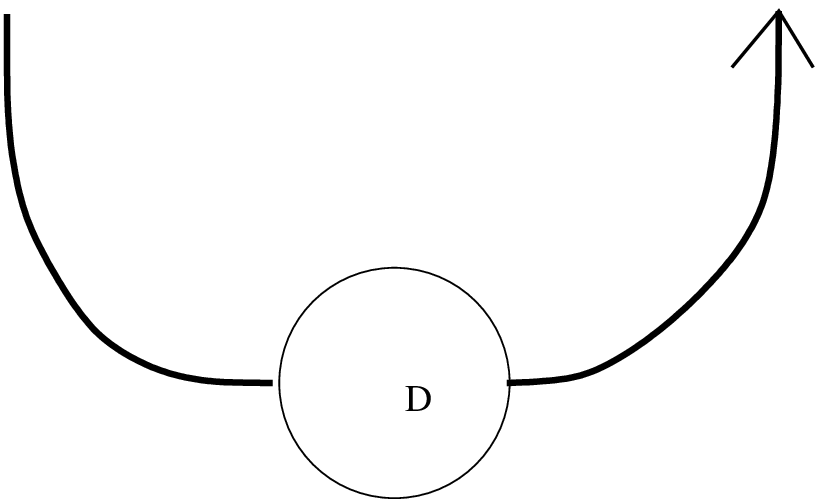}
\adjustrelabel <-0.05cm,-0cm> {D}{${}_{\nu}$}
\endrelabelbox} \end{array}  \ \in \A\left(\cupright\right).
$$
In the sequel, the \emph{Kontsevich integral} will refer to this tensor-preserving functor
$$
Z: \mathcal{T}_q \longrightarrow \A, \ \gamma \longmapsto Z(\gamma).
$$

\subsection{The Kontsevich--LMO invariant $Z$}

\label{subsec:Kontsevich-LMO}

We now construct from the Kontsevich integral a tensor-preserving functor 
$$
Z: \qTCub \longrightarrow \A
$$
with the following properties:
\begin{itemize}
\item For all $q$-tangle $\gamma$ in $[-1,1]^3$, 
$Z([-1,1]^3,\gamma)$ coincides with the Kontsevich integral $Z(\gamma)$, as normalized in \S \ref{subsec:Kontsevich}.
\item For all homology cube $B$, $Z(B,\varnothing)$ coincides with the Le--Murakami--Ohtsuki invariant
$\Omega(\hat{B})$ of the homology sphere $\hat{B}$, as defined in \cite{LMO}.
\end{itemize}
The fact that the LMO invariant and the Kontsevich integral can be unified into a single invariant
of $q$-tangles in homology cubes is well-known to experts. 
We do this below using the Aarhus formalism  \cite{BGRT1,BGRT2}.

For this, we fix a few notations. Given a $q$-tangle $L\cup \gamma$ in $[-1,1]^3$ whose connected components
are split into two parts, $L$ and $\gamma$, we set
$$
Z(L^\nu \cup \gamma) :=  \nu^{\otimes \pi_0(L)} \sharp_{\pi_0(L)} Z(L\cup \gamma) 
\quad \in \A(L\cup \gamma)
$$
which means that a copy of $\nu$ is summed along each connected component of $L$ in $Z(L\cup \gamma)$.
Also, we associate to the $(\pm 1)$-framed unknot the following quantity:
$$
U_\pm = \int \chi^{-1}\left( \nu \sharp Z(\circlearrowleft_{\pm 1})\right)
\quad \in \A(\varnothing).
$$
Note that $U_\pm$ is group-like (since, by Theorem \ref{th:JMM}, 
the formal integral of a non-degenerate Gaussian that is group-like, is group-like as well) 
and so is invertible.

\begin{definition}
Let $(B,\gamma)$ be a tangle in a homology cube. 
A \emph{surgery presentation} of $(B,\gamma)$ is a couple $(L,\gamma)$,
where $L$ is a framed oriented link in $[-1,1]^3$, $\gamma$ is a tangle in $[-1,1]^3$ disjoint from $L$
and  surgery along $L$ transforms $([-1,1]^3,\gamma)$ to $(B,\gamma)$.
\end{definition}

\begin{definition}
The \emph{Kontsevich--LMO invariant} of a $q$-tangle $\gamma$ in a homology cube $B$ is 
\begin{equation}
\label{eq:LMO}
Z(B,\gamma) := U_+^{-\sigma_+(L)} \sqcup U_-^{-\sigma_-(L)} \sqcup
\int_{\pi_0(L)} \chi^{-1}_{\pi_0(L)} Z(L^\nu \cup \gamma)  \quad \in \A(\gamma)
\end{equation}
where $(L,\gamma)$ is a surgery presentation of $(B,\gamma)$,
$(\sigma_+(L),\sigma_-(L))$ denotes the signature of the linking matrix of $L$, and
the action $\sqcup$ of  $\A(\varnothing)$ on $\A(\gamma)$ is given by the disjoint union operation.
\end{definition}

\noindent
The fact that $Z(B,\gamma)$ does not depend on the choice of the surgery presentation $(L,\gamma)$
of $(B,\gamma)$ follows from Kirby's theorem by adapting the arguments in \cite[\S 3 \& \S 5.1]{BGRT2}. 
See \cite{Moffatt} for a similar construction. 

One easily checks that, just as the Kontsevich integral, 
the Kontsevich--LMO invariant is functorial and tensor-preserving. 
By construction, $Z$ contains the Kontsevich integral and the LMO invariant as required.\\

In fact, we will only need to consider the Kontsevich--LMO invariant 
for bottom-top $q$-tangles in homology cubes. 
In this case, we can add the following statement:

\begin{lemma}
\label{lem:linking_matrix}
For all bottom-top $q$-tangle $\gamma$ in a homology cube $B$,
$$
\chi^{-1} Z(B,\gamma) \ \in \A(\pi_0(\gamma))
$$ 
is group-like and its $s$-reduction is $\left[\Lk_B(\gamma)/2\right]$.
\end{lemma}

\noindent
This group-like property of the Kontsevich--LMO invariant $Z(B,\gamma)$ 
is well-known when $B=[-1,1]^3$ or when $\gamma=\varnothing$.

\begin{proof}[Proof of Lemma \ref{lem:linking_matrix}]
Let $(L,\gamma)$ be a surgery presentation of $(B,\gamma)$.
Then, $Z(B,\gamma)$ is given by formula (\ref{eq:LMO}). 
The lemma is well-known to hold true when $B=[-1,1]^3$, so
$$
\chi^{-1}_{\pi_0(L\cup \gamma)} Z(L^\nu \cup \gamma)
= \left[\Lk_{[-1,1]^3}(L)/2 + \Lk_{[-1,1]^3}(\gamma)/2 + \Lk_{[-1,1]^3}(\gamma,L) \right] \sqcup [T]
$$
for a certain $T\in \A^{Y,c}(\pi_0(L\cup \gamma))$. Next, one integrates:  
\begin{eqnarray*}
&&\chi^{-1}_{\pi_0(\gamma)} \int_{\pi_0(L)} \chi^{-1}_{\pi_0(L)} Z(L^\nu \cup \gamma) \\
&=& \left\langle \left[ - \Lk_{[-1,1]^3}(L)^{-1} /2 \right],
\left[ \Lk_{[-1,1]^3}(\gamma)/2 + \Lk_{[-1,1]^3}(\gamma,L) \right] \sqcup [T] \right\rangle_{\pi_0(L)}\\
&=& \left[\Lk_{[-1,1]^3}(\gamma)/2 -
\Lk_{[-1,1]^3}(\gamma,L) \cdot \Lk_{[-1,1]^3}(L)^{-1} \cdot \Lk_{[-1,1]^3}(L,\gamma) /2\right] \sqcup [T']
\end{eqnarray*}
where the last identity follows from Theorem \ref{th:JMM} and involves a $T'\in \A^{Y,c}(\pi_0(\gamma))$.

\begin{claim}
\label{claim:lk}
Let $K$ be an oriented framed link in $[-1,1]^3$ whose linking matrix $\Lk_{[-1,1]^3}(K)$ is non-degenerate. 
For any two oriented knots $U$ and $V$ in $[-1,1]^3\setminus K$, one has 
$$
\Lk_{[-1,1]^3_K}(U,V) = \Lk_{[-1,1]^3}(U,V) -
\Lk_{[-1,1]^3}(U,K) \cdot {\Lk_{[-1,1]^3}(K)}^{-1} \cdot \Lk_{[-1,1]^3}(K,V).
$$ 
\end{claim}

This claim is easily proved using the homological definition of linking numbers. Thus, we obtain that
$$
\chi^{-1}_{\pi_0(\gamma)} \int_{\pi_0(L)} \chi^{-1}_{\pi_0(L)} 
Z(L^\nu \cup \gamma) = \left[\Lk_B(\gamma)/2 \right] \sqcup [T'].
$$
Since $U_\pm$ is a group-like element of $\A(\varnothing)$, 
we conclude that $\chi^{-1} Z(B,\gamma)$ is group-like of the form 
$\left[ \Lk_{B}(\gamma)/2 \right] \sqcup \left[T''\right]$
for a certain $T''\in \A^{Y,c}(\pi_0(\gamma))$.
\end{proof}

\begin{remark}
The Kontsevich--LMO invariant of $q$-tangles in $\Q$-homology cubes 
is defined exactly in the same way.
Lemma \ref{lem:linking_matrix} works in the rational case as well.
\end{remark}

\vspace{0.5cm}

\section{The functorial LMO invariant of Lagrangian cobordisms}

In this section, the Kontsevich--LMO invariant of bottom-top tangles in homology cubes is used to construct
an invariant of Lagrangian cobordisms.
After normalization, this invariant gives rise to a functor, which we call the \emph{LMO functor}.

\subsection{The category $\LqCob$ of Lagrangian $q$-cobordisms}

In this subsection, we define the domain of the LMO functor.

\begin{definition}
A \emph{Lagrangian $q$-cobordism} is a Lagrangian cobordism $(M,m)$ from $F_g$ to $F_f$ together 
with non-associative words $w_t(M)$ of length $g$ and $w_b(M)$ of length $f$ in the single letter $\bullet$.
\end{definition}

Given two Lagrangian $q$-cobordisms $M$ and $N$ such that $w_t(M)=w_b(N)$, 
one can form the new Lagrangian cobordism $M \circ N$
(by composition in $\Cob$) and equip it with the non-associative words 
$w_t(N)$ and $w_b(M)$. Thus, one obtains a category $\LqCob$ whose objects
are non-associative words in the single letter $\bullet$ 
and whose morphisms are Lagrangian $q$-cobordisms.
There is a tensor product $\otimes$ given by horizontal juxtaposition in the $x$-direction:
Thus, the category $\LqCob$ is monoidal (in the non-strict sense).

Similarly, we define the monoidal subcategory $\sLqCob$ of $\LqCob$
formed by \emph{special Lagrangian $q$-cobordisms}.

\begin{remark}
The category of \emph{$\Q$-Lagrangian $q$-cobordisms} is defined in the same way,
and is denoted by $\QLqCob$.
\end{remark}

\subsection{The category $\tsA$ of top-substantial Jacobi diagrams}

\label{subsec:top-substantial_diagrams}

In this subsection, we define the codomain of the LMO functor.

\begin{definition}
\label{def:top-substantial_diagram}
Let $f,g \geq 0$ be integers.
An element of $\A(\set{g}^+ \cup \set{f}^-)$ is \emph{top-substantial} if it is $\set{g}^+$-substantial.
\end{definition}

For all integers $f,g \geq 0$, we denote by 
$$
\tsA(g,f)
$$
the subspace of top-substantial elements of $\A(\set{g}^+ \cup \set{f}^-)$.
For all integers $f,g,h \geq 0$, we define a bilinear map
$$
\xymatrix{
{\tsA(g,f) \times \tsA(h,g)} \ar[rr]^-{\centerdot\ \circ\ \centerdot }& & {\tsA(h,f)}
}
$$
by the formula
$$
x\circ y := \left\langle\ (x/ i^+ \mapsto i^*)\ , 
\ (y / i^- \mapsto i^*)\  \right\rangle_{\set{g}^*}
$$
where $\set{g}^* = \{1^*, \dots, g^*\}$ is an extra set of variables,
$(y \left/ i^- \mapsto i^* \right.)$ denotes the Jacobi diagram obtained from $y$
by the change of variables $i^- \mapsto i^* $ for all $i=1,\dots,g$, 
and $(x\left/ i^+ \mapsto i^* \right.)$ has the similar meaning.
An equivalent formula for $\circ$ is
$$
x \circ y := 
\left(\begin{array}{c}
\hbox{sum of all ways of gluing the $i^+$-colored vertices of $x$}\\
\hbox{to the $i^-$-colored vertices of $y$, for all $i=1,\dots,g$}
\end{array}\right).
$$
It follows from the next lemma that
$$
\left\{
\begin{array}{ll}
\forall x \in \tsA(g,f), & \ x \circ \Id_g = x \\
\forall y \in \tsA(h,g), &  \Id_g \circ y = y
\end{array} \right.
\quad \hbox{ where } \quad
\Id_g := \left[ \sum_{i=1}^g \strutgraph{i^-}{i^+}\right].
$$
Thus, the above discussion defines a category $\tsA$.
The disjoint union operation of Jacobi diagrams
gives $\tsA$ the structure of a monoidal category (in the strict sense).

\begin{lemma}
\label{lem:exponential_shift}
Let $x \in \tsA(g,f)$ and let $y \in \tsA(h,g)$.
Then, for all $\set{h}^+ \times \set{g}^- $ matrix $D$, one has that
$$
x \circ \left( y \sqcup \left[ D \right] \right)  
= \left\langle\ \left( x / i^+ \mapsto i^* + D \cdot i^- \right)\ ,\ 
\left(y /i^-\mapsto i^*\right)\ \right\rangle_{\set{g}^*}
\ \in \tsA(h,f).
$$
Similarly, for all $\set{f}^- \times \set{g}^+$ matrix $C$, one has that
$$
\left( \left[C\right] \sqcup x \right) \circ y
= \left\langle\ \left( x / i^+ \mapsto i^* \right)\ ,\
\left(y /i^-\mapsto i^* + C \cdot i^+ \right) \ \right\rangle_{\set{g}^*}
\ \in \tsA(h,f).
$$
\end{lemma}

\noindent
In the above statement, the matrix $D$ is regarded as a linear map
$\Q \cdot \set{g}^- \to  \Q \cdot \set{h}^+$: 
Thus, $D \cdot i^-$ denotes $\sum_{j=1}^{h} d_{j^+,i^-} \cdot j^+$.
The matrix $C$ has the similar meaning.

\begin{proof}
We prove the first statement:
\begin{eqnarray*}
&& x \circ \left( y \sqcup \left[ D \right] \right) 
\\
&=& x \circ \left( y \sqcup \left[ \sum_{k=1}^g  \strutgraph{k^-}{D\cdot k^-}\quad \right] \right)
\\
&=& \left\langle (x / i^+\mapsto i^*) , 
(y/ i^- \mapsto i^*) \sqcup \bigsqcup_{k=1}^g
\left[ \strutgraph{k^*}{D\cdot k^-} \quad \right] \right\rangle_{\set{g}^*}
\\
&=& \sum_{n_1,\dots,n_g \geq 0} \frac{1}{n_1!\cdots n_g!} 
\left\langle (x / i^+\mapsto i^*) , 
(y/ i^- \mapsto i^*) \sqcup \bigsqcup_{k=1}^g
\left( \strutgraph{k^*}{D\cdot k^-}\quad\ \right)^{\sqcup n_k} \right\rangle_{\set{g}^*} 
\\
&=& \sum_{n_1,\dots,n_g \geq 0} \left\langle  
\left(\begin{array}{c} \hbox{\footnotesize sum of all ways of replacing $n_k$ times the color $k^*$ } \\
\hbox{\footnotesize by the color $D \cdot k^-$ in $(x/i^+ \mapsto i^*)$, for all $k=1,\dots,g$} 
\end{array}\right) ,
(y/ i^- \mapsto i^*) \right\rangle_{\set{g}^*}.
\end{eqnarray*}
We conclude since this is equal to 
$\left\langle\ \left( x / i^+ \mapsto i^* + D \cdot i^- \right)\ ,\ 
\left(y /i^-\mapsto i^*\right)\ \right\rangle_{\set{g}^*}$.
\end{proof}

Next lemma (which will be used later) describes how the composition law $a \circ b$ of $\tsA$ decomposes 
into an ``$s$-part'' and a ``$Y$-part'' if $a$ and $b$ can themselves be decomposed that way.

\begin{lemma}
\label{lem:composition_ts_diagrams}
Let $a \in \tsA(g,f)$ and  $b \in \tsA(h,g)$.
Assume that they can be decomposed as 
$$
a = \left[ A/2 \right] \sqcup a^Y
\quad \hbox{and} \quad
b = \left[ B/2 \right] \sqcup b^Y
$$
where $A$ is a symmetric 
$\left(\set{g}^+ \cup \set{f}^-\right) \times \left(\set{g}^+ \cup \set{f}^-\right)$ matrix 
and where $B$ is a symmetric 
$\left(\set{h}^+ \cup \set{g}^-\right) \times \left(\set{h}^+ \cup \set{g}^-\right)$ 
matrix of the form
\begin{equation}
\label{eq:matrices}
A =  \left( \begin{array}{cc}
0 & A^{+-} \\ A^{-+} & A^{--}
\end{array} \right)
\quad \hbox{and} \quad B = \left( \begin{array}{cc}
0 & B^{+-} \\ B^{-+} & B^{--}
\end{array} \right).
\end{equation}
Then, $a \circ b$ is also decomposable as
$$
a \circ b = \left[ \frac{1}{2} \left( \begin{array}{cc}
0  &  B^{+-} A^{+-} \\ A^{-+} B^{-+} & A^{--} + A^{-+} B^{--} A^{+-} 
\end{array} \right) \right] \sqcup \left( a^Y \Star{A}{B}\ b^Y\right)
$$
where $a^Y \Star{A}{B}\ b^Y$ belongs to $\AY(\set{h}^+ \cup \set{f}^-)$ and is defined below.
Moreover, if $a$ and $b$ are group-like, then $a\circ b$ is group-like as well. 
\end{lemma}

\noindent
To complete the previous statement, we associate to all pair of matrices $(A,B)$ 
of the form (\ref{eq:matrices}) a bilinear pairing 
$$
\xymatrix{
{\AY(\set{g}^+ \cup \set{f}^-) \times \AY(\set{h}^+\cup\set{g}^-)}
\ar[rr]^-{\centerdot \Star{A}{B}\ \centerdot} && {\AY(\set{h}^+\cup\set{f}^-)}
}
$$
defined by the formula
$$
x \Star{A}{B}\  y :=
\left\langle \left(x/ i^+ \mapsto i^* + B^{+-}\cdot i^- + A^{-+} B^{--} \cdot i^-\right)\ ,\ 
\begin{array}{r}
\left([B^{--}/2]/ i^- \mapsto i^*\right) \sqcup \\ 
\left(y / i^- \mapsto i^* + A^{-+} \cdot i^+ \right) 
\end{array}
\right\rangle_{\set{g}^*}.
$$

\begin{example}
\label{ex:star}
Consider the special case when $f=g=h$ and 
$$
A = B = \left( \begin{array}{cc}
0 & I_g^{+-} \\ I_g^{-+} & 0
\end{array} \right)
$$
where $I_g^{+-}$ denotes the ``identity'' matrix 
$\left( \delta_{i,j} \right)_{i^+\in \set{g}^+, j^-\in \set{g}^-}$
and $I_g^{-+}$ is its transpose. Then, the above product is simply denoted by $\star$ and the formula is 
$$
x \star y :=
\left\langle\ (x / i^+ \mapsto i^* + i^+)\ ,\ (y / i^- \mapsto i^* + i^-)\ \right\rangle_{\set{g}^*}.
$$
\end{example}

\begin{proof}[Proof of Lemma \ref{lem:composition_ts_diagrams}]
The last statement is an application of Theorem \ref{th:JMM}. 
The first statement is proved using Lemma \ref{lem:exponential_shift} as follows: 
\begin{eqnarray*}
&& a\circ b
\\
&=&\left\langle \begin{array}{l} 
\left([A^{--}/2]/ i^+\mapsto i^* \right) \sqcup \\
\left([A^{+-}]/i^+ \mapsto i^* \right) \sqcup \left( a^Y /i^+ \mapsto i^* \right) \end{array},
\begin{array}{r} 
\left([B^{--}/2]/i^- \mapsto i^* \right) \sqcup \\
\left( b^Y /i^- \mapsto i^* \right) \sqcup \left([B^{+-}]/i^- \mapsto i^* \right)\end{array}
\right\rangle_{\set{g}^*} 
\\
&=& [A^{--}/2] \sqcup \left\langle \begin{array}{r} 
 \left([A^{+-}]/i^+ \mapsto i^* + B^{+-}\cdot i^- \right) \\
\sqcup \left( a^Y /i^+ \mapsto i^*  + B^{+-}\cdot i^-  \right) \end{array},
\begin{array}{l}
\left([B^{--}/2]/i^- \mapsto i^* \right) \\
\sqcup  \left( b^Y /i^- \mapsto i^* \right) \end{array}
\right\rangle_{\set{g}^*}
\\
&=& [A^{--}/2] \sqcup [B^{+-} A^{+-} ] \sqcup \left\langle \begin{array}{r} 
 \left([A^{+-}]/i^+ \mapsto i^* \right) \\
\sqcup \left( a^Y /i^+ \mapsto i^*  + B^{+-}\cdot i^-  \right) \end{array},
\begin{array}{l}
\left([B^{--}/2]/i^- \mapsto i^* \right) \\
\sqcup  \left( b^Y /i^- \mapsto i^* \right) \end{array}
\right\rangle_{\set{g}^*}
\end{eqnarray*}
\begin{eqnarray*}
&=& [A^{--}/2] \sqcup [B^{+-} A^{+-} ] \sqcup 
\\
&& \left\langle \left(a^Y /i^+ \mapsto i^*  + B^{+-}\cdot i^-  \right),
\begin{array}{l}
\left([B^{--}/2]/i^- \mapsto i^* + A^{-+}\cdot i^+\right) \\
\sqcup  \left( b^Y /i^- \mapsto i^* + A^{-+}\cdot i^+ \right) \end{array}
\right\rangle_{\set{g}^*}
\\
&=& [A^{--}/2] \sqcup [B^{+-} A^{+-} ] \sqcup [A^{-+} B^{--} A^{+-}/2]  \sqcup
\\
&&\left\langle \left( a^Y / i^+ \mapsto (i^* + A^{-+}B^{--} \cdot i^-) + B^{+-}\cdot i^-  \right) ,
\begin{array}{l} \left([B^{--}/2]/i^- \mapsto i^*\right) \sqcup \\ 
\left( b^Y /i^- \mapsto i^* + A^{-+}\cdot i^+ \right) \end{array}
\right\rangle_{\set{g}^*}
\\
&=&
\left[ \frac{1}{2} \left( \begin{array}{cc}
0  &  B^{+-} A^{+-} \\ A^{-+} B^{-+} & A^{--} + A^{-+} B^{--} A^{+-} 
\end{array} \right) \right] \sqcup
\\
&& \left\langle \left(a^Y/ i^+ \mapsto i^* + B^{+-}\cdot i^- + A^{-+} B^{--} \cdot i^-\right), 
\begin{array}{l} \left([B^{--}/2]/i^- \mapsto i^*\right) \sqcup \\ 
\left( b^Y /i^- \mapsto i^* + A^{-+}\cdot i^+ \right) \end{array}
\right\rangle_{\set{g}^*}.
\end{eqnarray*}
\end{proof}

\subsection{The unnormalized LMO invariant $Z$}

Each Lagrangian cobordism corresponds to a unique bottom-top tangle in a homology cube (Lemma \ref{lem:bh_bt}). 
Thus, we merely define the LMO invariant of the former to be the Kontsevich--LMO invariant of the latter.
Taking into account parenthesizings, this gives the following

\begin{definition}
\label{def:unnormalized_LMO}
Let $M$ be a Lagrangian $q$-cobordism from $F_g$ to $F_f$.
The \emph{unnormalized LMO invariant}  of $M$ is 
$$
Z(M) := Z(B,\gamma) \ \in \A(\gamma) = \A(\gamma^+ \cup \gamma^-)  = 
\A\left( \cupright^{\set{g}} \capleft^{\set{f}} \right)
$$
where $(B,\gamma)$ is the bottom-top tangle presentation of $M$. 
More precisely, $\gamma$ is equipped here with the non-associative words 
$$
w_t(\gamma) := \left(w_t(M)/\bullet \mapsto (+-)\right) \quad \hbox{and} \quad
w_b(\gamma) := \left(w_b(M)/\bullet \mapsto (+-)\right)
$$
and the connected components of $\gamma^+$ and  $\gamma^-$
are numbered increasingly along the $x$ direction, 
from $1$ to $g$ and from $1$ to $f$ respectively.
\end{definition} 

\noindent
We will work mainly with the symmetrized version of $Z(M)$, namely
$$
\chi^{-1} Z(M) \in \A(\set{g}^+ \cup \set{f}^-).
$$
It follows from Lemma \ref{lem:linking_matrix} and Lemma \ref{lem:bh_bt} 
that $\chi^{-1} Z(M)$ is top-substantial. 

To sum up, we have obtained so far a family of maps
$$
\left( \LqCob(w,v) \longrightarrow \tsA(|w|,|v|), \ M \longmapsto \chi^{-1} Z(M) \right)_{w,v}
$$
where $v$ and $w$ range over non-associative words in the single letter $\bullet$.
Those maps are easily seen to preserve the tensor product,
but, the next subsection reveals that they do \emph{not} define a functor.

\subsection{Normalization of the LMO invariant}

\label{subsec:normalization}

Let us now see how the unnormalized LMO invariant $Z$ of Lagrangian cobordisms behaves 
with respect to composition. 

First of all, we fix some notations. 
For all formal variables $x,y,r$, set 
$$
\lambda(x,y;r) := \chi^{-1}\left( \begin{array}{c} {\relabelbox \small
\epsfxsize 1.2truein \epsfbox{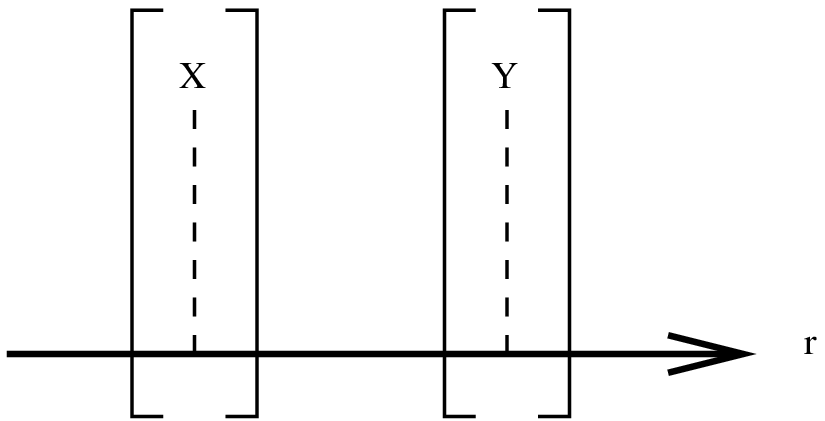}
\adjustrelabel <-0.1cm,-0.0cm> {X}{$x$}
\adjustrelabel <-0.1cm,-0.0cm> {Y}{$y$}
\adjustrelabel <-0.0cm,-0.0cm> {r}{$r$}
\endrelabelbox} \end{array}
\right) \ \in \A(\{x,y,r\})
$$
where the brackets denote the exponential map in $\A(\longrightarrow^r,\{x,y\})$
with respect to the natural multiplication.

\begin{remark}
As observed in \cite[Prop. 5.4]{BGRT2}, this formal series of Jacobi diagrams 
can be computed from the Baker--Campbell--Hausdorff series.
Indeed, the BCH series  
$$
\log\left(\exp(x)\cdot \exp(y)\right) \in  \Q[[x,y]]
$$
(where the variables $x$ and $y$ do not commute) belongs to the completed Lie $\Q$-algebra $\hbox{Lie}(x,y)$
freely generated by $x$ and $y$. Recall that $\hbox{Lie}(x,y)$ embeds into $\A^c(\{x,y,r\})$ 
by writing Lie commutators as $r$-rooted binary trees whose leaves are colored by $x$ and $y$, for example
$$
[x,[[x,y],y]] \longmapsto
\begin{array}{c}
{\relabelbox \small
\epsfxsize 0.8truein \epsfbox{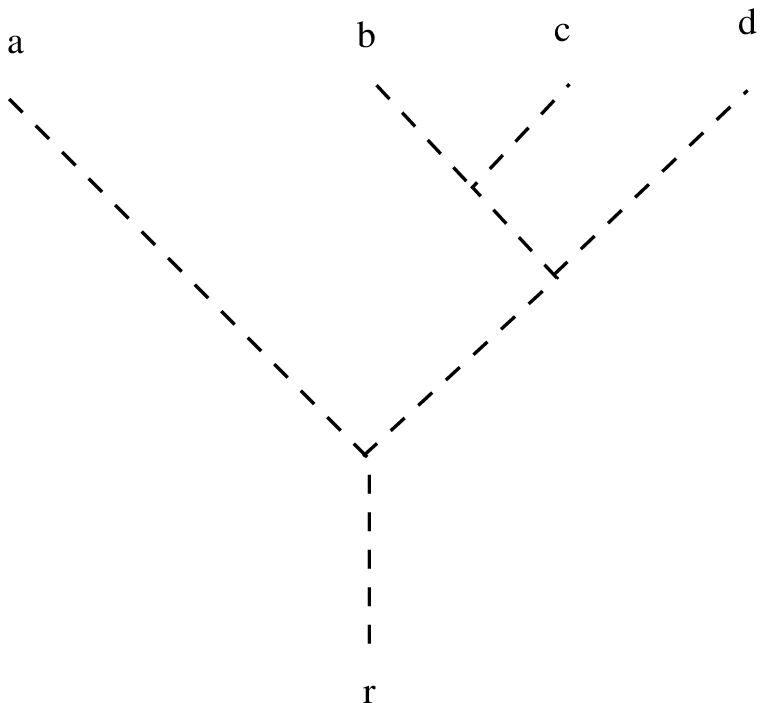}
\adjustrelabel <-0.05cm,-0.0cm> {a}{$x$}
\adjustrelabel <-0.cm,-0.0cm> {b}{$x$}
\adjustrelabel <-0.0cm,-0.0cm> {c}{$y$}
\adjustrelabel <-0.0cm,-0.0cm> {d}{$y$}
\adjustrelabel <-0.05cm,-0.1cm> {r}{$r$}
\endrelabelbox}
\end{array}.
$$
Thus, the BCH series defines a formal series of connected tree diagrams
$\Lambda(x,y;r) \in \mathcal{A}^c(\{x,y,r\})$.
It is easily seen that $\lambda(x,y;r) = \left[ \Lambda(x,y;r) \right]$.
\end{remark} 

Next, we  define an element of $\A\left(\{x_+,x_-\}\right)$ by
$$
\mathsf{T}(x_+,x_-) :=  U_+^{-1} \sqcup U_-^{-1} \sqcup \int_{r}  
\left\langle \lambda\left(x_-,y_-;r_-\right) \sqcup \lambda\left(x_+,y_+;r_+ \right) , 
\chi^{-1} Z\left(T_1^\nu \right) \right\rangle_{y}.
$$
In that formula,  $T_1$ denotes the bottom-top tangle of type $(1,1)$ shown on Figure \ref{fig:bt_tangle_T},
whose top and bottom components are labeled by $y_+$ and $y_-$ respectively.

\begin{lemma}
\label{lem:T_g}
$\mathsf{T}(x_+,x_-)$ is a group-like element of $\A\left(\{x_+,x_-\}\right)$  
with $s$-reduction $\left[\strutgraph{x_-}{x_+}\right]$.
\end{lemma}

\begin{proof}
The Kontsevich integral of a $q$-tangle in $[-1,1]^3$ is group-like 
and the series $\lambda(x,y;r)$ is clearly group-like. So, by Theorem \ref{th:JMM},
the integrand in the formula defining $\mathsf{T}(x_+,x_-)$ is group-like.
Since formal Gaussian integration transforms a group-like element to a group-like element 
(again, by  Theorem \ref{th:JMM}), we conclude that $\mathsf{T}(x_+,x_-)$ is group-like.
Furthermore, that integrand is equal to
\begin{eqnarray*}
&&\left\langle \lambda\left(x_-,y_-;r_-\right) \sqcup \lambda\left(x_+,y_+;r_+ \right) 
, \chi^{-1} Z\left(T_1^\nu \right) \right\rangle_{y} \\
&=& \left\langle \left[\strutgraph{r_-}{x_-} + \strutgraph{r_-}{y_-} +
\strutgraph{r_+}{x_+} + \strutgraph{r_+}{y_+} \right] \sqcup \left(\hbox{something in }\AY \right)
, \left[ - \strutgraph{y_-}{y_+} \right] \sqcup \left(\hbox{something in }\AY \right) \right\rangle_{y}\\
&=& \left[ - \strutgraph{r_-}{r_+} + \strutgraph{r_-}{x_-} + \strutgraph{r_+}{x_+}  \right] 
\sqcup \left(\hbox{something in }\AY \right).
\end{eqnarray*}
Thus, after formal Gaussian integration, one gets
$$
\mathsf{T}(x_+,x_-) = \left[ \strutgraph{x_-}{x_+} \right] \sqcup \left(\hbox{something in }\AY \right).
$$
\end{proof}

Finally, for all integer $g\geq 0$, we set
$$
\mathsf{T}_g := \mathsf{T}(1^+,1^-) \sqcup \cdots 
\sqcup \mathsf{T}(g^+,g^-)\ \in \A(\set{g}^+ \cup \set{g}^-).
$$
By the previous lemma, $\mathsf{T}_g$
is a group-like element of $\tsA(g,g)$ and its $s$-reduction is $\Id_g$. 

\begin{lemma}
\label{lem:functoriality}
Let $w$ be a non-associative word  of length $g$ in the single letter $\bullet$, 
and let $M$ and $N$ be two Lagrangian $q$-cobordisms such that $w_t(M)=w_b(N)=w$. 
Then, we have 
$$
\chi^{-1}Z(M \circ N) = 
\chi^{-1} Z(M) \circ \mathsf{T}_g \circ \chi^{-1} Z(N).
$$
\end{lemma}

\begin{proof}
Let $(B,\gamma)$ and $(C,\upsilon)$ be the bottom-top $q$-tangles corresponding to 
$M$ and $N$ respectively. Let also $(K,\gamma)$ 
and $(L,\upsilon)$ be surgery presentations of $(B,\gamma)$ and $(C,\upsilon)$ respectively.
We denote by $T$ the $2g$-component oriented framed link in $[-1,1]^3$
obtained by gluing the bottom-top tangle $T_g$ from
Figure \ref{fig:bt_tangle_T} ``between'' $\gamma^+$ and $\upsilon^-$.
Then, $(K \cup T \cup L,  \gamma^- \cup \upsilon^+)$ is a surgery presentation of 
$(B,\gamma) \circ (C,\upsilon)$ so that
$$
Z(M \circ N) = 
\frac{\displaystyle
\int_{\pi_0(K \cup T \cup L)} \chi^{-1}_{\pi_0(K \cup T \cup L)} 
Z((K \cup T \cup L)^\nu \cup (\gamma^- \cup \upsilon^+))}
{U_+^{\sigma_+(K\cup T \cup L)} \sqcup U_-^{\sigma_-(K\cup T \cup L)}}.
$$ 
By the functoriality of the Kontsevich integral 
$Z$ at the level of $q$-tangles in $[-1,1]^3$, we can write 
$$
Z\left((K \cup T \cup L)^\nu \cup (\gamma^- \cup \upsilon^+)\right)
= Z(K^\nu \cup \gamma) \circ Z(T_g^\nu) \circ Z(L^\nu\cup \upsilon)
$$
where $\circ$ denotes the composition in the category $\A$. This implies that
$$
\chi^{-1}_{\pi_0(K \cup T \cup L)} 
\left( Z((K \cup T \cup L)^\nu \cup (\gamma^- \cup \upsilon^+))\right)=
$$
$$  
\chi^{-1}_{\pi_0(T)} \left(
\chi_{\pi_{0}(K)}^{-1}\left(Z(K^\nu\cup \gamma)\right)
\circ  Z(T_g^\nu) \circ 
 \chi_{\pi_{0}(L)}^{-1}\left(Z(L^\nu\cup \upsilon)\right) \right).
$$
Since the matrices $\Lk_{[-1,1]^3}(K)$ and $\Lk_{[-1,1]^3}(L)$ are invertible, we can integrate by iteration 
(see \cite[Prop. 2.11]{BGRT2}) along $\pi_0(K)$, next along $\pi_0(L)$ and finally along $\pi_0(T)$.
Moreover, it is proved below that 
\begin{equation}
\label{eq:signature}
\sigma_\pm(K \cup T \cup L) = \sigma_\pm(K) + g + \sigma_{\pm}(L).
\end{equation}
Thus, we obtain that
$$
Z(M\circ N) = U_+^{-g} \sqcup U_-^{-g} \sqcup 
\int_{\pi_0(T)} \chi^{-1}_{\pi_0(T)}\left(Z(B,\gamma) 
\circ Z(T_g^\nu) \circ Z(C,\upsilon)\right)
$$
or, equivalently, that
$$
\chi^{-1} Z(M \circ N) =
U_+^{-g} \sqcup U_-^{-g} \sqcup 
\int_{\pi_0(T)}\! \chi^{-1}_{\pi_0(T)}\left(
\chi^{-1}_{\pi_0(\gamma^-)} Z(B,\gamma) 
\circ Z(T_g^\nu) \circ 
\chi^{-1}_{\pi_0(\upsilon^+)} Z(C,\upsilon)\right).
$$
Assume that $M$ is a cobordism from $F_g$ to $F_f$ and that $N$ is from $F_h$ to $F_g$.
We number connected components of $1$-manifolds as follows:
$$
\pi_0(\gamma^-) = \set{f}^- \ , \ \pi_0(\gamma^+) = \set{g}^\cup \quad \quad \quad 
\pi_0(\upsilon^-) = \set{g}^\cap \ , \ \pi_0(\upsilon^+) = \set{h}^+
$$
$$
\pi_0(T_g^-) = \set{g}^\vartriangle \ , \ \pi_0(T_g^+) = \set{g}^\triangledown \quad \quad \quad
\pi_0(T^-) = \set{g}^\bot  \ , \ \pi_0(T^+) = \set{g}^\top.
$$
The series $\lambda(x,y;r)$ is designed so that
\begin{eqnarray*}
&&\chi^{-1}_{\pi_0(T)}\left(
\chi^{-1}_{\pi_0(\gamma^-)} Z(B,\gamma) 
\circ Z(T_g^\nu) \circ \chi^{-1}_{\pi_0(\upsilon^+)} Z(C,\upsilon)\right)\\
&=&
\chi^{-1}_{\pi_0(T)}\left(
\chi_{\pi_0(\gamma^+)} \chi^{-1}_{\pi_0(\gamma)} Z(B,\gamma) 
\circ \chi_{\pi_{0}(T_g)} \chi_{\pi_{0}(T_g)}^{-1} Z(T_g^\nu) \circ 
\chi_{\pi_0(\upsilon^-)} \chi^{-1}_{\pi_0(\upsilon)} Z(C,\upsilon)\right)\\
&=&
\left\langle 
\bigsqcup_{i=1}^g \lambda(i^\vartriangle,i^\cup;i^\bot) 
\sqcup \bigsqcup_{i=1}^g \lambda(i^\cap,i^\triangledown; i^\top)  , 
\chi^{-1} Z(B,\gamma) \sqcup \chi^{-1} Z(C,\upsilon) \sqcup \chi^{-1} Z(T_g^\nu) 
\right\rangle_{\! \! \cap \triangledown \vartriangle \cup}\\
&=&
\left\langle 
\left\langle  \bigsqcup_{i=1}^g \lambda(i^\vartriangle,i^\cup;i^\bot) 
\sqcup \bigsqcup_{i=1}^g \lambda(i^\cap,i^\triangledown; i^\top)  , 
\chi^{-1} Z(T_g^\nu) \right\rangle_{\! \! \triangledown \vartriangle} ,
\chi^{-1} Z(B,\gamma) \sqcup \chi^{-1} Z(C,\upsilon) 
\right\rangle_{\! \! \cap\cup}.
\end{eqnarray*}
We deduce that
\begin{eqnarray*}
\chi^{-1} Z(M \circ N) &=& U_+^{-g} \sqcup U_-^{-g} \sqcup 
\int_{\top \bot}  \left\langle  \left\langle \cdots , \cdots \right\rangle_{\triangledown \vartriangle} ,
\chi^{-1} Z(M) \sqcup \chi^{-1} Z(N) \right\rangle_{\cap\cup} \\
& = & U_+^{-g} \sqcup U_-^{-g} \sqcup 
\left\langle \int_{\top \bot} \left\langle \cdots , \cdots \right\rangle_{\triangledown \vartriangle},
 \chi^{-1} Z(M) \sqcup \chi^{-1} Z(N) \right\rangle_{\cap \cup}\\
& = & \left\langle  \left(\mathsf{T}_w \left/ 
\begin{array}{c} i^- \mapsto i^\cup \\ i^+\mapsto i^\cap \end{array}\right.\right),
\chi^{-1}Z(M) \sqcup \chi^{-1}Z(N) \right\rangle_{\cap\cup}
\end{eqnarray*}
which involves the following element of $\A(\set{g}^+ \cup \set{g}^-)$:
$$
\mathsf{T}_w :=  U_+^{-g} \sqcup U_-^{-g} \sqcup
\int_{\top\! \bot} \! \!  \left\langle 
\bigsqcup_{i=1}^g \lambda\left(i^\vartriangle,i^-;i^\bot\right) 
\sqcup \bigsqcup_{i=1}^g \lambda\left(i^+,i^\triangledown;i^\top \right)
, 
\chi^{-1} Z\left(T_g^\nu \right) 
\right\rangle_{\! \! \! \triangledown \vartriangle}.
$$
Here,  the bottom-top tangle $T_g$ is equipped at the top and the bottom 
with the non-associative word obtained from $w$ by the rule ``$\bullet \mapsto (+-)$''. 
Since $T_g$ is the tensor product $g$ times of $T_1$, one sees that
$$
\forall \hbox{ word } w, \quad \chi^{-1} Z(T_g^\nu) = 
\underbrace{\chi^{-1} Z(T_1^\nu) \otimes \cdots \otimes \chi^{-1} Z(T_1^\nu)}_{g\ \hbox{\footnotesize times}} 
\ \in \A(\set{g}^\triangledown \cup \set{g}^\vartriangle).
$$
Thus, we conclude that $\mathsf{T}_w = \mathsf{T}_g$ so that 
$$
\chi^{-1} Z(M \circ N) = \chi^{-1} Z(M) \circ \mathsf{T}_g \circ \chi^{-1} Z(N).
$$

It now remains to prove identity (\ref{eq:signature}). 
The linking matrix of $L \cup T \cup K$ in $[-1,1]^3$ can be decomposed as
$$
\Lk(L \cup T \cup K) = \left(\begin{array}{c|c|c|c}
\Lk(L) & \Lk(L,\upsilon^-) & 0 & 0 \\
\hline \Lk(\upsilon^-,L) & \Lk(\upsilon^-) & - I_g & 0 \\
\hline 0 & -I_g & \Lk(\gamma^+) & \Lk(\gamma^+,K) \\
\hline 0 & 0 &  \Lk(K,\gamma^+) & \Lk(K)
\end{array}\right).
$$
Let $P$ be the non-degenerate matrix
$$
P:= \left(\begin{array}{c|c|c|c}
I_l & 0 & 0 & 0\\
\hline -\Lk(\upsilon^-,L) \Lk(L)^{-1} & I_g & 0 & 0\\
\hline 0 & 0 & I_g & -  \Lk(\gamma^+,K) \Lk(K)^{-1} \\
\hline 0 & 0 & 0 & I_k
\end{array}\right)
$$
where $l$ and $k$ are the number of connected components of $L$ and $K$ respectively.
The congruence $P \cdot \Lk(L \cup T \cup K) \cdot P^t$ gives
$$
\left(\begin{array}{c|c|c|c}
\Lk(L) & 0 & 0 & 0 \\
\hline 0 & \begin{array}{c} \Lk(\upsilon^-) -\\ \Lk(\upsilon^-,L) \Lk(L)^{-1} \Lk(L,\upsilon^-) \end{array}& -I_g & 0\\
\hline 0 & -I_g & \begin{array}{c} \Lk(\gamma^+) - \\ \Lk(\gamma^+,K) \Lk(K)^{-1} \Lk(K,\gamma^+) \end{array}& 0 \\
\hline 0 & 0 &  0 & \Lk(K)
\end{array}\right).
$$
Using Claim  \ref{claim:lk} and the fact that $\Lk_B(\gamma^+)=0$ (by Lemma \ref{lem:bh_bt}), 
we obtain that $\Lk(L \cup T \cup K)$ is congruent to
$$
\left(\begin{array}{c|c|c|c}
\Lk(L) & 0 & 0 & 0 \\
\hline 0 & \Lk_C(\upsilon^-) & -I_g & 0\\
\hline 0 & -I_g & 0 & 0 \\
\hline 0 & 0 &  0 & \Lk(K)
\end{array}\right)
$$
from which we deduce identity (\ref{eq:signature}).
\end{proof}

Lemma \ref{lem:functoriality} suggests the following normalization of the LMO invariant: 

\begin{definition}
\label{def:LMO}
The \emph{normalized LMO invariant} of a Lagrangian $q$-cobordism $M$ from $F_g$ to $F_f$ is 
$$
\Ztilde(M) :=  \chi^{-1} Z(M) \circ \mathsf{T}_{g}  \ \in \A\left(\set{g}^+ \cup \set{f}^-\right)
$$
where $Z(M)$ is the unnormalized LMO invariant from Definition \ref{def:unnormalized_LMO}.
\end{definition}

\noindent
According to the next lemma, $\Ztilde(M)$ splits as
$$
\Ztilde(M) = \left[\Lk(M)/2\right] \sqcup \ZtildeY(M)
$$
where $\ZtildeY(M) \in \AY\left(\set{g}^+ \cup \set{f}^-\right)$ 
denotes the $Y$-reduction of $\Ztilde(M)$.

\begin{lemma}
\label{lem:linking_matrix_L_cobordisms}
For all Lagrangian $q$-cobordism $M$ from $F_g$ to $F_f$, 
$\Ztilde(M)$ is group-like and its $s$-reduction is $\left[\Lk(M)/2\right]$.
\end{lemma}

\begin{proof}
Let $(B,\gamma)$ be the bottom-top $q$-tangle in a homology cube corresponding to $M$. 
Then, the definition of $\Ztilde(M)$ writes
$$
\Ztilde(M) = \chi^{-1} Z(B,\gamma) \circ \mathsf{T}_{g}.
$$
Since $\chi^{-1} Z(B,\gamma)$ and $\mathsf{T}_g$ are both group-like 
(by Lemma \ref{lem:linking_matrix} and Lemma \ref{lem:T_g} respectively), 
we conclude thanks to Lemma \ref{lem:composition_ts_diagrams}.
\end{proof}

\subsection{The LMO functor $\Ztilde$}

\label{subsec:functor}

We can now state and prove the main result of this section:

\begin{theorem}
\label{th:LMO_functor}
The normalized LMO invariant defines a tensor-preserving functor
$$
\Ztilde: \LqCob \longrightarrow \tsA
$$
from the category of Lagrangian $q$-cobordisms to the category of top-substantial Jacobi diagrams.
\end{theorem}

\begin{proof}
By Lemma \ref{lem:functoriality}, $\Ztilde$ preserves the composition law and, just like $\chi^{-1} Z$, 
it respects the tensor product as well.
It remains to check that, for all non-associative word $w$ of length $g$,
$\Ztilde: \LqCob(w,w) \to \tsA(g,g)$ sends $\Id_w$ to $\Id_g$.

We know from Lemma \ref{lem:functoriality} that 
$\Ztilde(\Id_w)\circ \Ztilde(\Id_w) = \Ztilde(\Id_w)$.
Let  $\star$ be the product defined in Example \ref{ex:star}:
Lemma \ref{lem:composition_ts_diagrams} 
implies that $\ZtildeY(\Id_w) \star \ZtildeY(\Id_w) = \ZtildeY(\Id_w)$.
Since $\ZtildeY(\Id_w)$ is group-like (by Lemma \ref{lem:linking_matrix_L_cobordisms}),
it can be written as
$$
\ZtildeY(\Id_w) = \varnothing + T + (\ideg>k)
$$
where $k>0$ and $T$ has i-degree k. Then, we must have $2\cdot T =T$ i.e$.$ $T=0$,
so that $\ZtildeY(\Id_w)=\varnothing$. 
\end{proof}

\begin{remark}
\label{rem:rational_case}
More generally, we obtain a tensor-preserving functor
$$
\Ztilde : \QLqCob \longrightarrow \tsA
$$ 
since the arguments in the last two subsections work with rational coefficients as well.
\end{remark}

\vspace{0.5cm}

\section{Computation of the LMO functor by pieces}

\label{sec:computation}

In order to compute the LMO functor $\Ztilde$ on a Lagrangian $q$-cobordism $M$, 
it is enough to decompose $M$ into ``elementary pieces'' -- with respect to the
composition law $\circ$ and the tensor product $\otimes$ of the category $\LqCob$ --
and to know the values of $\Ztilde$ on those pieces. 
In this section, we describe this approach.

\subsection{Generators of $\LqCob$}

\begin{figure}[h]
\centerline{\relabelbox \small
\epsfxsize 5truein \epsfbox{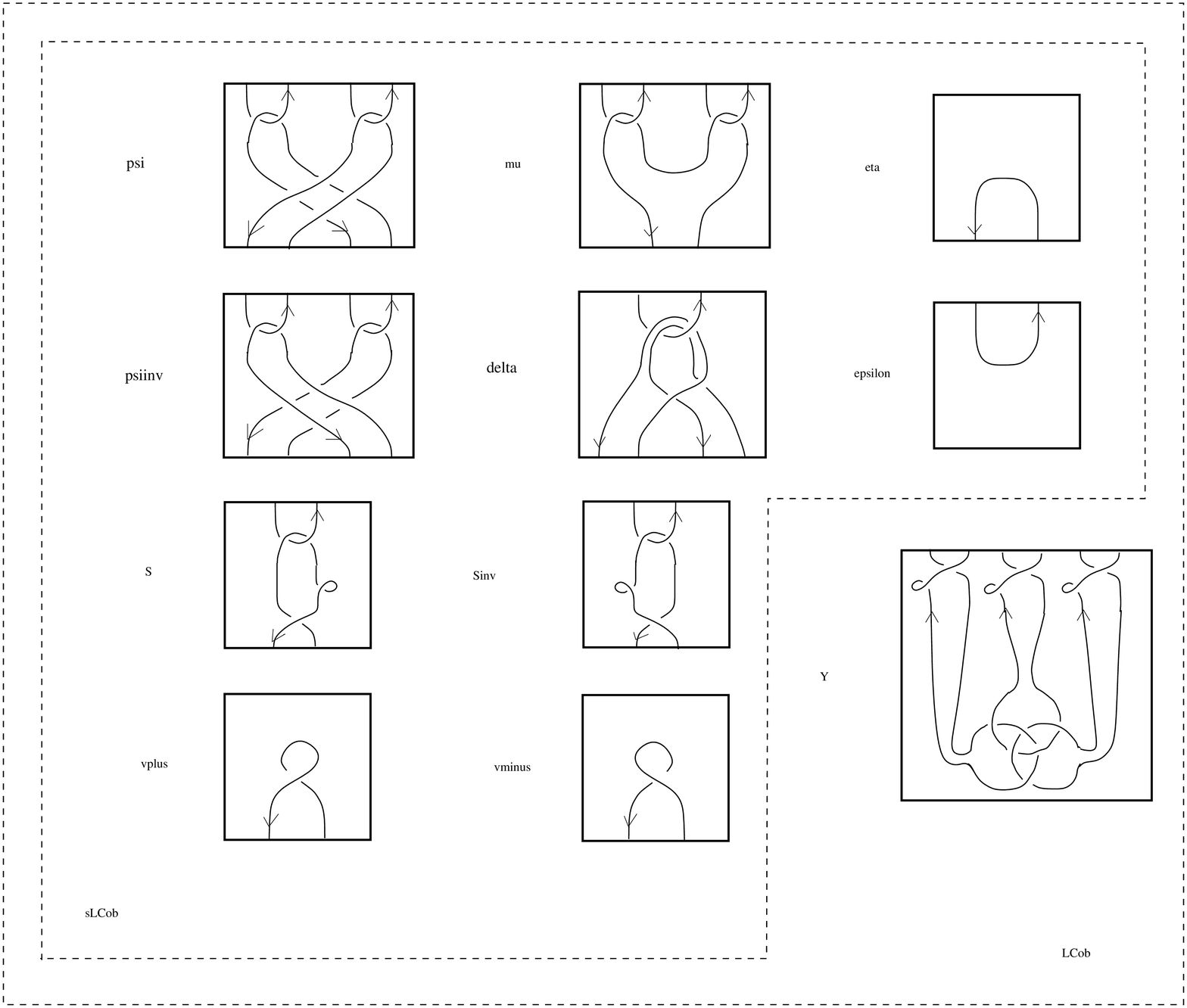}
\adjustrelabel <0cm,-0cm> {psi}{$\psi_{1,1}:=$}
\adjustrelabel <0cm,-0cm> {psiinv}{$\psi_{1,1}^{-1}:=$}
\adjustrelabel <0cm,-0cm> {eta}{$\eta:=$}
\adjustrelabel <0cm,-0cm> {mu}{$\mu:=$}
\adjustrelabel <0cm,-0cm> {epsilon}{$\epsilon:=$}
\adjustrelabel <0cm,-0cm> {delta}{$\Delta:=$}
\adjustrelabel <0cm,-0cm> {S}{$S:=$}
\adjustrelabel <0cm,-0cm> {Sinv}{$S^{-1}:=$}
\adjustrelabel <0cm,-0cm> {vplus}{$v_+:=$}
\adjustrelabel <0cm,-0cm> {vminus}{$v_-:=$}
\adjustrelabel <0cm,-0cm> {Y}{$Y:=$}
\adjustrelabel <0cm,-0cm> {LCob}{$\LCob$}
\adjustrelabel <0cm,-0cm> {sLCob}{$\sLCob$}
\endrelabelbox}
\caption{Generators of the monoidal categories $\sLCob$ and $\LCob$.}
\label{fig:generators}
\end{figure}

We indicate a system of generators for the monoidal category $\LqCob$.
For this, we recall from \cite[\S 14.5]{Habiro_BT}
that the monoidal category $\sLCob$ is generated by the morphisms
\begin{equation}
\label{eq:generators_sLCob}
\left(\psi_{1,1}^{\pm 1}, \mu,\eta,\Delta,\epsilon,S^{\pm 1}, v_\pm\right)
\end{equation}
shown on Figure \ref{fig:generators} in their bottom-top tangle presentations. 
For instance, observe that $\eta = C^0_1$ and that $\epsilon = C^1_0$.
Those generators of the monoidal subcategory $\sLCob$ of $\LCob \leq \Cob$ have the following categorical interpretation:
\begin{itemize}
\item The braiding $\psi_{1,1}$ extends in a unique way to braidings
  $\psi_{p,q}: p\otimes q \rightarrow q\otimes p$, defined for all $p,q\ge 0$, which give a braided category
  structure for $\Cob$, and hence for $\sLCob$ and $\LCob$.
\item $H:=(1,\mu ,\eta ,\Delta ,\epsilon ,S^{\pm 1})$ is a braided Hopf algebra with
  invertible antipode, as was first observed by Crane and Yetter \cite{CY}
  and Kerler \cite{Kerler} in the category $\Cob$.
\item The morphisms $v_\pm $ are ``ribbon elements'' of $H$ in the sense
  of Kerler \cite{Kerler_towards}.
\end{itemize}

Let also $Y: 3 \to 0$ be the Lagrangian cobordism
shown on Figure \ref{fig:generators} in its bottom-top tangle presentation.
This cobordism will be interpreted in \S \ref{subsec:claspers} as the result of a ``clasper'' surgery.
As will be explained in Remark \ref{rem:generation}, 
it follows from clasper calculus that the monoidal category $\LCob$ is generated by 
\begin{equation}
\label{eq:generators_LCob}
\left(\psi_{1,1}^{\pm 1}, \mu,\eta,\Delta,\epsilon,S^{\pm 1}, v_\pm , Y \right).
\end{equation}

\begin{example}
The Poincar\'e sphere is the result of surgery in  $S^3$ along the $(+1)$-framed right-handed trefoil.
The punctured Poincar\'e sphere decomposes as $Y \circ (v_+ \otimes v_+ \otimes v_+)$.
\end{example}

\begin{example}
Another cobordism of interest is the ``co-duality'' $c\in \sLCob(0,2)$.
As a bottom tangle, this is
$$
c:= \begin{array}{c}{\relabelbox \small
\epsfxsize 1truein \epsfbox{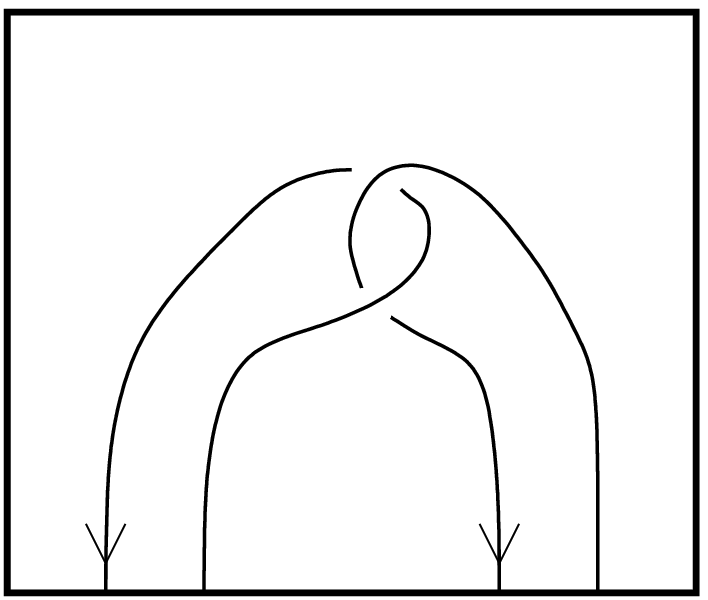}
\endrelabelbox} \end{array}.
$$
As observed in \cite{Kerler_towards}, $c$ decomposes as $\left(\mu \otimes \mu\right) \circ 
\left( \Id_1 \otimes \Delta \otimes \Id_1\right) \circ \left(v_- \otimes v_+ \otimes v_-\right)$.
\end{example}

Finally, we can deduce from the previous discussion a system of generators for the monoidal category $\LqCob$.
For this, equip the generators  $\psi_{1,1}^{\pm 1}, \mu,\eta,\Delta,\epsilon,S^{\pm 1}, v_\pm$ of $\sLCob$
with the only possible parenthesizings,
and lift $Y\in \LCob(3,0)$ to 
$$
Y\in \LqCob\left(((\bullet \bullet) \bullet),\varnothing\right).
$$
Also, for all non-associative words $u,v,w$ of total length $g:= |u| + |v| + |w|$, let 
$$
P_{u,v,w}: (u (v w)) \to ((u v) w) 
\quad \hbox{ and } \quad P_{u,v,w}^{-1}: ((u v) w) \to (u (v w))
$$ 
be the lifts of $\Id_g \in \LCob(g,g)$.
Then, the monoidal category $\LqCob$ is generated by the morphisms
\begin{equation}
\label{eq:generators_LqCob}
\left(\psi_{1,1}^{\pm 1}, \mu,\eta,\Delta,\epsilon,S^{\pm 1}, v_\pm , Y,
\left(P^{\pm 1}_{u,v,w}\right)_{u,v,w} \right).
\end{equation}

\subsection{Values of $\Ztilde$ on the generators}

Thus, it is important to compute the functor $\Ztilde$ 
on each of the morphisms listed in (\ref{eq:generators_LqCob}). 
Let us give some elements of computation, starting with the following

\begin{lemma}
\label{lem:values}
Here are the exact values of $\Ztilde$ on some of the generators of the category $\LqCob$:
$$
\begin{array}{rcll}
\Ztilde\left(\eta\right) &=& \varnothing  & \in \A(\{1^-\}) \\
\Ztilde\left(\epsilon \right) &=& \varnothing & \in \A(\{1^+\})\\
\Ztilde(v_+) &=& 
\chi^{-1} \exp\left(-\frac{1}{2}\!\figtotext{20}{10}{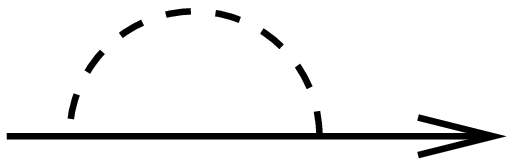}\!\!{}^{1^-}\right)
& \in \A(\{1^-\})  \\
\Ztilde(v_-) &=& \chi^{-1} \exp\left(\frac{1}{2}\!\figtotext{20}{10}{framing.eps}\!\!{}^{1^-}\right)
& \in \A(\{1^-\}).
\end{array}
$$
\end{lemma}

\begin{proof}
We have $Z(\eta) =\varnothing$ and
$Z(v_{\pm}) = \exp\left(\mp\frac{1}{2}\!\figtotext{20}{10}{framing.eps}\!\!{}^{1^-}\right)$
by our normalization of the Kontsevich integral:
We deduce the values of $\Ztilde$ on $\eta$ and $v_\pm$.
Also, one has by definition that
\begin{eqnarray*}
\Ztilde\left(\epsilon\right) &=& 
\chi^{-1} Z(\epsilon) \circ \mathsf{T}_{1}\\
&=& \left(\chi^{-1} Z(\Id_1) / i^- \mapsto 0\right)\circ \mathsf{T}_{1}\\
&=& \left(\left.\chi^{-1} Z(\Id_1) \circ \mathsf{T}_{1} \right/ i^- \mapsto 0 \right)\\
&=& \left( \left. \Ztilde(\Id_1) \right/ i^- \mapsto 0 \right)\\
&=& \left( \left. \left[\strutgraph{1^-}{1^+}\right] \right/ i^- \mapsto 0 \right)
\quad \quad   \quad  \quad = \ \varnothing.\\
\end{eqnarray*}
\end{proof}

We can already observe the following fact about special Lagrangian cobordisms:

\begin{corollary}
For all $M \in \sLCob_q(w,v)$, we have that
$\left(\left. \Ztilde(M) \right/ i^-\mapsto 0\right) = \varnothing$.
\end{corollary}

\begin{proof}
If $M \in \sLCob_q(w,v)$, then 
$\epsilon^{\otimes v} \circ M = \epsilon^{\otimes w}$.
We deduce that $\varnothing^{\otimes f} \circ \Ztilde(M) = \varnothing^{\otimes g}$
where $f:=|v|$ and $g:=|w|$, so that $\varnothing \circ \Ztilde(M) = \varnothing$.
\end{proof}

Lemma \ref{lem:values} is generalized by the following lemma, which reduces the computation of 
the LMO functor on a special Lagrangian $q$-cobordism to the Kontsevich integral.

\begin{lemma} 
\label{lem:special}
Let $M\in \sLqCob(w,v)$ where $w$ and $v$ are non-associative words 
of length $g$ and $f$ respectively.
Present $M$ as a bottom-top tangle in the following way:
$$
M =\begin{array}{c}
{\relabelbox \small
\epsfxsize 2truein \epsfbox{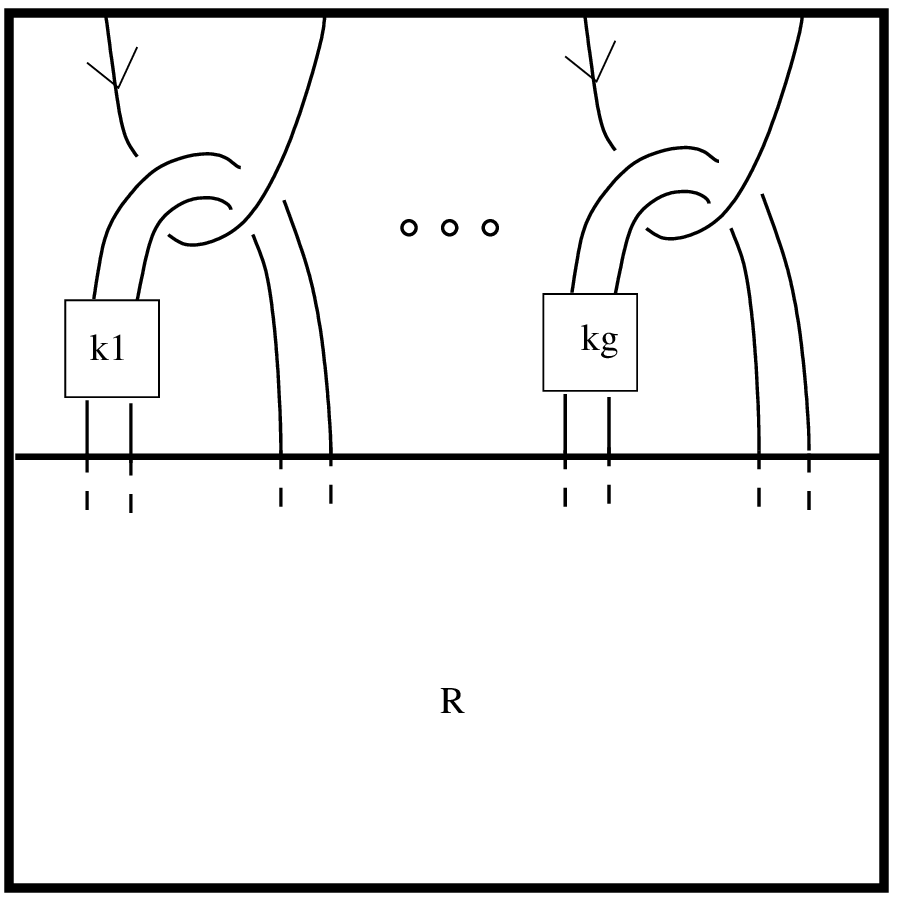}
\adjustrelabel <-0.15cm,-0cm> {k1}{$\downarrow^{w_1}$}
\adjustrelabel <-0.25cm,-0cm> {kg}{$\downarrow^{w_g}$}
\adjustrelabel <0cm,-0cm> {R}{$L$}
\endrelabelbox}
\end{array}
$$
where $L$ is a tangle in $[-1,1]^3$ and where $w_1,\dots,w_g$ are non-associative words in the letters $(+,-)$.
Equip $L$ with the non-associative words
$$
w_b(L) := \left(v/ \bullet \mapsto (+-)\right)
\quad \hbox{and} \quad
w_t(L) := \left(w/ \bullet_i \mapsto \left(w_i w_i^{op}\right) \right) 
$$
where the ``opposite'' word $w_i^{op}$ is obtained from $w_i$ by reading from right to left and 
by changing all the signs. 
Then, $\Ztilde(M)$ can be computed from $Z(L)$ as follows:
$$
\Ztilde(M) = \chi^{-1}\left(
\begin{array}{c}
{\relabelbox \small
\epsfxsize 2.2truein \epsfbox{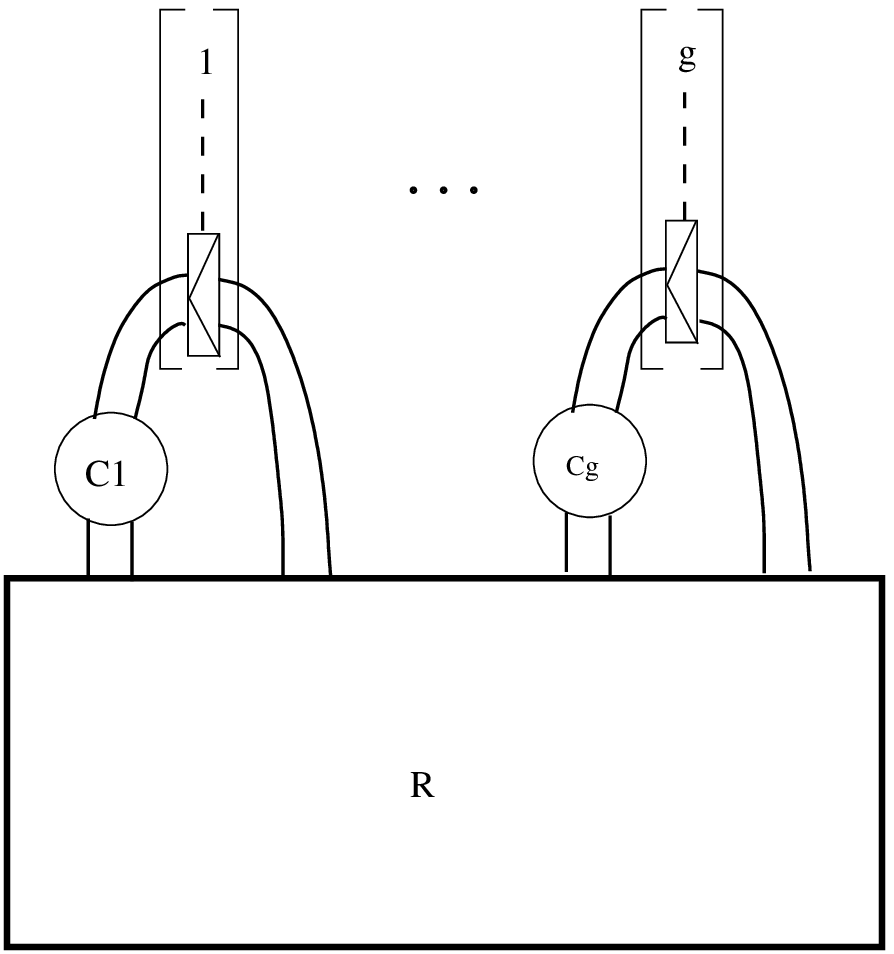}
\adjustrelabel <-0.1cm,-0cm> {1}{$1^+$}
\adjustrelabel <-0.1cm,-0cm> {g}{$g^+$}
\adjustrelabel <-0.1cm,0.05cm> {C1}{$C_{w_1}$}
\adjustrelabel <-0.15cm,0.05cm> {Cg}{$C_{w_g}$}
\adjustrelabel <-0.1cm,-0cm> {R}{$Z(L)$}
\endrelabelbox}
\end{array} \right) \quad \in \tsA(g,f).
$$
In that formula, the brackets denote exponentials,
a directed rectangle means
$$
{\relabelbox \small
\epsfxsize 4.2truein \epsfbox{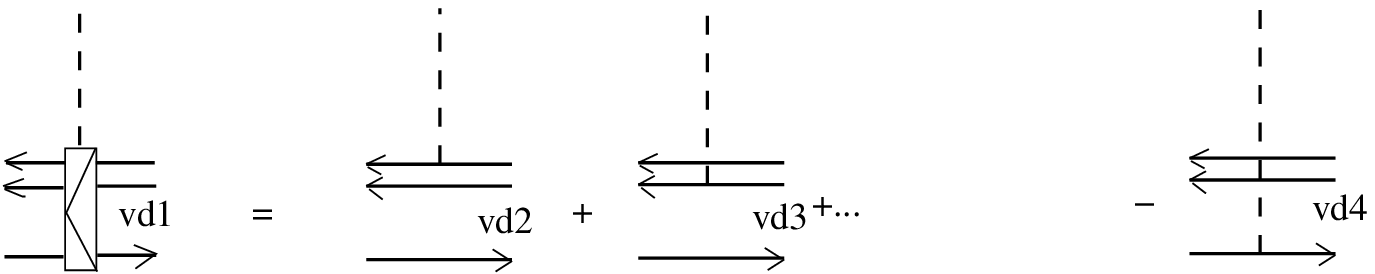}
\adjustrelabel <-0cm,-0cm> {=}{$:=$}
\adjustrelabel <-0cm,-0cm> {+}{$+$}
\adjustrelabel <0.5cm,-0cm> {+...}{$\pm \ \cdots$}
\adjustrelabel <-0cm,-0cm> {-}{$-$}
\adjustrelabel <-0cm,-0.1cm> {vd1}{$\vdots$}
\adjustrelabel <-0cm,-0.1cm> {vd2}{$\vdots$}
\adjustrelabel <-0cm,-0.1cm> {vd3}{$\vdots$}
\adjustrelabel <-0cm,-0.1cm> {vd4}{$\vdots$}
\endrelabelbox}
$$
and, for all non-associative word $u$ in the letters $(+,-)$,
$C_u \in \A\left(\downarrow^u\right)$ is the \emph{doubling anomaly} defined by the following axioms:
\begin{itemize}
\item[$(c_1)$] $C_\varnothing = \varnothing \ \in \A(\varnothing)$ and
$C_{(+)} = \downarrow \ \in \A(\downarrow)$.
\item[$(c_2)$] If $u$ is obtained from $u'$ by changing its $i$-th letter, then $C_u= S_i\left(C_{u'}\right)$
where $S_i$ is the ``orientation-reversal'' map on the $i$-th component of $\downarrow^{u'}$.
\item [$(c_3)$] If $u= (u_1 u_2)$ is the concatenation of two words $u_1$ and $u_2$, then
$$
C_u = \left( C_{u_1} \otimes C_{u_2} \right) 
\circ \Delta_{u_1,u_2}^{++} \Bigg( \underbrace{Z\left(\begin{array}{c} \figtotext{40}{20}{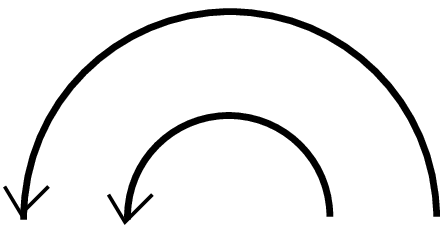}\\
\left((++)(--)\right)
\end{array}\right)}_{\in\ \A(\! \figtotext{15}{7}{doublecap.eps} \!)\ =\ \A(\downarrow \downarrow)}
\Bigg) \in \A(\downarrow^{u_1} \downarrow^{u_2}).
$$
\end{itemize}
In particular, if the Drinfeld associator $\Phi$ is assumed to be even, we have
$$
{\relabelbox \small
\epsfxsize 2.3truein \epsfbox{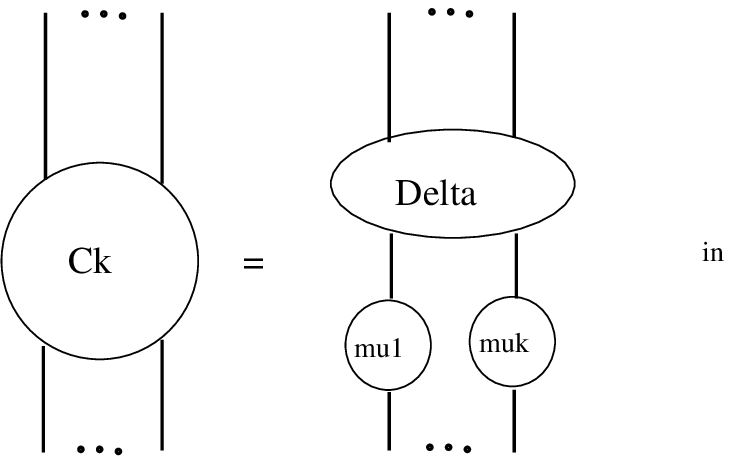}
\adjustrelabel <-0cm,-0cm> {Ck}{$C_u$}
\adjustrelabel <-0cm,-0cm> {=}{$=$}
\adjustrelabel <-0cm,-0cm> {in}{$\in \A\left(\downarrow^u\right).$}
\adjustrelabel <-0.1cm,0.1cm> {Delta}{$\Delta_u^+\left(\sqrt{\nu}\right)$}
\adjustrelabel <-0cm,-0cm> {mu1}{$\frac{1}{\sqrt{\nu}}$}
\adjustrelabel <-0cm,-0cm> {muk}{$\frac{1}{\sqrt{\nu}}$}
\adjustrelabel <-0cm,-0cm> {=}{$=$}
\endrelabelbox}
$$
\end{lemma}

\begin{proof}
Let $\gamma$ be the bottom-top $q$-tangle in $[-1,1]^3$ 
corresponding to $M$. Then, as suggested by the above figure, 
$\gamma$ decomposes as $L \cup U$ where $L$ is as described in the statement 
and where $U$ is the bottom-top $q$-tangle in $[-1,1]^3$ obtained from  $\Id_g\in \btT(g,g)$ 
(which is equipped with the word $\left(w/ \bullet \mapsto (+-)\right)$ at the top and at the bottom)
in the following manner: We double the bottom components and, if necessary, 
we reverse orientation of new components 
following the ``instructions'' given by $w_1,\dots,w_g$.
By definition, one has 
$$
\Ztilde(M) = \chi^{-1} Z([-1,1]^3,\gamma) \circ \mathsf{T}_g = 
\chi^{-1}\left( Z(L) \circ Z(U)\right) \circ \mathsf{T}_g  
$$
where the second $\circ$ denotes the composition in $\A$. Let
$$
\Delta: \A(\Id_g)=\A(\cupright^{\set{g}} \capleft^{\set{g}}) \simeq 
\A(\cupright^{\set{g}} \downarrow^{+ \cdots +}) \to \A(\cupright^{\set{g}} \downarrow^{w_1} \cdots \downarrow^{w_g})
\simeq \A(U)
$$
be the map $\Delta^{+ \cdots +}_{w_1,\dots,w_g}$ defined in Notation \ref{not:doubling}.
We also set 
$$
C:= (C_{w_1}' \otimes \Id_{w_1^{op}}) \otimes \cdots \otimes (C_{w_g}' \otimes \Id_{w_g^{op}})
\in \A(w_1 w_1^{op} \cdots w_g w_g^{op} , w_1 w_1^{op} \cdots w_g w_g^{op})
$$
where, for all non-associative word $u$ in the letters $(+,-)$, 
$$
C'_u:= Z\left(\Delta_u^+\left(\begin{array}{c}\capleft\\ {}^{(+-)} \end{array}\right)\right) 
\ \in \A\left( \Delta_u^+( \capleft) \right) = \A(\downarrow^u)
$$ 
is the Kontsevich integral of the $q$-tangle obtained by doubling $\capleft$ and by reversing
orientation of new components, in accordance with the ``instructions'' read in $u$. 
Then, by adapting the argument of \cite[Lemma 4.1]{LM3}, we see that 
$Z(U) = C \circ \Delta Z(\Id_g)$ in the category $\A$.
We deduce that
$$
\Ztilde(M) =
\chi^{-1}\left( Z(L) \circ C \circ \Delta Z(\Id_g) \right)  \circ \mathsf{T}_g. 
$$ 
Alternatively, we can write
$$
\Ztilde(M) = 
\chi^{-1}_{\pi_0(U^+)}\left(
\chi^{-1}_{\pi_0(L\cup U^-)}\left(
Z(L) \circ C \circ \Delta Z(\Id_g)
\right) \right)\circ \mathsf{T}_g
$$
where $U^+$ and $U^-$ denote the top and the bottom components of $U$ respectively.

\begin{claim}
\label{claim:gluing}
Let $X \cup Y$ be the gluing of two compact oriented $1$-manifolds $X$ and $Y$, 
such that all of the connected components of $X$, $Y$ and $X\cup Y$ are intervals.
Then, there exists a $\lambda \in \A\left(\pi_0(X \cup Y) \cup \pi_0(X) \cup \pi_0(Y)\right)$ such that
$$
\forall a \in \A(\pi_0(X)), \ \forall b \in \A(\pi_0(Y)), \quad
\chi^{-1} \left(\chi(a) \cup \chi(b) \right) = 
\langle  \lambda , a \sqcup b \rangle_{\pi_0(X) \cup \pi_0(Y)}.
$$
\end{claim}

\begin{proof}[Proof of Claim \ref{claim:gluing}] 
The argument appears in \cite[Prop. 5.4]{BGRT2}. For instance, assume that
$$
X \cup Y = \begin{array}{c} {\relabelbox \small
\epsfxsize 2.1truein \epsfbox{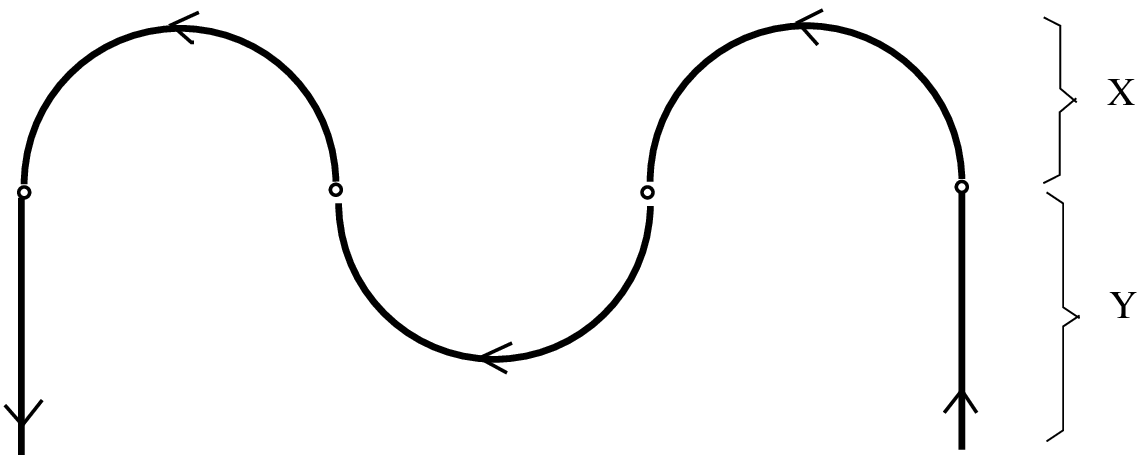}
\adjustrelabel <-0cm,-0.1cm> {X}{$Y$}
\adjustrelabel <-0cm,-0.1cm> {Y}{$X$}
\endrelabelbox}\end{array}
$$
and define 
$$
\lambda' :=
\begin{array}{c} {\relabelbox \small
\epsfxsize 2.1truein \epsfbox{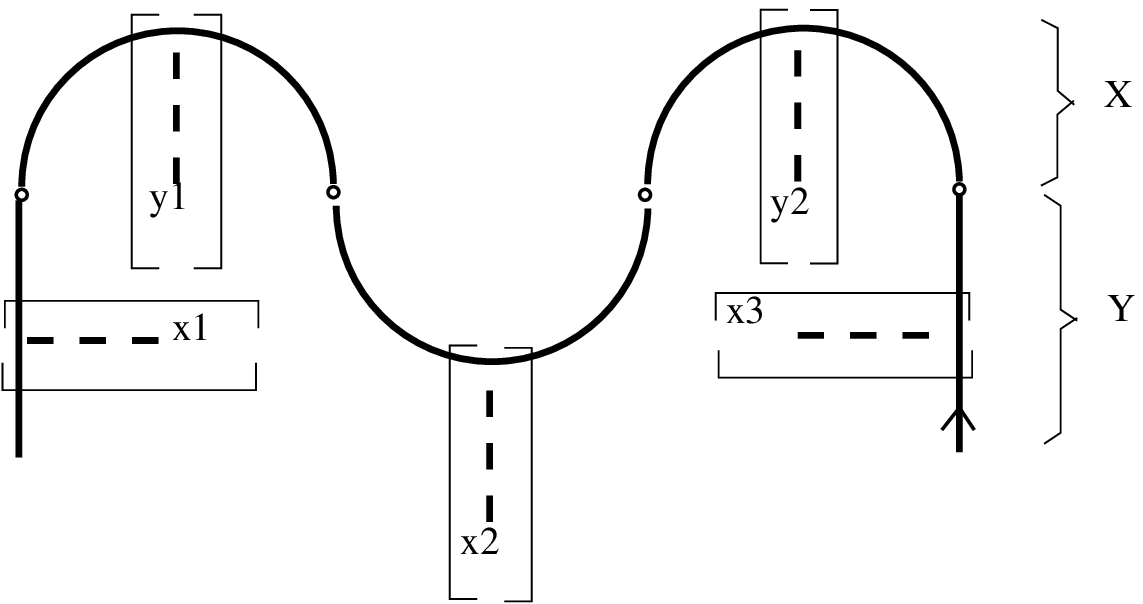}
\adjustrelabel <-0cm,-0.1cm> {X}{$Y$}
\adjustrelabel <-0cm,-0.1cm> {Y}{$X$}
\adjustrelabel <-0cm,-0.1cm> {x1}{$x_1$}
\adjustrelabel <-0cm,-0.1cm> {x2}{$x_2$}
\adjustrelabel <-0cm,-0.1cm> {x3}{$x_3$}
\adjustrelabel <-0cm,-0.1cm> {y1}{$y_1$}
\adjustrelabel <-0cm,-0.1cm> {y2}{$y_2$}
\endrelabelbox}\end{array} \quad \in \A\left(X\cup Y, \pi_0(X) \cup \pi_0(Y)\right) 
$$
where $x_1,x_2,x_3$ and $y_1,y_2$ index the connected components of $X$ and $Y$ respectively. 
Then, observe that
$$
\chi(a) \cup \chi(b) = \langle \lambda' , a\sqcup b \rangle_{\{x_1,x_2,x_3\} \cup \{y_1,y_2\}}
$$
and set $\lambda := \chi^{-1}(\lambda')$.
\end{proof}

\noindent
So, there exists a certain $\lambda$ such that
\begin{eqnarray*}
\Ztilde(M) &=& \left\langle \chi^{-1}_{\pi_0(L)} 
\left(Z(L) \circ C \right) \sqcup \chi^{-1}_{\pi_0(U)}  \Delta Z(\Id_g), 
\lambda \right\rangle_{\pi_0(L) \sqcup \pi_0(U^-)} \circ \mathsf{T}_g\\
 &=& \left\langle \chi^{-1}_{\pi_0(L)}  \left(Z(L) \circ C\right)
\sqcup \left(\chi^{-1}_{\pi_0(U)} \Delta Z(\Id_g) 
\circ \mathsf{T}_g \right), \lambda \right\rangle_{\pi_0(L) \cup \pi_0(U^-)}.
\end{eqnarray*}
Recall that the composite map $\chi^{-1}_{\pi_0(U^-)} \Delta \chi_{\pi_0(\Id_g^-)}$ 
can be described as an operation at the level of the external vertices colored
with bottom components of $\Id_g$. So, we obtain
\begin{eqnarray*}
 \chi^{-1}_{\pi_0(U)} \Delta Z(\Id_g) 
\circ \mathsf{T}_g &=& 
\chi^{-1}_{\pi_0(U^-)} \Delta \chi_{\pi_0(\Id_g^-)} \chi^{-1}_{\pi_0(\Id_g)} Z(\Id_g) \circ \mathsf{T}_g \\
&=& 
\chi^{-1}_{\pi_0(U^-)} \Delta \chi_{\pi_0(\Id_g^-)}
\left( \chi^{-1}_{\pi_0(\Id_g)}  Z(\Id_g)  \circ \mathsf{T}_g \right)\\
&=& 
\chi^{-1}_{\pi_0(U^-)} \Delta \chi_{\pi_0(\Id_g^-)}
\left( \left[ \sum_{i=1}^g \strutgraph{i^-}{i^+}\right] \right)
\end{eqnarray*}
\begin{eqnarray*}
&=& \chi^{-1}_{\pi_0(U^-)} \left(\begin{array}{c}
{\relabelbox \small
\epsfxsize 2truein \epsfbox{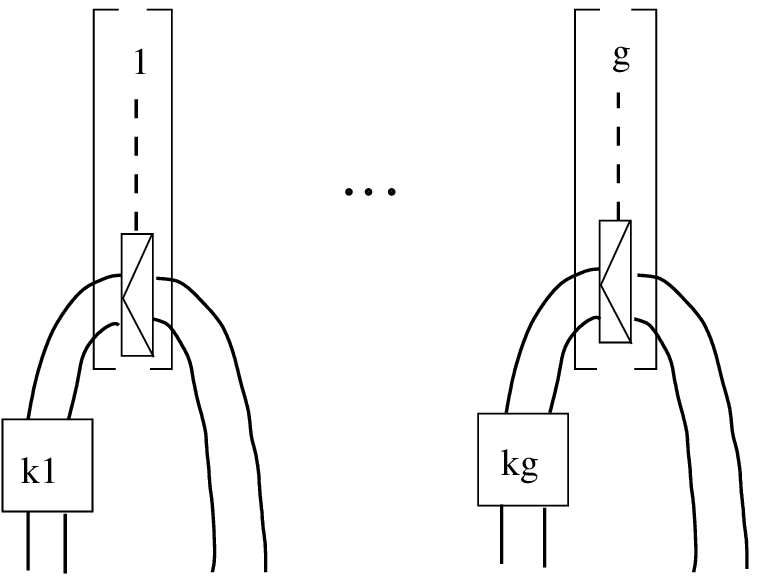}
\adjustrelabel <-0.1cm,-0cm> {1}{$1^+$}
\adjustrelabel <-0.1cm,-0cm> {g}{$g^+$}
\adjustrelabel <-0.1cm,-0cm> {k1}{$\downarrow^{w_1}$}
\adjustrelabel <-0.1cm,-0cm> {kg}{$\downarrow^{w_g}$}
\endrelabelbox}
\end{array} \right).
\end{eqnarray*}
A reverse application of Claim \ref{claim:gluing}
proves the lemma if each $C_{w_i}$ is replaced by $C_{w_i}'$.

Thus, we should now prove that $C_u =C'_u$ for all non-associative word $u$. 
The family $(C'_u)_u$ clearly satisfies axioms $(c_1)$ and $(c_2)$. 
Compare the identity
\begin{eqnarray*}
C'_{(u_1 u_2)} &=& \Delta^{+++}_{u_1,u_2, u_2^{op} u_1^{op}}(\Phi) 
\circ \left( \Id_{u_1} \otimes \Delta^{+++}_{u_2,u_2^{op},u_1^{op}}(\Phi^{-1}) \right)
\circ \left( \Id_{u_1} \otimes C'_{u_2} \otimes \Id_{u_1^{op}}\right) \circ C'_{u_1}\\
&=& \left( C'_{u_1} \otimes C'_{u_2} \otimes \Id_{u_2^{op}} \otimes \Id_{u_1^{op}} \right) 
\circ \Delta^{+++}_{u_1,u_2, u_2^{op} u_1^{op}}(\Phi) 
\circ \left( \Id_{u_1} \otimes \Delta^{+++}_{u_2,u_2^{op},u_1^{op}}(\Phi^{-1}) \right)\\
&&
\circ \left( \Id_{u_1} \otimes \Delta_{u_2}^+\left(\capleft\right) 
\otimes \Id_{u_1^{op}}\right) \circ \Delta_{u_1}^+\left(\capleft \right)
\end{eqnarray*}
in the category $\A$, to the identity
$$
C'_{(++)} = \Delta^{+++}_{+,+, --}(\Phi) 
\circ \left( \Id_{+} \otimes \Delta^{+++}_{+,-,-}(\Phi^{-1}) \right)
\circ \left( \Id_+ \otimes \capleft \otimes \Id_{-}\right) \circ \capleft
$$
to conclude that the family $(C'_u)_u$ has the property $(c_3)$ as well.

Finally, let us assume that the associator $\Phi$ is even. 
Let $\widehat{Z}_{f}^0$ denote the Kontsevich integral \emph{as normalized} 
by Le and Murakami in \cite{LM3}.
Observe that it differs from our $Z$ and from the invariant $\widehat{Z}_f$ introduced in \cite{LM1,LM2}, 
by the values it takes on the ``cap'' and the ``cup'':
$$
\widehat{Z}_{f}^0\left(\begin{array}{c}\capleft\\ {}^{(+-)} \end{array}\right) 
= \! \begin{array}{c}{\relabelbox \small
\epsfxsize 0.6truein \epsfbox{diagram_on_cap.eps}
\adjustrelabel <-0.15cm,0.05cm> {D}{${}_{\sqrt{\nu}}$}
\endrelabelbox} \end{array}  \in  \A\left(\capleft\right)
\quad  \quad  \quad 
\widehat{Z}_{f}^0\left(\begin{array}{c} {}_{(+-)} \\ \cupright  \end{array}\right) 
= \! \begin{array}{c}{\relabelbox \small
\epsfxsize 0.6truein \epsfbox{diagram_on_cup.eps}
\adjustrelabel <-0.25cm,0.05cm> {D}{${}_{\sqrt{\nu}}$}
\endrelabelbox} \end{array} \in \A\left(\cupright\right).
$$
Also, recall that $\widehat{Z}_{f}^0$ commutes with the ``doubling'' maps whatever the type of the $q$-tangle is. 
Thus, by denoting $\pi_0:=\pi_0(\Delta_u^+(\capleft))$, we obtain that
\begin{eqnarray*}
C'_u  = C_u &=& \widehat{Z}_f^0 
\left(\Delta_u^+\left(\begin{array}{c}\capleft\\ {}^{(+-)} \end{array}\right)\right) 
\sharp_{\pi_0} \left(1/\sqrt{\nu}\right)^{\otimes \pi_0}\\
&=& \Delta_u^+ \widehat{Z}_f^0 \left(\begin{array}{c}\capleft\\ {}^{(+-)} \end{array}\right) 
\sharp_{\pi_0} \left(1/\sqrt{\nu}\right)^{\otimes \pi_0}\\
&=& 
\Delta_u^+(\sqrt{\nu}) \sharp_{\pi_0} \left(1/\sqrt{\nu}\right)^{\otimes \pi_0}.
\end{eqnarray*}
\end{proof}

We can apply Lemma \ref{lem:special} to obtain the exact values of 
$\Ztilde$ on $\Psi_{1,1}^{\pm 1}$ and $S^{\pm 1}$. The same can be done for $\mu, \Delta$
and $P^{\pm 1}_{u,v,w}$, but the answers then depend on the associator $\Phi$.
This completes our discussion on the values of $\Ztilde$ on the generators of the monoidal category $\sLqCob$.\\

\subsection{Low-degree computations of $\Ztilde$}

We conclude this computational section by some explicit computations in low degree. 

\begin{lemma}
\label{lem:T}
Assume that the elected Drinfeld associator $\Phi$ is even. Then, one has 
$$
\mathsf{T}_{1}  = \left[\strutgraph{1-}{1+}\right]
\sqcup \left( \varnothing - \frac{1}{8}\cdot \phigraph{1^-}{1^+} - 
\frac{1}{48}\cdot \phigraphtop{1^+}{1^+} 
+ \frac{1}{8}\cdot \Hgraph{1^-}{1^-}{1^+}{1^+}
+ (\ideg >2) \right) \ \in \tsA(1,1).
$$
\end{lemma}

\begin{proof}
It follows from Theorem \ref{th:LMO_functor} that 
$\mathsf{T}_1$ is a unit of the algebra $\left(\tsA(1,1),\circ\right)$
with inverse $\chi^{-1}Z(\Id_{1})$, 
where $\Id_{1}$ denotes the bottom-top tangle from Figure \ref{fig:bt_tangle_id}.
The Kontsevich integral of this tangle can easily be computed in low degree,
using the fact that the Drinfeld associator is equal to
$$
\Phi = 1 + \frac{1}{24} \cdot \figtotext{30}{30}{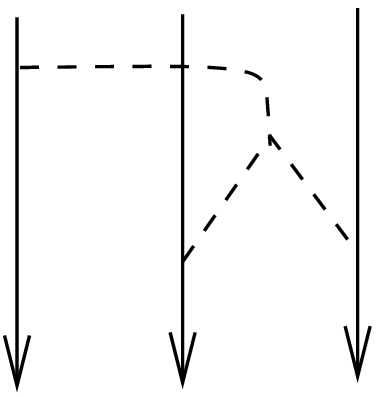} + (\deg > 3)
$$
since it has been assumed to be even. 
The reader may check that such a computation gives the following result:
$$
\chi^{-1} Z( \Id_1) = \left[\strutgraph{1^-}{1^+}\right]
\sqcup \left( \varnothing + \frac{1}{8}\cdot \phigraph{1^-}{1^+} + 
\frac{1}{48}\cdot \phigraphtop{1^+}{1^+} 
- \frac{1}{8}\cdot \Hgraph{1^-}{1^-}{1^+}{1^+}
+ (\ideg >2) \right).
$$
\end{proof}

Next proposition computes the functor $\Ztilde$ up to internal degree $2$:

\begin{proposition}
\label{prop:low_degree} 
Assume that the elected Drinfeld associator $\Phi$ is even.
Then, the functor $\Ztilde$ reduced modulo $(\ideg >2)$ takes the following values:
$$
\begin{array}{|c|c|c|}
\hline
\begin{array}{c} \\ M \\ \\\end{array} & 
\ \frac{\Lk(M)}{2} = \log_\sqcup \Ztilde^s(M) \ & \log_\sqcup \ZtildeY(M) \ {\hbox{mod } (\ideg>2)} \\[0.1cm]
\hline 
\eta & 0 & 0 \\
\hline 
\epsilon & 0 & 0 \\
\hline
P^{\pm 1}_{u,v,w} & \strutgraph{1^-}{1^+} + \strutgraph{2^-}{2^+} + \strutgraph{3^-}{3^+} & 0 \\
\hline
v_\pm & \mp \frac{ 1}{2} \cdot \strutgraphbot{1^-}{1^-} &  \frac{1}{48} \cdot \phigraphbot{1^-}{1^-}\\
\hline
S^{\pm 1} & - \strutgraph{1^-}{1^+} &  \mp \frac{1}{4} \cdot \phigraph{1^-}{1^+}
\mp \frac{1}{4} \cdot \Hgraph{1^-}{1^-}{1^+}{1^+}\\
\hline
\psi_{1,1}^{\pm 1} &  \strutgraph{2^-}{1^+} + \strutgraph{1^-}{2^+} &  
\mp \frac{1}{2} \cdot  \Hgraph{2^-}{1^-}{1^+}{2^+}\\
\hline 
\mu & \strutgraph{1^-}{1^+} + \strutgraph{1^-}{2^+} &  -\frac{1}{2} \cdot \Ygraphbottoptop{1^-}{1^+}{2^+}
+ \frac{1}{12} \cdot \Hgraphtopright{1^-}{1^+}{1^+}{2^+} + \frac{1}{12} \cdot \Hgraphtopleft{1^-}{1^+}{2^+}{2^+}\\
\hline
\Delta & \strutgraph{1^-}{1^+} + \strutgraph{2^-}{1^+} &  \frac{1}{2} \cdot \Ygraphbotbottop{1^-}{2^-}{1^+}
+ \frac{1}{12} \cdot \Hgraphbotleft{1^-}{2^-}{2^-}{1^+} + \frac{1}{12} \cdot \Hgraphbotright{1^-}{1^-}{2^-}{1^+}
- \frac{1}{4} \cdot \Hgraph{1^-}{2^-}{1^+}{1^+}\\
\hline
Y & 0 &  - \Ygraphtop{1^+}{2^+}{3^+} + \frac{1}{2} \cdot \Hgraphcup{1^+}{1^+}{2^+}{3^+}
+ \frac{1}{2} \cdot \Hgraphcup{2^+}{2^+}{3^+}{1^+} + \frac{1}{2} \cdot \Hgraphcup{3^+}{3^+}{1^+}{2^+}\\
\hline
c & -\strutgraphbot{1^-}{2^-} &  \frac{1}{8} \cdot \phigraphbot{1^-}{2^-} +
\frac{1}{8} \cdot \Hgraphcap{1^-}{1^-}{2^-}{2^-} \\
\hline
\end{array}
$$
\end{proposition}

\begin{proof}
Except for $\Ztilde(Y)$, all the above values can be derived 
from Lemma \ref{lem:special} by computing the appropriate Kontsevich integral.
Computation details are left to the interested reader.
As for $\Ztilde(Y)$, it needs a special treatment. First, because 
$$
Y \circ \left( \eta \otimes \Id_2 \right)  
= Y \circ \left( \Id_1 \otimes \eta \otimes \Id_1 \right) = 
Y \circ \left( \Id_2 \otimes \eta \right) = \epsilon \otimes \epsilon,
$$
one must have
$$
\log_\sqcup \Ztilde(Y) =  a \cdot \Ygraphtop{1^+}{2^+}{3^+} + 
b_1 \Hgraphcup{1^+}{1^+}{2^+}{3^+}
+ b_2 \Hgraphcup{2^+}{2^+}{3^+}{1^+} + b_3 \Hgraphcup{3^+}{3^+}{1^+}{2^+} + (\ideg >2)
$$
for some $a,b_1,b_2,b_3 \in \Q$. 
Using Le's observation (see \cite[Prop. 1.3]{Le_Grenoble} or \cite[Lemma 11.22]{Ohtsuki}),
we obtain that 
$$
\chi^{-1}Z\left(\epsilon^{\otimes ((\bullet \bullet) \bullet)}\right) - \chi^{-1}Z(Y) 
= \Ygraphtop{1^+}{2^+}{3^+} + (\ideg >1).
$$
By composing on the right with $\mathsf{T}_3$, we obtain that
\begin{equation}
\label{eq:Le}
\Ztilde\left(\epsilon^{\otimes ((\bullet \bullet) \bullet)}\right)
- \Ztilde(Y) = \Ygraphtop{1^+}{2^+}{3^+}
+ (\ideg >1).
\end{equation}
Since $\Ztilde\left(\epsilon^{\otimes ((\bullet \bullet) \bullet)}\right)
= \Ztilde(\epsilon)^{\otimes 3}= \varnothing$, we obtain that $a=-1$.
Next, because 
$$
Y \circ (\Id_1 \otimes c ) = Y \circ (c \otimes \Id_1) = 
Y \circ (\Id_1 \otimes \psi_{1,1}^{-1}) \circ (c \otimes \Id_1) = \epsilon,
$$ 
one can deduce that $b_1=b_2=b_3 =1/2$.
\end{proof} 

\vspace{0.5cm}

\section{Functorial invariants of cobordisms between closed surfaces}

\label{sec:closed_surfaces}

In this section, we show that a suitable reduction of the LMO functor factors through 
the category of Lagrangian cobordisms between \emph{closed} surfaces. 
In particular, we recover the TQFT-like functor constructed by Le and the first author in \cite{CL1}.

\subsection{A reduction of the LMO functor}

Instead of considering cobordisms between surfaces with \emph{one boundary component} 
as we did in \S \ref{subsec:cob}, we can consider cobordisms between \emph{closed} surfaces. 
More precisely, we replace the surface $F_g$ by 
$$
\widehat{F_g}:= F_g \cup D
$$
where $D$ is a $2$-disk such that $\partial D= - \partial ([-1,1]^2\times 0)$,
and we replace the cube with tunnels and handles $C_{g_-}^{g_+}$  by
$$
\widehat{C_{g_-}^{g_+}}:= 
C_{g_-}^{g_+} \cup \left( D\times[-1,1] \right).
$$
Thus, one obtains a category of cobordisms $\widehat{\Cob}$ including a subcategory of Lagrangian cobordisms 
$\widehat{\LCob}$, which can further be enhanced to $\widehat{\LqCob}$. 
There is the functor
$$
\widehat{\centerdot} : \mathcal{C} \longrightarrow \widehat{\mathcal{C}}
\quad \quad \forall \mathcal{C} = \Cob, \LCob, \LqCob
$$
that consists in gluing the $2$-handle $D\times[-1,1]$ along the ``vertical'' boundary of cobordisms:
This functor is bijective at the level of objects and surjective at the level of morphisms.\\

On the diagrammatic side, we define a quotient category $\widehat\tsA$ of $\tsA$ as follows.
Let $\mathcal I(g,f)$ denote the closure of the subspace of $\tsA(g,f)$ 
spanned by the Jacobi diagrams $D$ that contain
\begin{gather*}
{\relabelbox \small
\epsfxsize 5.5truein \epsfbox{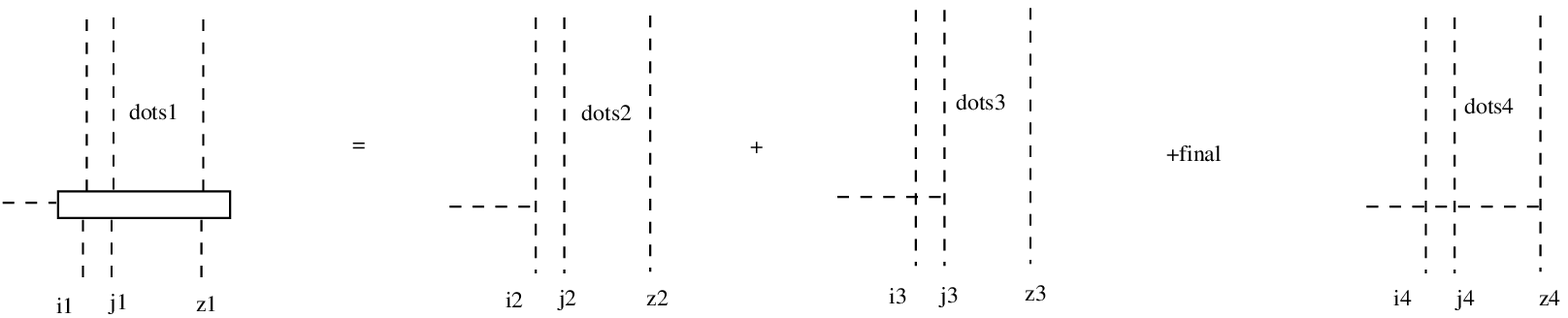}
\adjustrelabel <-0cm,-0.cm> {dots1}{$\dots$}
\adjustrelabel <-0cm,-0.cm> {dots2}{$\dots$}
\adjustrelabel <-0cm,-0.cm> {dots3}{$\dots$}
\adjustrelabel <-0cm,-0.cm> {dots4}{$\dots$}
\adjustrelabel <-0cm,-0.cm> {i1}{$i^-$}
\adjustrelabel <-0cm,-0.cm> {j1}{$j^-$}
\adjustrelabel <-0cm,-0.cm> {z1}{$z^-$}
\adjustrelabel <-0cm,-0.cm> {i2}{$i^-$}
\adjustrelabel <-0cm,-0.cm> {j2}{$j^-$}
\adjustrelabel <-0cm,-0.cm> {z2}{$z^-$}
\adjustrelabel <-0cm,-0.cm> {i3}{$i^-$}
\adjustrelabel <-0cm,-0.cm> {j3}{$j^-$}
\adjustrelabel <-0cm,-0.cm> {z3}{$z^-$}
\adjustrelabel <-0cm,-0.cm> {i4}{$i^-$}
\adjustrelabel <-0cm,-0.cm> {j4}{$j^-$}
\adjustrelabel <-0cm,-0.cm> {z4}{$z^-$}
\adjustrelabel <-0cm,-0.cm> {=}{$:=$}
\adjustrelabel <-0cm,-0.cm> {+}{$+$}
\adjustrelabel <-0cm,-0.cm> {+final}{$+\cdots +$}
\endrelabelbox}
\end{gather*}
In the part of $D$ that is drawn, one should see all the external vertices of $D$ 
colored with elements from  $\set{f}^-$: 
The labels are then denoted by $i^-,j^-,\dots ,z^-$ and may contain repetitions. 

\begin{remark}
\label{rem:top=bot}
By the IHX relation, the relation $\mathcal{I}$ can be re-written 
replacing labels from $\set{f}^-$ by labels from $\set{g}^+$.
\end{remark}
\noindent
Then, the quotient spaces
\begin{gather*}
  \widehat \tsA(g,f) := \tsA(g,f)/\mathcal I(g,f)
\end{gather*}
with $g,f\ge 0$, form a quotient category of $\tsA$. 
We note that, in contrast with $\tsA$, 
there is no natural monoidal structure for $\widehat\tsA$.

\begin{theorem}
\label{th:closed}
The LMO functor $\Ztilde$, reduced modulo $\mathcal{I}$, factors through $\widehat{\LqCob}$:
$$
\xymatrix{
{\LqCob} \ar[d]_{\widehat{\centerdot}} \ar[r]^{\Ztilde}  & {\tsA} \ar[d]^{\widehat{\centerdot}}\\
{\widehat\LqCob} \ar@{-->}[r]_{\Ztilde}  & {\widehat\tsA}
}
$$
\end{theorem}

\begin{figure}[h]
\centerline{\relabelbox \small
\epsfxsize 3truein \epsfbox{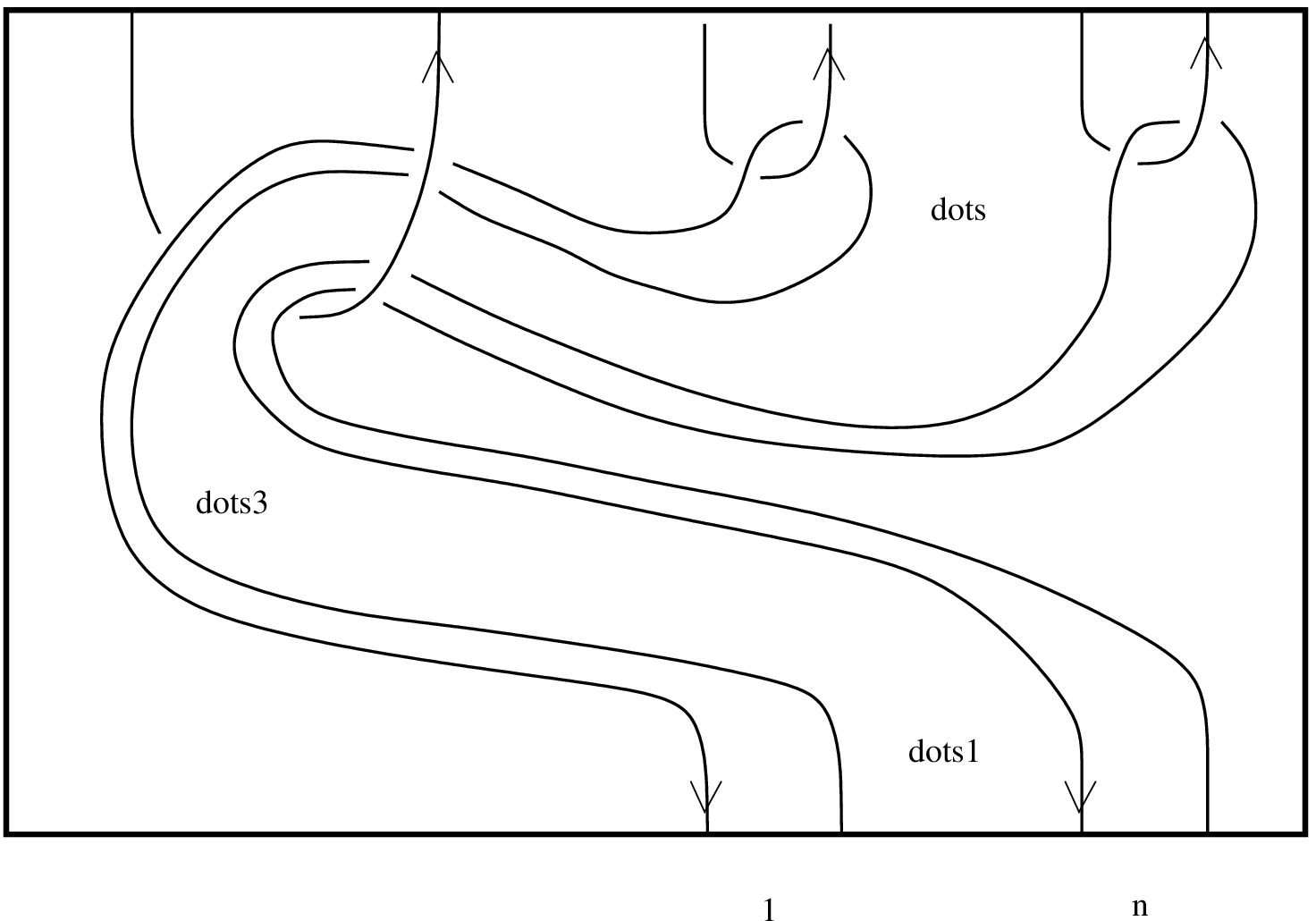}
\adjustrelabel <-0cm,-0.cm> {dots1}{$\dots$}
\adjustrelabel <-0cm,-0.cm> {dots}{$\dots$}
\adjustrelabel <-0cm,-0.cm> {dots3}{$\dots$}
\adjustrelabel <-0cm,-0.cm> {1}{$1$}
\adjustrelabel <-0cm,-0.cm> {n}{$n$}
\endrelabelbox}
\label{fig:adjoint_action}
\caption{The  cobordism $\alpha_n: n+1 \to n$.}
\end{figure}

In order to prove Theorem \ref{th:closed}, we need some further results. 
For all $n\geq 0$, let $\alpha_n \in \sLCob(n+1,n)$ be the special cobordism drawn
on Figure \ref{fig:adjoint_action} as a bottom-top tangle.
Let $\sim$ be the congruence relation on the category $\Cob$ generated by 
$$
\alpha_n \sim \epsilon \otimes \Id_n \quad \quad \forall n\geq 0.
$$ 
The quotient category is denoted by $\Cob/\sim$.
The following result gives an algebraic characterization of the ``kernel'' 
of the canonical functor $\Cob \to \widehat{\Cob}$.

\begin{proposition}
\label{prop:closing_cobordisms}
The functor $ \ \widehat{\centerdot}: \Cob \to \widehat{\Cob}$ induces an isomorphism
$$
\widehat{\centerdot} : \Cob/\sim\  \longrightarrow \widehat{\Cob}.
$$ 
A similar statement holds for $\LCob$ and $\sLCob$.
\end{proposition}

\begin{proof}
Let us consider ``bottom-top tangles'' $(\widehat{B},\gamma)$ 
which are defined as in Definition \ref{def:bt_tangle}, 
except that we are now considering cobordisms $\widehat{B}$ from $\widehat{F_0}$ to $\widehat{F_0}$
(instead of  cobordisms $B$ from $F_0$ to $F_0$). 
Such bottom-top tangles form a category $\widehat{\btT}$ and, again, there is an obvious functor
$\ \widehat{\centerdot} : \btT \to \widehat{\btT}$. 
Because Theorem \ref{th:iso} holds in this closed case as well
 -- which provides an isomorphism $\hbox{D}: \widehat{\btT} \to \widehat{\Cob}$ --
and because the $\hbox{D}$ construction obviously commutes with the $\widehat{\centerdot}$ operation,
the proposition is equivalent to showing that the functor $\widehat{\centerdot} : \btT \to \widehat{\btT}$
determines the congruence relation $\sim$.

So, let $(B_1,\gamma_1)$ and $(B_2,\gamma_2) \in \btT(g,f)$ be such that 
$(\widehat{B_1},\gamma_1)= (\widehat{B_2},\gamma_2)$. 
This precisely means that there exists a homeomorphism $\phi:\widehat{B_1} \to \widehat{B_2}$
which preserves the orientations and the boundary parameterizations,
and which sends $\gamma_1$ to $\gamma_2$.
Let $X_1$ and $X_2$ be the framed  string knots in $\widehat{B_1}$ and $\widehat{B_2}$ respectively,
corresponding to the $2$-handles $D\times [-1,1]$ attached to $B_1$ and $B_2$ respectively.
If $\phi(X_1)= X_2$, then $\phi$ restricts to a homeomorphism $B_1 \to B_2$, 
so that $(B_1,\gamma_1)=(B_2,\gamma_2)$. 
Otherwise, we perform a ``slam-dunk'' move (shown on Figure \ref{fig:slam-dunk})
along a framed $2$-component link $L:=(R,V) \subset \widehat{B_2} \setminus X_2$,
where $V$ is the ``veering'' $X_2$ must follow to become $\phi(X_1)$ and 
where $R$ rings $X_2$ and $V$: This move provides 
a homeomorphism $ \psi : \widehat{B_2} \to \left(\widehat{B_2}\right)_L$ sending $\phi(X_1)$ to $X_2$.
Thus, $\psi \circ \phi$ restricts to a homeomorphism $B_1 \to (B_2)_L$ which sends $\gamma_1$ to $(\gamma_2)_L$,
said more concisely: $(B_1,\gamma_1) = (B_2,\gamma_2)_L$. 
By adding an extra bottom component to $(B_2,\gamma_2)$ (which follows $V$ in $B_2$), 
we obtain $(B_2,\gamma_2')$ such that $\alpha_f \circ (B_2,\gamma'_2) = (B_2,\gamma_2)_L$ and,
of course, $(\epsilon \otimes \Id_f) \circ (B_2,\gamma'_2) = (B_2,\gamma_2)$.
We conclude that $(B_1,\gamma_1) \sim (B_2,\gamma_2)$.
\end{proof}

\begin{proof}[Proof of Theorem \ref{th:closed}]
For all non-associative word $v$ of length $n$ in the single letter $\bullet$, 
we lift $\alpha_n \in \LCob(n+1,n)$ to $\alpha_v \in \LqCob\left((\bullet v),v\right)$.
It suffices to prove that
\begin{claim}
\label{claim:in_I}
$\Ztilde(\alpha_v) - \Ztilde(\epsilon \otimes \Id_v) \in \mathcal{I}(n+1,n)$.
\end{claim}
Lemma \ref{lem:special} shows that $\Ztilde(\alpha_v)$ equals
$$
\chi^{-1} \left(\begin{array}{c} 
{\relabelbox \small
\epsfxsize 3.8truein \epsfbox{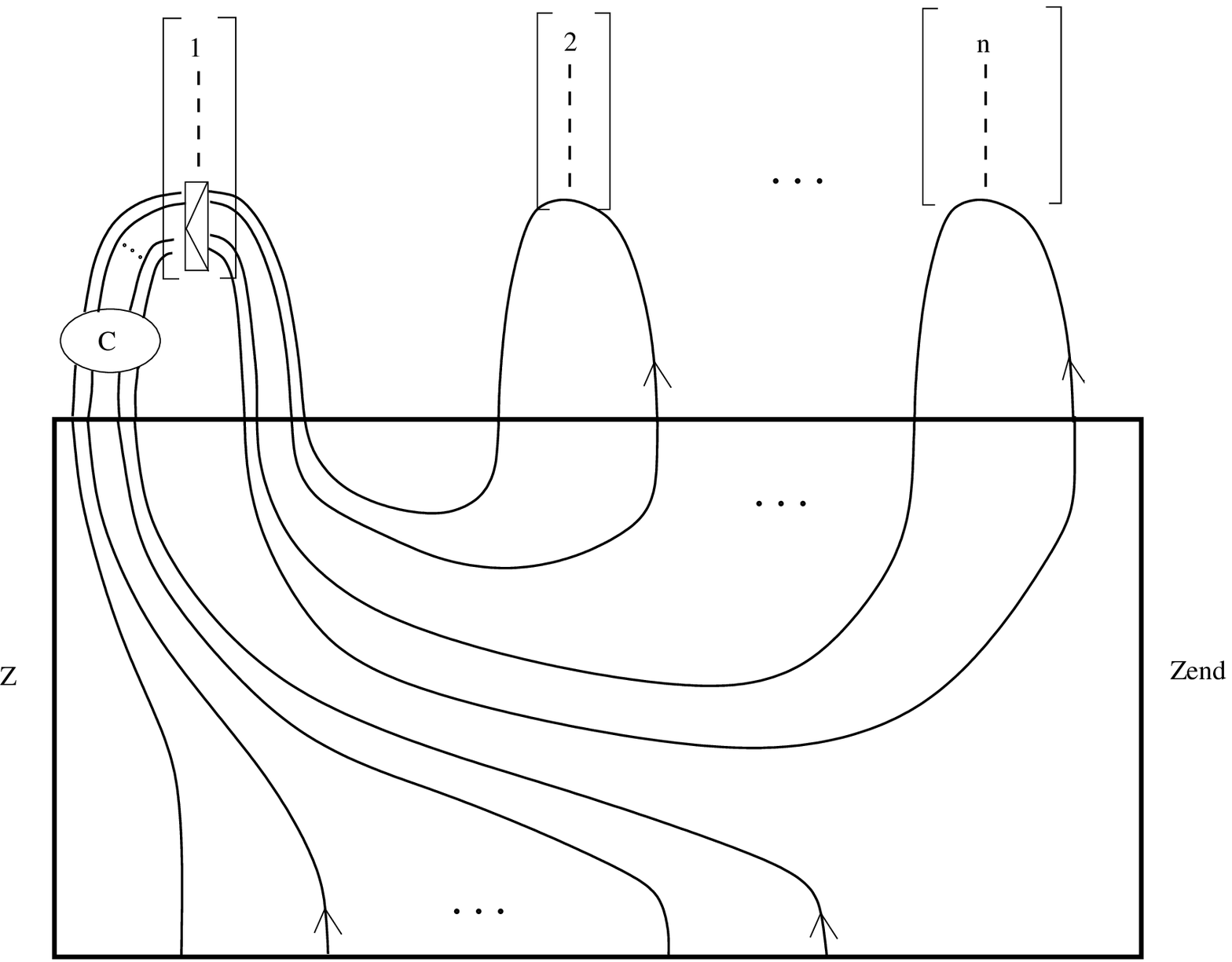}
\adjustrelabel <-0cm,-0.1cm> {2}{${}^{2^+}$}
\adjustrelabel <-0.3cm,-0.05cm> {n}{${}^{(n+1)^+}$}
\adjustrelabel <-0cm,-0.1cm> {1}{${}^{1^+}$}
\adjustrelabel <-0.1cm,-0.cm> {C}{$C_w$}
\adjustrelabel <-0.2cm,-0.cm> {Z}{$Z\Bigg($}
\adjustrelabel <-0cm,-0.1cm> {Zend}{$\Bigg)$}
\endrelabelbox}
\end{array} \right)
$$
where $w:= \left(v/ \bullet \mapsto (+-) \right)$. 
Let $L$ be the tangle encapsulated in the box, 
then $L$ is equipped with $w_b(L) = w $ and $w_t(L)= \left( (w w^{op}) w \right)$.
Thus, we have 
$$
\Ztilde(\alpha_v) = \chi^{-1}\left(Z(L) \circ U\right)
$$
where $U$ is the part outside the box, and $Z(L)$ is easily computed:
$$
Z(L) =\left( \Id_{2n} \otimes Z\left(\Delta^+_{w^{op}}(\cupright)\right) \right)
\circ \Delta^{+++}_{w,w^{op},w}\left(\Phi^{-1}\right).
$$
The $U$ part  decomposes as a series over $k \geq0$ 
when the exponential of the directed rectangle is expanded: $U= \sum_{k\geq 0} U_k/k!$. 
The very first term $U_0$ of this series contributes in $\Ztilde(\alpha_v)$ to 
$$
\chi^{-1} \left( Z(L) \circ U_0 \right) 
= \exp_\sqcup\left(\sum_{i=1}^{n} \strutgraph{i^-}{(i+1)^+}\quad \quad \right) 
= \Ztilde(\epsilon \otimes \Id_v),
$$
as is easily deduced from the fact that 
$$
C_w = Z\left(\Delta_w^+\left(\begin{array}{c}\capleft\\ {}^{(+-)} \end{array}\right)\right) 
\ \in \A\left( \Delta_w^+( \capleft) \right) = \A(\downarrow^w).
$$
So, to justify Claim \ref{claim:in_I},
it is enough to show that, for all $k>0$, $Z(L) \circ U_k$ belongs to $\chi\left(\mathcal{I}(n+1,n)\right)$.
Using Remark \ref{rem:top=bot}, we see that
$$
Z(L) \circ U_k = 
\left( \Id_{2n} \otimes Z\left(\Delta^+_{w^{op}}(\cupright)\right) \right) \circ  U_k 
\quad \hbox{mod } \ \chi(\mathcal{I}(n+1,n))
$$
by ``pushing'' directed rectangles arising from the ``doubling'' of the leftmost component of $\Phi^{-1}$ 
``towards the top'': This ``pushing'' means many application of the IHX and STU relations.
Finally, by ``pushing'' the directed rectangles in $U_k$ ``towards'' the struts labelled by $2^+,\dots,(n+1)^+$, 
we see that
$$
\left( \Id_{2n} \otimes Z\left(\Delta^+_{w^{op}}(\cupright)\right) \right) \circ  U_k  = 0
\quad \hbox{mod } \ \chi(\mathcal{I}(n+1,n)).
$$
We conclude that $Z(L) \circ U_k$ is trivial modulo $\chi(\mathcal{I}(n+1,n))$.
\end{proof}

\subsection{Hom-duals of the LMO functor}

Let us now explain how to recover from the LMO functor $\Ztilde$
the functor $\tau$ constructed by Le and the first author in \cite{CL1}.
For this, we need to derive from $\Ztilde$ other functors
whose target is the category of $\mathcal{A}(\varnothing)$-modules.

In general, for all integer $k\geq 0$, we denote by
$$
\ZZ_k : \LqCob \longrightarrow \mathcal{M}od_{\A(\varnothing)}
$$
the ``Hom-dual'' of the functor $\Ztilde$ defined by the object $k$ of $\tsA$, namely
$$
\ZZ_k := \tsA\left(k, \Ztilde(\centerdot) \right).
$$
Explicitly, for all non-associative word $w$ in the single letter $\bullet$, we set
$$
\ZZ_k(w) := \tsA(k,|w|)
$$
and, for all $M\in \LqCob(w,v)$, we set
$$
\ZZ_k(M) := \left( \tsA(k,|w|) \longrightarrow \tsA(k,|v|),\ x \longmapsto \Ztilde (M) \circ x \right).
$$
Similarly, we denote by
$$
\ZZ^k : \LqCob \longrightarrow \mathcal{M}od_{\A(\varnothing)}
$$
the contravariant ``Hom-dual'' of the functor $\Ztilde$ defined by the object $k$ of $\tsA$, namely
$$
\ZZ^k := \tsA\left(\Ztilde(\centerdot),k \right).
$$
In the closed case, there are also functors
$$
\ZZ_k : \widehat{\LqCob} \longrightarrow \mathcal{M}od_{\A(\varnothing)}
\quad \hbox{ and } \quad
\ZZ^k : \widehat{\LqCob} \longrightarrow \mathcal{M}od_{\A(\varnothing)}, 
$$
which are covariant and contravariant respectively, and are defined by 
$$
\ZZ_k := \widehat{\tsA}\left(k, \Ztilde(\centerdot) \right)
\quad \hbox{ and } \quad
\ZZ^k := \widehat{\tsA}\left(\Ztilde(\centerdot),k \right).
$$
All those functors  may be useful in the study of the categories $\LCob$ and $\widehat{\LCob}$, 
but we are particularly interested in the functors $\ZZ_0$ and $\ZZ^0$:

\begin{proposition} The following diagram is commutative:
$$
\xymatrix{
{\LqCob} \ar[r]^-{\ZZ_0} \ar@{->>}[d]_-{\widehat{\centerdot}} 
& {\mathcal{M}od_{\mathcal{A}(\varnothing)}}.\\
{\widehat{\LqCob}} \ar@{->}[ru]_-{\ZZ_0}
}
$$
A similar statement holds for the contravariant functor $\ZZ^0$.
\end{proposition}

\begin{proof} Remark \ref{rem:top=bot} implies that, for all integer $g\geq 0$,
$$
\widehat{\tsA}(0,g) = \tsA(0,g) = \A( \set{g}^-).
$$
So, the two functors $\ZZ_0$'s agree at the level of objects. 
For all $k\geq 0$ and for all $M \in \LqCob(w,v)$, the diagram
$$
\xymatrix{
{\tsA(k,|w|)} \ar@{->>}[d] \ar[rr]^-{\ZZ_k(M)} & & {\tsA(k,|v|)} \ar@{->>}[d]\\
{\widehat{\tsA}(k,|w|)}  \ar[rr]_-{\ZZ_k(\widehat{M})} & & {\widehat{\tsA}(k,|v|)} 
}
$$
is commutative since $\Ztilde(\widehat{M}) = \widehat{\Ztilde(M)}$, by definition.
Taking $k=0$ shows that the two functors $\ZZ_0$'s agree at the level of morphisms.
A similar argument works for $\ZZ^0$.
\end{proof}

We now aim at showing that the functor $\ZZ^0$ determines the 
functor $\tau$ defined in \cite{CL1}. 
The later is defined on a restricted class of Lagrangian cobordisms:

\begin{definition}
\label{def:LLCob}
A cobordism $(M,m)$ from $F_{g_+}$ to $F_{g_-}$ is \emph{doubly Lagrangian} 
if the following two conditions are satisfied:
\begin{enumerate}
\item $H_1(M) = m_{-,*}(A_{g_-}) + m_{+,*}(B_{g_+})$,
\item $m_{+,*}(A_{g_+}) \subset m_{-,*}(A_{g_-})$ and $m_{-,*}(B_{g_-}) \subset m_{+,*}(B_{g_+})$
as subgroups of $H_1(M)$.
\end{enumerate}
\end{definition}

\noindent
Doubly Lagrangian cobordisms define a monoidal subcategory of $\LCob$, 
which we denote by $\dLCob$. A subcategory $\widehat{\dLCob}$
of $\widehat{\LCob}$ is defined in a similar way. 
By ``reversion'' of cobordisms (i.e$.$ exchange of top and bottom, and reversion of the orientation),
$\dLCob$ and $\widehat{\dLCob}$ are isomorphic to their respective opposite categories.
The opposite of $\widehat{\dLCob}$ is
the category $\mathfrak{Z}$ of ``semi-Lagrangian cobordisms'' introduced in \cite{CL2}.
So, we wish to compare our contravariant functor
$$
\ZZ^0|_{\widehat{\dLqCob}}: \widehat{\dLqCob} \longrightarrow \mathcal{M}od_{\mathcal{A}(\varnothing)}
$$
to the covariant functor 
$$
\tau: \mathfrak{Z} \longrightarrow \mathcal{M}od_{\mathcal{A}(\varnothing)}
$$ 
constructed in \cite{CL1}. 

Thus, we consider an arbitrary cobordism $\widehat{M} \in \widehat{\dLCob}(g,f)$ 
and we equip it at the top and at the bottom with the left-handed parenthesizing 
$\left(\cdots \left( \left(\bullet \bullet\right) \bullet \right) \cdots \bullet\right)$.
On the one hand, observe that $\Ztilde(\widehat{M}) \in \widehat{\tsA}\left(g,f\right)$ 
is not only $\set{g}^+$-substantial but also $\set{f}^-$-substantial. So,
the map $\ZZ^0(\widehat{M}): \AY(\set{f}^+) \to \AY(\set{g}^+)$ being defined by 
$$
x \longmapsto 
\left\langle (x/i^+ \mapsto i^*) ,
 \left(\Ztilde(\widehat{M}) / i^-\mapsto i^*\right) \right\rangle_{\set{f}^*},
$$
it extends by the same formula to a map $\A(\set{f}^+) \to \A(\set{g}^+)$. 
On the other hand, there is the operator 
$\tau(\widehat{M}): \A(\cupright^{\set{f}}) \to \A(\cupright^{\set{g}})$ 
associated to $\widehat{M} \in \widehat{\dLCob}(g,f)=\mathfrak{Z}(f,g)$ in \cite{CL1}. 
We deduce from the definitions that the following diagram is commutative:
$$
\xymatrix{
{\A\left(\set{f}^+\right)}  \ar[d]_-{\ZZ^0(\widehat{M})} &  
 {\A\left(\set{f}^+\right)} 
\ar[l]_-{ \langle \centerdot\ , \mathsf{T}_f \rangle}^\simeq 
\ar[r]^-{\chi}_-\simeq  &
{\A\left(\cupright^{\set{f}}\right)} \ar[rr]_-{\simeq}^-{ \sqrt{\nu}^{\otimes \set{f}}\ \sharp_{\set{f}}\ \centerdot} &  &
{\A\left(\cupright^{\set{f}}\right)} \ar[d]^-{\tau\left(\widehat{M}\right)} \\
{\A\left(\set{g}^+\right)} & {\A\left(\set{g}^+\right)} 
\ar[l]^-{\langle \centerdot\ , \mathsf{T}_g \rangle}_\simeq 
\ar[r]_-{\chi}^-\simeq  & 
{\A\left(\cupright^{\set{g}}\right)} \ar[rr]^-{\simeq}_-{\sqrt{\nu}^{\otimes \set{g}}\ \sharp_{\set{g}}\ \centerdot} &  &
{\A\left(\cupright^{\set{g}}\right)}
}
$$
where $\langle \centerdot\ , \mathsf{T}_f \rangle$ denotes
$\left \langle (\centerdot/i^+\mapsto i^*), (\mathsf{T}_f/i^-\mapsto i^*) \right\rangle_{\set{f}^*}$.

\begin{remark}
The construction of the functor $\tau$ in \cite{CL1} needs
to extend the Kontsevich integral to trivalent graphs embedded in $S^3$.
Our approach avoids this extension.
\end{remark}

\vspace{0.5cm}

\section{The LMO functor as universal finite-type invariant}

In this section, we prove the universality of the LMO functor 
among finite-type invariants of Lagrangian cobordisms. 
This is achieved by clasper calculus.

\subsection{Clasper calculus}

\label{subsec:claspers}

This theory has been introduced independently by Goussarov and the second author --
see \cite{Habiro_claspers,GGP,Goussarov_clovers,Ohtsuki}.
For the reader's convenience, we recall some definitions following \cite{Habiro_claspers}. 

First of all, recall that clasper calculus applies to \emph{$3$-manifolds with tangles}, namely pairs $(M,\gamma)$
where $M$ is a compact oriented $3$-manifold whose boundary (if any) 
is identified with an abstract surface and where $\gamma\subset M$ is a framed oriented tangle whose
boundary (if any) corresponds to marked points on that surface.

A \emph{graph clasper} is a compact surface $G$ embedded into the interior of $M\setminus \gamma$ and decomposed in a certain way.
Precisely, $G$ should be decomposed between disks (the \emph{nodes}) 
and annuli (the \emph{leaves}) which are connected by means of bands (the \emph{edges}) as follows:
Each leaf should be connected to a single node and to no leaf, and each node should be connected to exactly three nodes or leaves.
Thus, if non-empty, each connected component of $G$ should contain at least one node.

The \emph{shape} of $G$ is the uni-trivalent graph which encodes the incidence rules of leaves and nodes in $G$,
where univalent and trivalent vertices correspond to leaves and nodes respectively.
The \emph{internal degree} of $G$ is the internal degree of its shape, i.e$.$ the number of its nodes.

\begin{example} If a graph clasper $G$ is shaped like a ``Y'', then it is called  a \emph{$Y$-graph}. 
It has one single node which is connected to three leaves by three edges. 
With the blackboard framing convention, it is pictured as
\begin{figure}[h!]
\begin{center}
\includegraphics[height=3cm,width=3cm]{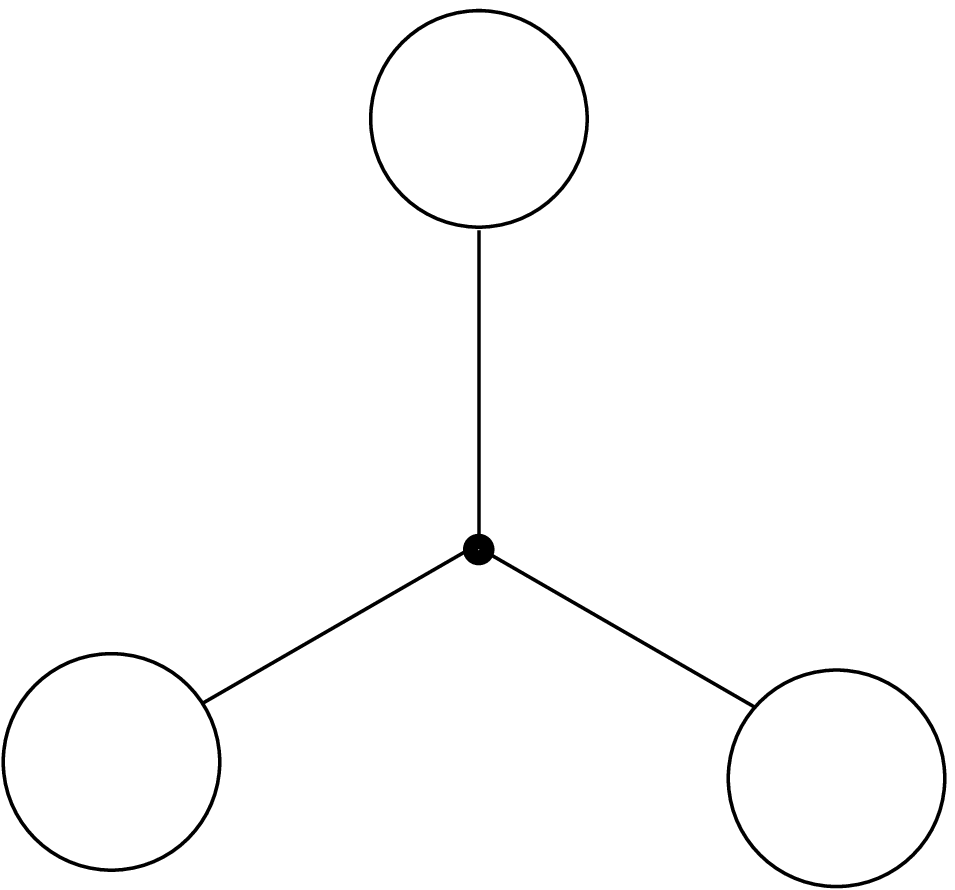}.
\end{center}
\end{figure} 
\end{example}

\begin{remark}
Graph claspers, as defined here, correspond to ``allowable graph claspers'' 
in the terminology of \cite{Habiro_claspers}.
\end{remark}

A graph clasper carries surgery instructions.
First, suppose that $G$ is a $Y$-graph in $(M,\gamma)$ and let N$(G)$ be its
regular neighborhood in $M\setminus \gamma$. This is a genus $3$ handlebody, 
where one can perform the surgery along the following six-component framed link $B$: 
\centerline{\relabelbox \small
\epsfxsize 2truein \epsfbox{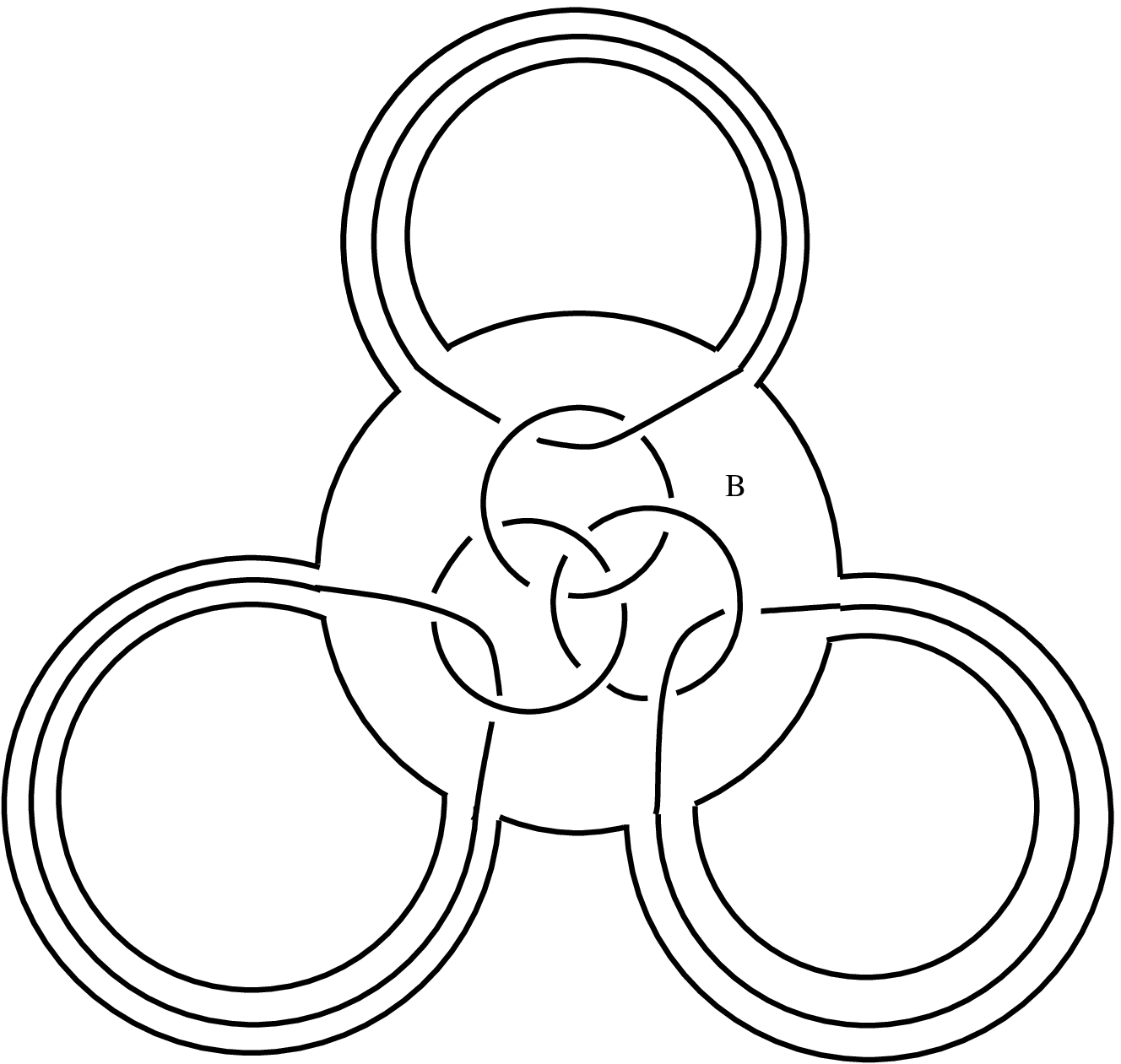}
\adjustrelabel <0cm,-0cm> {B}{$B$}
\endrelabelbox}
The manifold with tangle obtained from $(M,\gamma)$ by \emph{surgery along} $G$ is $(M_G,\gamma_G)$ where
$$
M_G := \left(M\setminus \hbox{int}\ \hbox{N}(G)\right) \cup \hbox{N}(G)_B
$$
and $\gamma_G$ is the trace in $M_G$ of $\gamma \subset \left( M\setminus \hbox{int}\ \hbox{N}(G)\right)$.
\begin{example}
\label{ex:Y}
If $M$ is the standard cube $[-1,1]^3$ and $\gamma$ is  the trivial $3$-component top tangle,
then we can consider the following $Y$-graph $G$:
$$
\includegraphics[width=3cm,height=3cm]{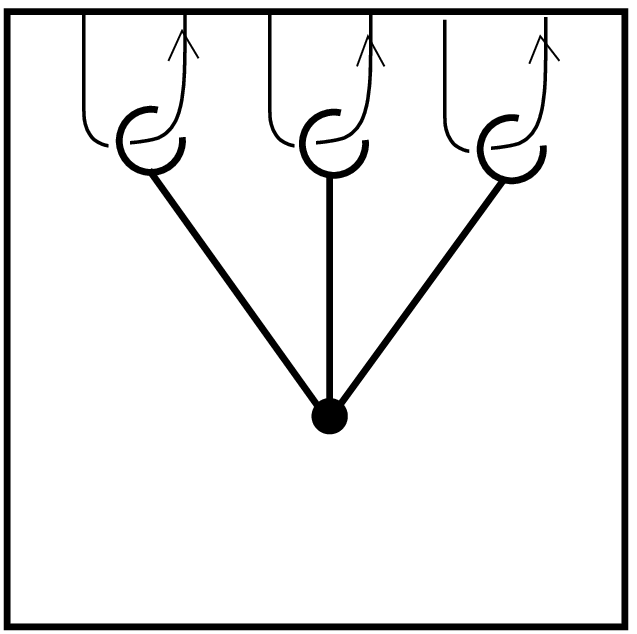}.
$$
Then, surgery along $G$ produces the top tangle $Y$ shown on Figure \ref{fig:generators}.
\end{example}
In general, suppose that $G$ is a graph clasper in $(M,\gamma)$ of internal degree $k$. 
By applying the ``fission'' rule
\begin{center}
\includegraphics[height=1.5cm,width=8cm]{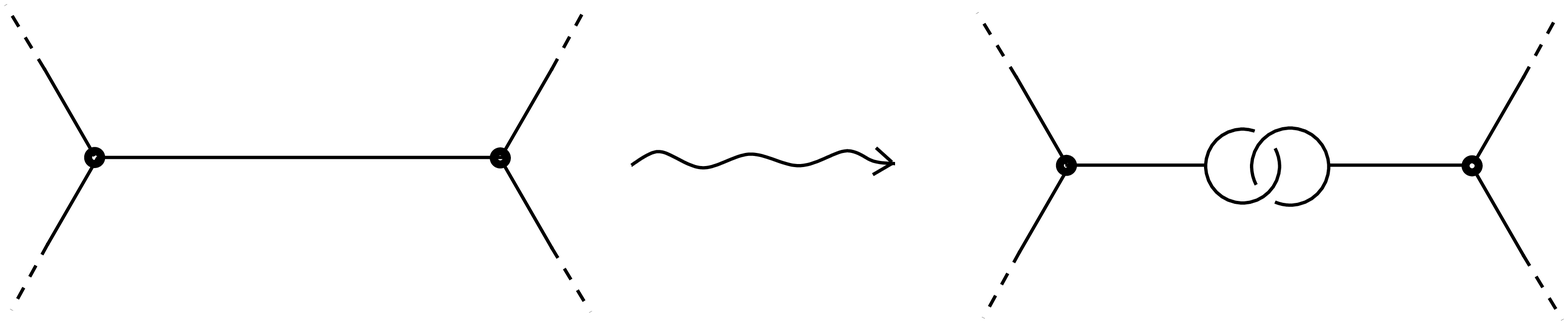} 
\end{center}
as many times as necessary, 
$G$ can be transformed to a disjoint union $Y(G)$ of $k$ $Y$-graphs in $M$.
The $3$-manifold with tangle obtained from $(M,\gamma)$ by \emph{surgery along} $G$
is the pair $(M,\gamma)$ surgered along each component of $Y(G)$.

\begin{definition}
\label{def:Y_k}
Let $k>0$ be an integer.  The \emph{$Y_k$-equivalence} is the equivalence relation 
among $3$-manifolds with tangles generated 
by surgeries along connected graph claspers of internal degree $k$
(and by homeomorphisms that preserve the orientations and the boundary parameterizations).
\end{definition}

\begin{remark}
The $Y_k$-equivalence is called ``$A_k$-equivalence'' in \cite{Habiro_claspers}.
\end{remark}

With this terminology recalled, let us characterize the $Y_1$-equivalence relation
in the case of bottom-top tangles in homology cubes:

\begin{theorem}[Matveev \cite{Matveev}, Murakami--Nakanishi \cite{MN}]
\label{th:Y1}
Two bottom-top tangles in homology cubes $(B,\gamma)$ and $(B',\gamma')$
are $Y_1$-equivalent if, and only if, they have the same linking matrix.
\end{theorem}

\begin{proof}
It is well-known  (and easily checked from the definition) that
surgery along a $Y$-graph preserves any information of homological-type, e.g$.$ linking numbers.
Assume that $(B,\gamma)$ has the same linking matrix as $(B',\gamma')$. 
By a result due to Matveev \cite{Matveev}, 
any homology sphere is $Y_1$-equivalent to $S^3$. So, $(B,\gamma)$ and $(B',\gamma')$
are $Y_1$-equivalent to some  bottom-top tangles $([-1,1]^3,\gamma_0)$ and $([-1,1]^3,\gamma'_0)$ respectively. 
Since $\gamma_0$ and $\gamma_0'$ have the same linking matrix in $[-1,1]^3$,
a theorem due to H. Murakami and Nakanishi \cite{MN} guarantees that $\gamma_0$ and $\gamma'_0$
are connected one to the other by isotopies and ``$\Delta$-moves.''
Those moves being realizable by surgeries along $Y$-graphs,
we deduce that $([-1,1]^3,\gamma_0)$ is $Y_1$-equivalent to $([-1,1]^3,\gamma_0')$.
\end{proof}

\begin{corollary}
\label{cor:Y1_L_cobordisms}
Two Lagrangian cobordisms $(M,m)$ and $(M',m')$ from $F_g$ to $F_f$ are $Y_1$-equivalent
if, and only if, there exists an isomorphism $\psi: H_1(M) \to H_1(M')$ 
such that the following diagram commutes:
$$
\xymatrix{
& {H_1(M)} \ar[dd]^-\simeq_-\psi & \\
{H_1(F_f)} \ar[ru]^-{m_{-,*}} \ar[rd]_-{m'_{-,*}} & & 
{H_1(F_g).} \ar[lu]_-{m_{+,*}} \ar[dl]^-{m'_{+,*}} \\
&H_1(M')&
}
$$
\end{corollary}

\begin{proof}
Let $(B,\gamma)$ and $(B',\gamma')$ be the  
bottom-top tangles in homology cubes corresponding to $M$ and $M'$ respectively.
Then, by Theorem \ref{th:Y1}, $M$ is $Y_1$-equivalent to $M'$ if and only if 
$\Lk_B(\gamma)=\Lk_{B'}(\gamma')$. This amounts to say that
$$
\left\{\begin{array}{ll}
\forall i \in\{1,\dots,f\}, & m_{-,*}(\beta_i) = m'_{-,*}(\beta_i) \in H_1(M) = H_1(M')\\
\forall i \in\{1,\dots,g\}, &m_{+,*}(\alpha_i) = m'_{+,*}(\alpha_i) \in H_1(M) = H_1(M')
\end{array}\right.
$$
where the identification $H_1(M) = H_1(M')$ is defined
by $m_{-,*}(\alpha_i) = m'_{-,*}(\alpha_i)$ (for all $i=1,\dots,f$)
and $m_{+,*}(\beta_i) = m'_{+,*}(\beta_i)$ (for all $i=1,\dots,g$).
\end{proof}

\begin{remark}
\label{rem:generation}
According to Corollary \ref{cor:Y1_L_cobordisms}, 
any Lagrangian cobordism is $Y_1$-equivalent to a special Lagrangian cobordism.
It follows from Example \ref{ex:Y} that
the monoidal category $\LCob$ is generated by the subcategory $\sLCob$
together with the morphism $Y: 3 \to 0$.
\end{remark}

\subsection{Universality of the LMO functor}

\label{subsec:universality}

We now prove that the LMO functor is universal among rational-valued finite-type invariants.
Let us start by recalling what a ``finite-type invariant'' is.

In general, we fix a $Y_1$-equivalence class $\mathcal{M}^0$ of $3$-manifolds, 
and  we denote by $\Q\cdot \mathcal{M}^0$ the $\Q$-vector space generated by its elements. 
Consider the following subspace:
$$
\mathcal{F}_{k}\left(\mathcal{M}^0\right)
:= \left\langle\ [M,G]\ \left\vert\ 
 M \in \mathcal{M}^0 \hbox{ and $G$ is a graph clasper in $M$ of i-degree $k$} 
\right.\ \right\rangle_{\Q}.
$$
Here, $[M,G]$ denotes the \emph{surgery bracket} defined by
$$
[M,G] := \sum_{G'\subset G} (-1)^{|G'|} \cdot M_{G'} \ \in \Q\cdot \mathcal{M}^0
$$ 
where the sum is taken over all subsets $G'\subset G$ of the set of connected components of $G$.
The \emph{$Y$-filtration} is the series
\begin{equation}
\label{eq:Y-filtration}
\Q\cdot \mathcal{M}^0 \geq \mathcal{F}_1\left(\mathcal{M}^0\right)
\geq \mathcal{F}_2\left(\mathcal{M}^0\right) \geq \mathcal{F}_3\left(\mathcal{M}^0\right) \geq \cdots
\end{equation}
whose associated graded vector space is denoted by
$$
\mathcal{G}\left(\mathcal{M}^0\right) 
= \prod_{i\geq 0} \mathcal{G}_i\left(\mathcal{M}^0\right) 
:= \prod_{i\geq 0} \frac{\mathcal{F}_i\left(\mathcal{M}^0\right)}{\mathcal{F}_{i+1}\left(\mathcal{M}^0\right)}.
$$

\begin{definition}
\label{def:fti}
Let $V$ be a $\Q$-vector space. A map $f:\mathcal{M}^0 \to V$ is a \emph{finite-type invariant} 
of \emph{degree} at most $d$ if its linear extension to $\Q\cdot \mathcal{M}^0$ 
vanishes on $\mathcal{F}_{d+1}\left(\mathcal{M}^0\right)$.
\end{definition}

We now come back to the case of Lagrangian cobordisms.
Lemma \ref{lem:linking_matrix_L_cobordisms} together with Corollary \ref{cor:Y1_L_cobordisms}
tell us that the $Y_1$-equivalence class of a Lagrangian cobordism $M$ is encoded by the $s$-reduction
of $\Ztilde(M)$. We are now going to see that rational finite-type invariants of positive degree 
correspond to the $Y$-reduction of $\Ztilde(M)$.
For this, we need an intermediate  result, 
which is proved by clasper calculus \cite{GGP,Habiro_claspers,Goussarov_clovers,Ohtsuki}.

\begin{theorem}[Garoufalidis \cite{Garoufalidis}]
\label{th:surgery}
Let $\mathcal{M}^0$ be a $Y_1$-equivalence class of Lagrangian cobordisms from $F_g$ to $F_f$.
Then, there exists a surjective graded linear map
$$
\mathfrak{S}: \AY({\set{g}}^+ \cup {\set{f}}^-) \longrightarrow \mathcal{G}\left(\mathcal{M}^0\right)
$$
which realizes each Jacobi diagram as a graph clasper in a representative of $\mathcal{M}^0$, 
and performs the surgery bracket along it.
\end{theorem}

\noindent
Let us specify what is meant by the ``realization'' of a Jacobi diagram.
For this, we first choose a representative $(M^0,m^0)$ in the class $\mathcal{M}^0$. Let 
$$
D\in \AY({\set{g}}^+ \cup {\set{f}}^-)
$$ 
be a Jacobi diagram without strut of i-degree $i$.
As a vertex-oriented graph, $D$ can be thickened in a unique way to
an oriented surface (vertices are thickened to disks, and
edges to bands). Cut a smaller disk in the interior of each disk
that has been produced from an external vertex of $D$.
One gets an oriented compact surface $S(D)$,
decomposed between disks, bands and annuli (corresponding to internal
vertices, edges and external vertices of $D$ respectively).
Use the induced orientation on $\partial S(D)$ to orient the cores of the annuli.
Embed $S(D)$ into the interior of $M^0$ in such a way that
each annulus of $S(D)$, as a framed oriented knot, is a push-off into the interior of 
$M^0$ of the framed curve $m^0_-(\alpha_j) \subset \partial M^0$ 
(respectively $m^0_+(\beta_j)  \subset \partial M^0$) 
if the color of the corresponding external vertex of $D$ is $j^-$ (respectively $j^+$). 
We obtain a graph clasper $C(D) \subset M^0$ of shape $D$, which
we call a \emph{topological realization} of the Jacobi diagram $D$.
Then, we set 
$$
\mathfrak{S}(D):= \left\{ \left[M^0,C(D)\right]\right\} \ \in \mathcal{G}_i\left(\mathcal{M}^0\right).
$$
The map $\mathfrak{S}$ is well-defined, surjective and does not depend 
on the initial choice of $(M^0,m^0)$ in the class $\mathcal{M}^0$. 
(See \cite[Theorem 1]{Garoufalidis} for a proof.)\\

The universality of the LMO functor is formulated as follows:

\begin{theorem}
\label{th:universality}
Let  $\mathcal{M}^0$ be a $Y_1$-equivalence class of Lagrangian cobordisms from $F_g$ to $F_f$.
Let also $w$ and $v$ be non-associative words in the single letter $\bullet$ of length $g$ and $f$ respectively.
Then, the i-degree $i$ part of 
$$
\ZtildeY: \mathcal{M}^0 \longrightarrow \AY({\set{g}}^+ \cup {\set{f}}^-), 
\ M \longmapsto \ZtildeY\left(
\begin{array}{c} M  \hbox{\footnotesize equipped with the words}\\
\hbox{\footnotesize $w_t(M) = w$, $w_b(M)=v$}
\end{array}\right),
$$ 
is a finite-type invariant of degree $i$. Moreover, the induced map 
$$
\Gr\ \ZtildeY : \mathcal{G}\left(\mathcal{M}^0\right) \longrightarrow \AY({\set{g}}^+ \cup {\set{f}}^-)
$$
gives, up to an explicit sign, a right-inverse to the surgery map $\mathfrak{S}$ from Theorem \ref{th:surgery}. 
In particular, $\mathfrak{S}$ and $\Gr\ \ZtildeY$ are both isomorphisms.
\end{theorem}

\noindent
Of course, when $\mathcal{M}^0=\{[-1,1]^3\}$ is the class of homology cubes, we recover Le's result 
on the universality of the LMO invariant for homology spheres \cite{Le_universality,Le_Grenoble,Ohtsuki}.

\begin{proof}[Proof of Theorem \ref{th:universality}]
Let $M \in \mathcal{M}^0$ and let $C$ be a graph clasper in $M$. 
Let $i$ be the internal degree of $C$, and let $c$ be the number of its connected components.
By ``fission'', $C$ can be transformed to a disjoint union $Y(C) = C_1 \cup  \cdots \cup C_i$ of $i$ $Y$-graphs. Then, 
$$
[M,C] = (-1)^{i-c} \cdot [M,Y(C)].
$$
We can write $M\in \LqCob(w,v)$ as
$$
M = \left( \Id_v \otimes \left( \cdots \left( {\epsilon^{\otimes ((\bullet \bullet) \bullet)}}_1 
\otimes {\epsilon^{\otimes ((\bullet \bullet) \bullet)}}_2 \right) \otimes \cdots \cdots \otimes 
{\epsilon^{\otimes ((\bullet \bullet) \bullet)}}_{i} \right)\right) \circ R
$$
where $\Id_v \in \LqCob(v,v)$ is the identity of the object $v$, 
where ${\epsilon^{\otimes ((\bullet \bullet) \bullet)}}_j$ is 
a copy of $\epsilon^{\otimes ((\bullet \bullet) \bullet)}$ 
and corresponds to a regular neighborhood of the $Y$-graph $C_j$, and where $R$ belongs to 
$\LqCob\left(w,v \otimes ((\bullet \bullet) \bullet)^{\otimes i}\right)$. 
More precisely, we assume that the picture of $C_j$ in ${\epsilon^{\otimes 3}}_j$ is exactly
$$
\begin{array}{c}
\includegraphics[width=3cm,height=3cm]{Ygen.eps}
\end{array}
$$
so that  ${\epsilon^{\otimes ((\bullet \bullet) \bullet)}}_j$ surgered along $C_j$ 
is a copy of the generator $Y$ of $\LqCob$, 
which we denote by $Y_j$. Since the functor $\Ztilde$ preserves the tensor product, we deduce that
$$
\Ztilde\left([M,C] \right) = (-1)^{i+c} \cdot
\left( \Ztilde(\Id_v) \otimes \bigotimes_{j=1}^i
\left( \Ztilde\left({\epsilon^{\otimes ((\bullet \bullet) \bullet)}}_{j} \right) -  \Ztilde\left(Y_j\right)\right) 
\right) \circ \Ztilde(R).
$$
Proposition \ref{prop:low_degree} implies that
\begin{equation}
\label{eq:many_Ys}
\Ztilde\left([M,C] \right) = (-1)^{i+c} \cdot
\left( \left[\sum_{k=1}^f \strutgraph{k^-}{k^+} \right] \otimes \bigotimes_{j=1}^i
\left( \Ygraphtop{1^+}{2^+}{3^+} \right) 
\right) \circ \left[\Lk(R)/2\right] + (\hbox{i-filter} > i).
\end{equation}
In particular, $\Ztilde\left([M,C] \right)$ has i-filter at least $i$,
which implies that $\ZtildeY\left([M,C] \right)$ starts in i-degree $i$.
This proves that the i-degree $(i-1)$ part of $\ZtildeY$ is a finite-type invariant of degree $i-1$.

In order to prove the second statement, we must carry on our argument further.
We now assume that the graph clasper $C$ is the topological realization 
of a Jacobi diagram $D \in \AY\left(\set{g}^+ \cup \set{f}^- \right)$: 
$$
C =C(D)
$$ 
as described just after Theorem \ref{th:surgery}.
Let $(B,\gamma)$ be a bottom-top tangle presentation of $R$: 
Then, a bottom-top tangle presentation of $M$ is obtained  
by deleting the last $3i$ bottom components of $\gamma$, which correspond to the leaves of $Y(C)$
and are labeled with the finite set 
$$
\mathcal{L}:=\left\{ (f+1)^-, (f+2)^-, (f+3)^-,\dots, (f+3i-2)^-, (f+3i-1)^-, (f+3i)^- \right\}.
$$
A leaf of $Y(C)$ can be of two types: Either it was already a leaf of $C$, 
or it has bornt from the ``fission'' process in which case it has a twin. 
Thus, there is a distinguished subset $\mathcal{O}$ of $\mathcal{L}$
which correspond to ``old leaves'', 
and a distinguished symmetric subset $\mathcal{T}$ of $\mathcal{L} \times \mathcal{L}$ 
corresponding to ``twin leaves''.
To each $l\in \mathcal{O}$, we can associate an unique $c(l) \in \set{g}^+ \cup \set{f}^-$,
which is the color of the corresponding external vertex of $D$.
Then, the identity
\begin{equation}
\label{eq:strut_combinations}
\Lk(R) = \Lk(M) + 2 \cdot \sum_{l\in \mathcal{O}} \strutgraph{l}{c(l)} 
 - \sum_{(t_1,t_2) \in \mathcal{T}} \strutgraph{t_1}{t_2}
\end{equation}
between linear combinations of struts is satisfied.
The minus sign in this identity is seen by comparing the ``fission'' process $C\leadsto Y(C)$
with the ``realization'' process $D \leadsto C(D)$.
We inject (\ref{eq:strut_combinations}) into (\ref{eq:many_Ys}) 
and we use Lemma \ref{lem:composition_ts_diagrams} to obtain 
$$
\Ztilde\left([M,C] \right) = (-1)^{i+c+e} \cdot \left[\Lk(M)/2\right] \sqcup D + (\hbox{i-filter}>i)
$$
where $e$ is the number of internal edges of $C$. We conclude that
\begin{equation}
\label{eq:universality}
\left(\Gr\ \ZtildeY\right) \circ \mathfrak{S}(D) = (-1)^{i+c+e} \cdot D 
\end{equation}
where the explicit sign is determined by 
$$
\left\{
\begin{array}{l}
i = \hbox{internal degree of } D\\
c = \hbox{number of connected components of }D\\
e= \hbox{number of internal edges of }D.
\end{array}
\right.
$$
\end{proof}

\begin{remark}
Lescop has proved that the Kontsevich--Kuperberg--Thurston invariant satisfies a ``splitting formula''
with respect to $\Q$-homology handlebody replacements:  See \cite{Lescop} for a precise statement.
Using the same kind of arguments than those used in the proof of Theorem \ref{th:universality},
we can prove that the LMO invariant satisfies the same formula. 
\end{remark}

\begin{corollary}
Let $f,g \geq 0$ be integers.
For all connected Jacobi diagram $D \in \AY\left(\set{g}^+ \cup \set{f}^-\right)$ of internal degree $i$,
and for all non-associative words $v$ and $w$ of length $f$ and $g$ respectively,
there exists an $M\in \LqCob(w,v)$ such that
$$
\ZtildeY(M) = \varnothing + D + (\ideg  >i) \ \in \AY\left(\set{g}^+ \cup \set{f}^-\right).
$$
\end{corollary}

\begin{proof}
This is a direct application of Theorem \ref{th:universality} to the $Y_1$-equivalence class
of the cube with tunnels and handles $C^g_f = \eta^{\otimes f} \circ \epsilon^{ \otimes g}$.
\end{proof}

\vspace{0.5cm}

\section{The LMO homomorphism of homology cylinders}

In this final section, we apply the LMO functor $\Ztilde$ to the study of homology cylinders, 
whose definition we first recall:

\subsection{The monoid $\Cyl(F_g)$ of homology cylinders}

Homology cylinders have been introduced in \cite{Goussarov_note,Habiro_claspers}
by Goussarov and the second author, 
in connection with surgery equivalence relations and finite-type invariants.

\begin{definition}
A \emph{homology cylinder} over $F_g$ is a cobordism $(M,m)$ from $F_g$ to $F_g$ 
such that $m_{\pm,*}: H_1(F_g) \to H_1(M)$ is an isomorphism and $m_{+,*}=m_{-,*}$.
\end{definition}

\noindent
It is easily seen that the set $\Cyl(F_g)$ of homology cylinders over $F_g$
is a monoid for the composition law $\circ$ of cobordisms and that
$$
\Cyl(F_g) \subset \LCob(g,g).
$$
Note that, for all Lagrangian cobordism $M$ from $F_g$ to $F_g$,
$$
M\in \Cyl(F_g) \ \ \Longleftrightarrow \ \
\Lk(M) = \left( \begin{array}{cc}
0 & I_g^{+-} \\ I_g^{-+} & 0
\end{array} \right).
$$
It follows from Corollary \ref{cor:Y1_L_cobordisms} that $\Cyl(F_g)$ 
is the $Y_1$-equivalence class of the trivial cobordism $\Id_g = F_g \times [-1,1]$.

\begin{example}
Let $\I(F_g)$ denote the \emph{Torelli group} of the surface $F_g$, namely
the kernel of the canonical homomorphism $\M(F_g) \to \Aut(H_1(F_g)), h \mapsto h_*$.
Then, the mapping cylinder construction restricts to an inclusion of monoids
$$
\I(F_g) \hookrightarrow \Cyl(F_g).
$$
\end{example}

\subsection{The LMO homomorphism $\ZtildeY$}

We now restrict the $Y$-part of the LMO functor to the monoid of homology cylinders.
Precisely, we assign to each homology cylinder $M$ over $F_g$ the series of Jacobi diagrams
$$
\ZtildeY(M) \in \AY(\set{g}^+ \cup \set{g}^-)
$$ 
where $M$ is seen as a Lagrangian $q$-cobordism 
with $w_t(M)=w_b(M)= \left(\cdots \left( \left(\bullet \bullet\right) \bullet \right) \cdots \bullet\right).$
Then, Theorem \ref{th:LMO_functor} specializes to the following

\begin{corollary}
\label{cor:LMO_h_cylinders}
The LMO invariant of homology cylinders defines a monoid homomorphism
$$
\ZtildeY: \left(\Cyl(F_g),\circ\right) \to \left(\AY(\set{g}^+ \cup \set{g}^-),\star\right).
$$
\end{corollary}

\noindent
Here, the multiplication $\star$  is the product defined in Example \ref{ex:star}. 
Recall that the formula for any $D, E \in \AY(\set{g}^+ \cup \set{g}^-)$ is
$$
D \star E = \left\langle\ 
\left(D \left/ i^+ \mapsto i^+ + i^* \right. \right)\ ,\ 
\left(E \left/ i^- \mapsto i^- +i^* \right.\right)\ \right\rangle_{\set{g}^*}.
$$
This multiplication $\star$ essentially coincides with the operation defined by Garoufalidis
and Levine under the same symbol in \cite{GL}. See Remark \ref{rem:*=*} below.\\

\begin{figure}[h!]
\centerline{\relabelbox \small
\epsfxsize 3truein \epsfbox{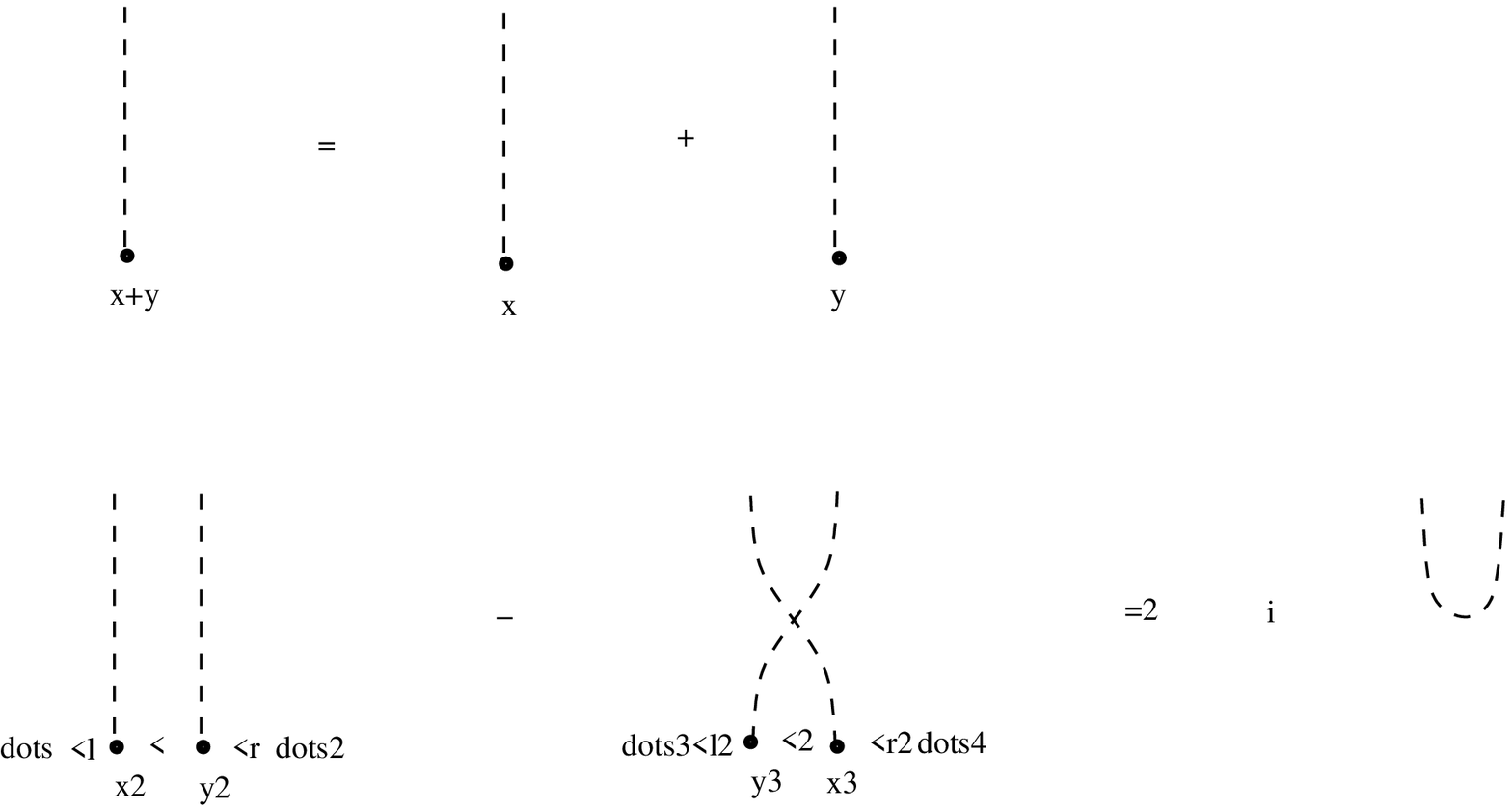}
\adjustrelabel <-0.05cm,-0.05cm> {x}{$x$}
\adjustrelabel <-0.05cm,-0.05cm> {y}{$y$}
\adjustrelabel <-0.35cm,-0.05cm> {x+y}{$x+y$}
\adjustrelabel <-0.1cm,-0.1cm> {x2}{$x$}
\adjustrelabel <-0.05cm,-0.05cm> {y2}{$y$}
\adjustrelabel <-0.05cm,-0.05cm> {x3}{$x$}
\adjustrelabel <-0.1cm,-0.05cm> {y3}{$y$}
\adjustrelabel <-0.1cm,-0.05cm> {<}{$<$}
\adjustrelabel <-0.05cm,-0.05cm> {<2}{$<$}
\adjustrelabel <-0.1cm,-0cm> {<l}{$<$}
\adjustrelabel <-0.05cm,-0cm> {<l2}{$<$}
\adjustrelabel <-0.05cm,-0.05cm> {<r}{$<$}
\adjustrelabel <-0.05cm,-0.05cm> {<r2}{$<$}
\adjustrelabel <-0.2cm,-0cm> {dots}{$\cdots$}
\adjustrelabel <0.05cm,-0.05cm> {dots2}{$\cdots$}
\adjustrelabel <-0.2cm,-0cm> {dots3}{$\cdots$}
\adjustrelabel <0.05cm,-0.05cm> {dots4}{$\cdots$}
\adjustrelabel <-0.2cm,-0cm> {-}{$-$}
\adjustrelabel <-0.1cm,-0.05cm> {+}{$+$}
\adjustrelabel <-0.3cm,-0cm> {i}{$(y \bullet x)$}
\adjustrelabel <0cm,-0cm> {=}{$=$}
\adjustrelabel <-0.2cm,-0cm> {=2}{$=$}
\endrelabelbox}
\caption{The \emph{multilinearity} and the \emph{STU-like} relations.
(Here, $x,y$ belong to $H_1(F_g)$ and $y\bullet x$ denotes their homological intersection.)}
\label{fig:relations}
\end{figure}

For future use, let us give an alternative description of the algebra $\left(\AY(\set{g}^+ \cup \set{g}^-),\star\right)$.
We consider the vector space
$$
\A\left(F_g\right) := 
\frac{\Q\cdot \left\{ \begin{array}{c} \hbox{vertex-oriented uni-trivalent graphs without strut and whose}\\
\hbox{external vertices are colored with } H_1(F_g) \hbox{ and are totally ordered} \end{array} \right\}}
{\hbox{AS, IHX, STU-like, multilinearity}}
$$
introduced by the second author in \cite[\S 8.5]{Habiro_claspers}. 
Here, the AS and IHX relations are as usual, while the multilinearity
and STU-like relations are defined on Figure \ref{fig:relations}. 
The space $\A\left(F_g\right)$ is graded by the internal degree of uni-trivalent graphs,
and its degree completion is denoted the same way.

For all diagrams $D,E \in \A(F_g)$, let $D \osqcup E$ denote the disjoint union of $D$ and $E$
where the external vertices are totally ordered in such a way that the vertices of $E$ 
are upper than those of $D$ and the given order of the vertices of $E$ (respectively $D$) is preserved.
By linear extension, one obtains an operation 
$\osqcup: \A(F_g) \times \A(F_g) \to \A(F_g)$, called the \emph{ordered disjoint union}.

\begin{lemma}
There is a $\Q$-algebra isomorphism
$$
\varphi: \left(\AY(\set{g}^+ \cup \set{g}^-),\star\right) \longrightarrow \left(\A(F_g), \osqcup\right)
$$
defined by declaring that ``each $i^-$-colored vertex should be lower than any $i^+$-colored vertex''
and by changing the colors of external vertices according to the rules 
$(i^- \mapsto \alpha_i)$ and $(i^+ \mapsto \beta_i)$.
\end{lemma}

\begin{proof}
First of all, the map $\varphi$ is well-defined for the following reason:
For any diagram $D \in \A(F_g)$ whose external vertices are colored by the
$\alpha_j$'s and $\beta_k$'s (rather than linear combinations of these),
only the order of the $\alpha_i$-colored vertices with respect 
to the $\beta_i$-colored vertices (for each $i=1,...,g$) is relevant.

Let $D \in \A(F_g)$ be a diagram whose external vertices are colored by the
$\alpha_j$'s and $\beta_k$'s: Let $n(D)$ be the number of couples $(a,b)$
where $a$ is an $\alpha_i$-colored vertex of $D$, $b$ is a $\beta_i$-colored vertex of $D$ and $a>b$.
If $n(D)=0$, then $D$ belongs to Im$(\varphi)$. If not, the STU-like relation allows one to write $D$ as $D'+D''$
with $n(D') = n(D)-1$ and $n(D'') < n(D)$. Thus, by induction on $n(D)$,
one oncludes that 
$D$ belongs to Im$(\varphi)$. It follows that $\varphi$ is surjective.

Let $D \in \A(F_g)$ be a diagram whose external vertices are colored by the
$\alpha_j$'s and $\beta_k$'s. Let $\psi(D)$ be 
$$
\left( \left. 
\begin{array}{c}
\hbox{sum of all ways of connecting some $\alpha_i$-colored}\\
\hbox{vertices of $D$ to some \emph{lower} $\beta_i$-colored vertices}
\end{array} 
\right/
\begin{array}{l} \alpha_i \mapsto i^- \\ \beta_i \mapsto i^+\end{array}
\right)_{\hbox{\& forget the order}}.
$$
Then, by linear extension, one gets a homomorphism $\psi: \A(F_g) \to \AY(\set{g}^+ \cup \set{g}^-)$.
Since $\psi \circ \varphi$ is the identity, we conclude that $\varphi$ is injective.

Let $D,E \in \AY(\set{g}^+ \cup \set{g}^-)$ be Jacobi diagrams. Then, $D\star E$ is the sum of all ways of gluing
some $i^-$-colored vertices of $E$ to some $i^+$-colored vertices of $D$.
Thus, $D \star E$ is equal to $\psi\left(\varphi(D) \osqcup \varphi(E)\right)$.
\end{proof}

Let $\Delta$ be the usual coproduct of $\AY(\set{g}^+ \cup \set{g}^-)$.
There is a similar coproduct on $\A(F_g)$ defined by the formula 
$$
\Delta(D) := \sum_{D = D' \sqcup D''} D' \otimes D''
$$
where, for all decomposition $D' \sqcup D''$ of $D$, 
the total ordering of the external vertices of $D'$ 
(respectively $D''$) in $D'\otimes D''$ is induced by the given order on $D$. 

Let $\varepsilon$ be the usual augmentation of $\AY(\set{g}^+ \cup \set{g}^-)$.
There is also a linear map $\varepsilon: \A(F_g) \to \Q$ defined by a similar formula,
namely $\Delta(D)= \delta_{D,\varnothing}$ for all Jacobi diagram $D$.

\begin{proposition}
\label{prop:Hopf_algebras}
$\left(\A(F_g), \osqcup, \varnothing, \Delta, \varepsilon \right)$
and $\left(\AY(\set{g}^+ \cup \set{g}^-), \star, \varnothing, \Delta, \varepsilon \right)$
are cocommutative graded Hopf algebras and are isomorphic via $\varphi$.
\end{proposition}

\begin{proof}
One easily sees that $\left(\A(F_g), \osqcup, \varnothing, \Delta, \varepsilon \right)$
is a Hopf algebra and that $\Delta \circ \varphi = (\varphi \otimes \varphi) \circ \Delta$.
Thus, $\left(\AY(\set{g}^+ \cup \set{g}^-), \star, \varnothing, \Delta, \varepsilon \right)$ 
is a Hopf algebra as well and is isomorphic to the previous one.
Those two algebras are graded by the internal degree.
\end{proof}

\subsection{The algebra dual to finite-type invariants of homology cylinders}

We now formulate the universality of the LMO homomorphism among 
rational-valued finite-type invariants of homology cylinders.
In other words, we apply the results from \S \ref{subsec:universality} to the $Y_1$-equivalence class
$$
\mathcal{M}^0 := \Cyl(F_g).
$$
In this case, the $Y$-filtration (\ref{eq:Y-filtration}) defines a filtration 
$$
\Q\left[\Cyl(F_g)\right] \geq \mathcal{F}_1\left( \Cyl(F_g)\right) 
\geq \mathcal{F}_2\left( \Cyl(F_g)\right) \geq \mathcal{F}_3\left( \Cyl(F_g)\right)  \geq \cdots
$$
of the \emph{monoid algebra} $\Q\left[\Cyl(F_g)\right]$ in the sense that
$$
\mathcal{F}_i\left( \Cyl(F_g)\right) \circ \mathcal{F}_j\left( \Cyl(F_g)\right) 
\subset \mathcal{F}_{i+j}\left( \Cyl(F_g)\right).
$$
The associated graded algebra  
$\mathcal{G}\left(\Cyl(F_g)\right)$
is called the \emph{algebra dual to finite-type invariants of homology cylinders}.

The LMO homomorphism from Corollary \ref{cor:LMO_h_cylinders} 
extends in the natural way to an algebra homomorphism
$$
\ZtildeY: \Q\left[\Cyl(F_g)\right] \longrightarrow \AY({\set{g}}^+ \cup {\set{g}}^-).
$$
Then, Theorem \ref{th:universality} specializes to the following

\begin{corollary}
\label{cor:universality_h_cylinders}
The algebra homomorphism 
$\ZtildeY: \Q\left[\Cyl(F_g)\right] \to \AY({\set{g}}^+ \cup {\set{g}}^-)$
sends the $Y$-filtration to the i-degree filtration, 
and induces an isomorphism at the graded level. Moreover, the inverse map of 
$$
\Gr\ \ZtildeY : \mathcal{G}\left(\Cyl(F_g)\right) \longrightarrow \AY({\set{g}}^+ \cup {\set{g}}^-)
$$ 
is, up to an explicit sign, the surgery map $\mathfrak{S}$ that realizes each Jacobi diagram as a graph clasper 
in $F_g \times[-1,1]$ and performs the surgery bracket along it.
\end{corollary}

\begin{remark}
\label{rem:*=*}
Let us mention how Corollary \ref{cor:universality_h_cylinders} connects to prior results: 
\begin{enumerate}
\item Together with Proposition \ref{prop:Hopf_algebras}, it proves that the algebra dual 
to finite-type invariants of homology cylinders is isomorphic to $\left(\A(F_g), \osqcup\right)$.
This had been announced by the second author in \cite[\S 8.5]{Habiro_claspers}.
\item 
Garoufalidis and Levine have shown in \cite{GL}, by means of clasper calculus, 
that the surgery map $\mathfrak{S}$ is surjective and sends (a sign-modification of) 
the multiplication $\star$ of Jacobi diagrams
to the composition $\circ$ of homology cylinders.
\item The fact that the domain and the codomain of $\Gr\ \ZtildeY$ are isomorphic 
as vector spaces has already been proved by Habegger \cite[Theorem 2.7]{Habegger}.
He used a different construction of the LMO invariant for homology cylinders,
but he did not consider the multiplicativity issue. See Remark \ref{rem:Habegger} below.
\end{enumerate}
\end{remark}

\subsection{The Lie algebra of homology cylinders}

The sequence of $Y_k$-equivalence relations ($k\geq 1$) filters the monoid $\Cyl(F_g)$ if one defines
$$
\Cyl_k\left(F_g\right)
:= \left\{M \in \Cyl\left(F_g\right) :
M \hbox{ is $Y_k$-equivalent to } F_g \times [-1,1] \right\}.
$$
According to \cite{Goussarov_note,Habiro_claspers}, the quotient monoid
$\Cyl_k\left(F_g\right)/Y_l$ is a group for all $l\geq k \geq 1$ and,
furthermore, the following inclusion holds  for all $k,k'\geq 1$ and $l\geq k+k'$:
$$
\left[\ \Cyl_{k}\left(F_g\right)/Y_l\ ,\
\Cyl_{k'}\left(F_g\right)/Y_l\ \right] \subset
\Cyl_{k+k'}\left(F_g\right)/Y_l.
$$
The graded \emph{Lie algebra of homology cylinders} is the graded vector space
$$
\overline{\Cyl}(F_g) := 
\prod_{i\geq 1} \frac{\Cyl_i(F_g)}{Y_{i+1}} \otimes \Q
$$
with the Lie bracket defined by taking commutators in the groups $\Cyl\left(F_g\right)/Y_l$ 's.
This Lie algebra has been introduced by the second author in \cite[\S 8.5]{Habiro_claspers}.

On the diagrammatic side, let $\mathcal{A}^{Y,c}({\set{g}}^+ \cup {\set{g}}^-)$
be the subspace of $\AY({\set{g}}^+ \cup {\set{g}}^-)$ spanned 
by non-empty $c$onnected Jacobi diagrams. 
This is the subspace of primitive elements and, so, 
is a Lie algebra with bracket $[x,y] := x \star y - y \star x$.

\begin{theorem}
\label{th:Lie_algebras}
The LMO homomorphism of homology cylinders induces a graded Lie algebra isomorphism
$$
\Gr\ \ZtildeY: \overline{\Cyl}(F_g) \longrightarrow \mathcal{A}^{Y,c}({\set{g}}^+ \cup {\set{g}}^-)
$$
that, for all $M\in \Cyl_i(F_g)$, sends $\{M\}\otimes 1$ 
to the i-degree $i$ part of $\ZtildeY(M)$.
\end{theorem}

\begin{proof}
Doing clasper calculus,
one can show in a way similar to the proof of Theorem \ref{th:surgery} 
that there is a surjective graded linear map
\begin{equation}
\label{eq:surg}
\mathcal{A}^{Y,c}({\set{g}}^+ \cup {\set{g}}^-) \to  \overline{\Cyl}(F_g)
\end{equation}
which ``realizes'' each connected Jacobi diagram as a graph clasper 
in $F_g \times[-1,1]$ and performs surgery along it. The linear map
\begin{equation}
\label{eq:alg}
\prod_{i\geq 1} \frac{\Cyl_i(F_g)}{Y_{i+1}} \otimes \Q \to 
\prod_{i\geq 1}\frac{\mathcal{F}_{i}\left(\Cyl(F_g)\right)}{\mathcal{F}_{i+1}\left(\Cyl(F_g)\right)},
\quad \{M\} \otimes 1 \mapsto \left\{M-(F_g \times[-1,1])\right\},
\end{equation}
preserves the Lie bracket because of the algebraic identity
$$
xyx^{-1}y^{-1} - 1 = [x-1,y-1] + \left((x-1)(y-1)-(y-1)(x-1)\right) (x^{-1}y^{-1}-1).
$$
Composing (\ref{eq:surg}) with (\ref{eq:alg}), one gets
the surgery map $\mathfrak{S}$ from Corollary \ref{cor:universality_h_cylinders} restricted to connected diagrams.
The conclusion follows since the inverse of $\mathfrak{S}$ is 
(up to an explicit sign) $\Gr\ \ZtildeY$.
\end{proof}

\begin{remark}
It should be emphasized that the LMO homomorphism $\ZtildeY$ strongly depends on the choice of
meridians and parallels $(\alpha,\beta)$ on the surface $F_g$.
Nevertheless, it can be proved by clasper calculus that the composite
$$
\xymatrix{
{\overline{\Cyl}(F_g)} \ar[r]^-{\Gr\ \ZtildeY}_-\simeq & 
{\mathcal{A}^{Y,c}({\set{g}}^+ \cup {\set{g}}^-)} \ar[r]^-\varphi_-\simeq &
{\A^c(F_g)}
}
$$
is independent of the choice of the basis $(\alpha,\beta)$ and, so, only depends on the surface $F_g$.
\end{remark}

We can deduce from the previous results that studying rational finite-type invariants of homology cylinders
is equivalent to studying the Lie algebra of homology cylinders: 

\begin{corollary}
\label{cor:Quillen}
The algebra dual to finite-type invariants of homology cylinders is 
the enveloping algebra of the Lie algebra of homology cylinders:
$$
\mathcal{G}\left( \Cyl(F_g) \right)
\simeq \hbox{U}\left(\overline{\Cyl}(F_g) \right).
$$
\end{corollary}

\begin{proof}
It appears from the proof of Theorem \ref{th:Lie_algebras} that the linear map
$$
\prod_{i\geq 1} \frac{\Cyl_i(F_g)}{Y_{i+1}} \otimes \Q \to 
\prod_{i\geq 1}\frac{\mathcal{F}_{i}\left(\Cyl(F_g)\right)}{\mathcal{F}_{i+1}\left(\Cyl(F_g)\right)},
\quad \{M\}\otimes 1 \mapsto \left\{M-(F_g \times[-1,1])\right\}
$$
is an embedding of Lie algebras, whose image is
$\mathfrak{S}\left(\mathcal{A}^{Y,c}({\set{g}}^+ \cup {\set{g}}^-)\right)$.

By Proposition \ref{prop:Hopf_algebras}, 
$\mathcal{A}^{Y}({\set{g}}^+ \cup {\set{g}}^-)$ is a co-commutative Hopf algebra isomorphic to $\A(F_g)$.
Using the STU-like relation, the latter is easily seen to be generated by its primitive elements: So is the former.
It follows from Milnor--Moore's theorem \cite{MM} that $\mathcal{A}^{Y}({\set{g}}^+ \cup {\set{g}}^-)$
is the enveloping algebra of $\mathcal{A}^{Y,c}({\set{g}}^+ \cup {\set{g}}^-)$.
\end{proof}

\begin{remark}
See \cite{Massuyeau} for an algebraic proof of Corollary \ref{cor:Quillen}.
\end{remark}

\subsection{The tree-reduction of the LMO invariant}

Following Habegger's approach \cite{Habegger}, which  we are going to recall, 
we now connect the tree-reduction of $\ZtildeY$ to Johnson homomorphisms.

\begin{figure}[h]
\centerline{\relabelbox \small
\epsfxsize 3truein \epsfbox{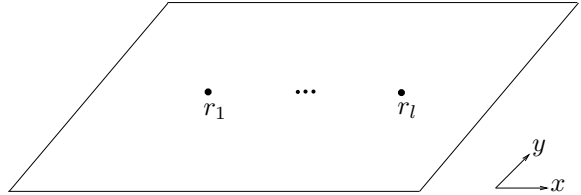}
\adjustrelabel <-0.05cm,-0.1cm> {r1}{$r_1$}
\adjustrelabel <-0.05cm,-0.1cm> {rn}{$r_l$}
\adjustrelabel <0cm,-0cm> {x}{$x$}
\adjustrelabel <0cm,-0cm> {y}{$y$}
\endrelabelbox}
\caption{The standard points $r_1,\dots,r_l$ on $[-1,1]^2$.}
\label{fig:single_points}
\end{figure}

For all $l\geq 1$, let $r_1,\dots,r_l$ be $l$ points
on $[-1,1]^2$ choosen uniformly in the $x$ direction as shown in Figure \ref{fig:single_points}.

\begin{definition}
By a \emph{string-link} on $l$ \emph{strands}, we mean an equivalence class of couples $(B,\gamma)$
where $B$ is a cobordism from $F_0$ to $F_0$ 
and $\gamma$ is an $l$-component framed oriented tangle in $B$ 
whose each component $\gamma_i$ runs from $r_i \times (-1)$ to $r_i \times 1$.
\end{definition}

Given two string-links $(C,\upsilon)$ and $(B,\gamma)$ on $l$ strands,
the  gluing of $\gamma$ to $\upsilon$ defines a string-link in $B \circ C$ (composition in $\Cob$).
Thus, the set of string-links on $l$ strands is a monoid 
and we are interested in the following submonoids:
$$
\begin{array}{ccl}
\mathcal{S}_n &:=& \left\{\hbox{string-link }(B,\gamma) : H_*(B)\simeq H_*([-1,1]^3) \right\}\\
\cup & &\\
\mathcal{S}^0_n &:=& \left\{\hbox{string-link }(B,\gamma) : H_*(B)\simeq H_*([-1,1]^3)
 \ \hbox{ and } \ \Lk_B(\gamma) = 0 \right\}.
\end{array}
$$

\begin{figure}[h]
\centerline{\relabelbox \small
\epsfxsize 4.5truein \epsfbox{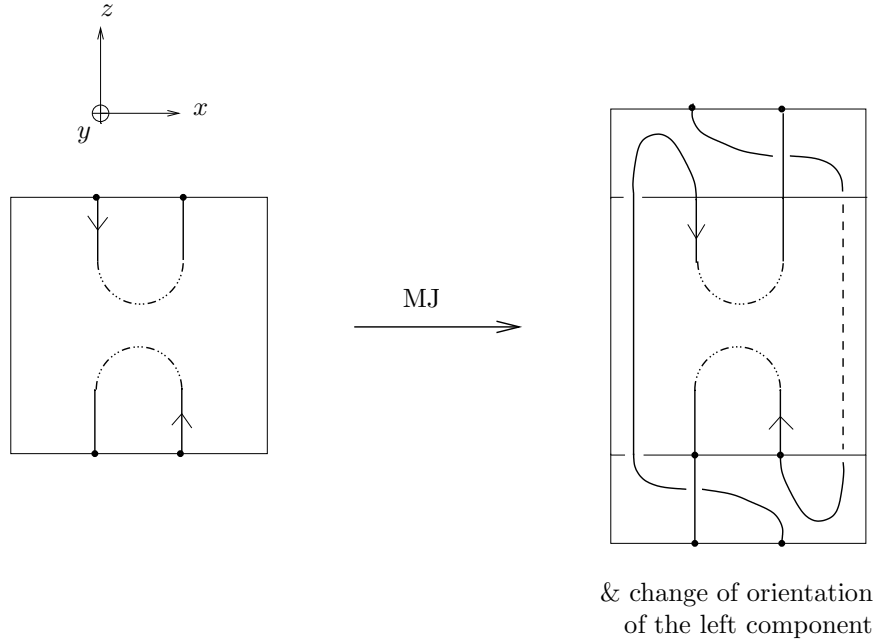}
\adjustrelabel <0cm,-0cm> {x}{$x$}
\adjustrelabel <-0.1cm,-0.1cm> {y}{$y$}
\adjustrelabel <0cm,-0cm> {z}{$z$}
\adjustrelabel <0cm,-0cm> {MJ}{$\hbox{MJ}$}
\adjustrelabel <-0.3cm,0.3cm> 
{inversion}{\begin{tabular}{r}\& change of orientation \\ of the left component \end{tabular}}
\endrelabelbox}
\caption{From a bottom-top tangle of type $(g,g)$
to a string-link on $2g$ strands. (Here, $g=1$; if $g>1$, 
the picture must be repeated $g$ times in the $x$ direction.)}
\label{fig:MJ}
\end{figure}

There are many ways to transform a bottom-top tangle of type $(g,f)$ 
to a string-link on $f+g$ components.
Figure \ref{fig:MJ} illustrates a way to do this in the special case when $f=g$.
This bijection transforms $\Id_g\in {\btT}(g,g)$ to
the trivial $(2g)$-component string-link in $[-1,1]^3$ and, so,
sends the $Y_1$-equivalence class of the former to
the $Y_1$-equivalence class of the latter.
According to \cite{Matveev,MN}, the $Y_1$-equivalence class of the trivial $n$-component
string-link is exactly $\mathcal{S}^0_n$. Hence a bijection
$$
\xymatrix{
\hbox{MJ}: {\Cyl\left(F_g\right)} \ar[r]^-\simeq &{\mathcal{S}^0_{2g}},
}
$$
which is called the \emph{Milnor--Johnson correspondence}.
This is essentially the  bijection introduced under the same name by Habegger in \cite{Habegger}. 

\begin{remark}
\label{rem:Habegger}
In his paper \cite{Habegger}, Habegger defines the LMO invariant of homology cylinders 
to be the composite map
$$
\xymatrix{
{\Cyl\left(F_g\right)} \ar[rr]_-\simeq^-{\hbox{MJ}} &&{\mathcal{S}^0_{2g}} 
\ar[rr]^-{\chi^{-1}Z} & & {\AY(\set{2g})} 
}
$$
where the letter $Z$ denotes the Kontsevich--LMO invariant of string-links in homology cubes
(as defined in \S \ref{subsec:Kontsevich-LMO}). 
But, since the map $\hbox{MJ}$ is not multiplicative, 
it does not seem  easy to understand how this invariant of homology cylinders 
behaves with respect to composition of cobordisms.
\end{remark}

\begin{theorem}[Habegger \cite{Habegger}]
\label{th:MJ}
The Milnor--Johnson correspondence relates
the Johnson homomorphisms of homology cylinders to the Milnor invariants
of string-links in homology cubes\footnote{See the proof below for the precise statement.}.
\end{theorem}

\begin{proof}

As mentioned previously, one could check that our version of the Milnor--Johnson correspondence
essentially coincides with Habegger's one (which is not obvious from the definitions) 
and appeal to his result \cite[Theorem 2.1]{Habegger}.
Instead, we prefer to repeat his arguments in our formulation of the MJ correspondence,
which gives the opportunity to review Milnor and Johnson invariants:\\

\noindent
\emph{Review of Johnson homomorphisms for homology cylinders \cite{GL}.} Let $\pi$ denote $\pi_1(F_g,*)$
which is free with basis $(\alpha_1,\dots,\alpha_g,\beta_1,\dots,\beta_g)$ (see Figure \ref{fig:surface}).

Let $(M,m) \in \Cyl(F_g)$.
Since $m_\pm:F_g \to M$ induces an isomorphism at the level of homology, it induces an isomorphism at the level
of each nilpotent quotient of the fundamental groups \cite{Stallings},
hence a monoid anti-homomorphism
\begin{equation}
\label{eq:monoid_homo}
\Cyl(F_g) \to \Aut\left(\pi/\pi_{n+1}\right), \ (M,m) \mapsto {m_{+,*}}^{-1}\circ m_{-,*}.
\end{equation}
Let $\Cyl(F_g)[n]$ be its kernel, whose elements are said to have \emph{Johnson filtration} $n$.
For instance, $\Cyl(F_g)[1]$ is $\Cyl(F_g)$.
There is the exact sequence of groups
$$
1 \to \Hom(\pi/\pi_2,\pi_{n+1}/\pi_{n+2}) \to \Aut(\pi/\pi_{n+2}) \to \Aut(\pi/\pi_{n+1})
$$
where  a group homomorphism $t: \pi/\pi_2 \to \pi_{n+1}/\pi_{n+2}$
goes to the automorphism of $\pi/\pi_{n+2}$ defined by $\{x\} \mapsto \{x \cdot t(\{x\})\}$.
So, the restriction of the map (\ref{eq:monoid_homo}) at the $(n+1)$-st level to $\Cyl(F_g)[n]$ 
gives a map
$$
\tau_n: \Cyl(F_g)[n] \to \Hom(\pi/\pi_2,\pi_{n+1}/\pi_{n+2})
\simeq  \pi/\pi_2 \otimes \pi_{n+1}/\pi_{n+2}
$$
where the last isomorphism is induced by the left-adjoint of the intersection pairing on $F_g$. 
The monoid homomomorphism $\tau_n$ is called the \emph{$n$-th Johnson homomorphism}, 
and is given by the formula
\begin{equation}
\label{eq:Johnson}
\tau_n(M,m) =
\sum_{i=1}^g \alpha_i \otimes m_{+,*}^{-1}\left((\beta_{i}^+)^{-1}\cdot \beta_i^-\right)
-\beta_i \otimes m_{+,*}^{-1}\left((\alpha_{i}^+)^{-1} \cdot \alpha_i^-\right)
\end{equation}
where $\beta_i^\pm := m_\pm (\beta_i)$ and $\alpha_i^\pm := m_\pm (\alpha_i)$.\\

\noindent
\emph{Review of Milnor invariants for string-links \cite{HL}.}
Consider the disk with $l$ holes $D_l := [-1,1]^2\setminus \hbox{N}(\{r_1,\dots,r_l\})$,
where $r_1,\dots,r_l$ are those points shown on Figure \ref{fig:single_points}. 
Let $\varpi$ denote $\pi_1(D_l,*)$, which is free with basis $(x_1,\dots,x_l)$ if $x_i$ is the loop winding around
$r_i$ in the trigonometric direction.

Let $(B,\gamma) \in \mathcal{S}_l$. The complement of $\gamma$ in $B$ is a compact oriented $3$-manifold $S$,
and the framing of $\gamma$  together with the given identification $b: \partial [-1,1]^3 \to \partial B$
define a homeomorphism $s: \partial(D_l \times [-1,1]) \to S$ onto the boundary of $S$: 
In other words, $(B,\gamma)$ can be regarded as a \emph{cobordism $(S,s)$ from $D_l$ to $D_l$}.
Since $B$ is a homology cube, this cobordism has the property that $s_\pm: D_l \to N$ is an isomorphism at the level of homology
and, so, at the level of each nilpotent quotient of the fundamental groups \cite{Stallings}.
Thus, one obtains a monoid anti-homomorphism
$$
\mathcal{S}_l \to \Aut\left(\varpi/\varpi_{n+1}\right), \ (B,\gamma)=(S,s) \mapsto {s_{+,*}}^{-1}\circ s_{-,*}.
$$
This is called the  $n$-th \emph{Artin representation}.
Let $\mathcal{S}_l[n]$ be its kernel, whose elements are said to have \emph{Milnor filtration} $n$.
For instance, $\mathcal{S}_l=\mathcal{S}_l[1]$ and $\mathcal{S}_l^0= \mathcal{S}_l[2]$.
Recall from \cite{HM} the well-known

\begin{claim}
\label{claim:free_commutators}
Let $\lambda$ belong to the free group $F(x_1,\dots,x_l)$ and let $i\in \{1,\dots,l\}$.
Then, $[\lambda,x_i] \in F_{n+1}$ if, and only if, $\lambda x_i^p \in F_n$ for some $p\in \Z$. 
\end{claim}

\noindent
Set $x_i^\pm := s_\pm(x_i)$. Observe that, modulo $\varpi_{n+1}$,
$$
x_i^{-1} \cdot s_{+,*}^{-1} \circ s_{-,*}(x_i) =
s_{+,*}^{-1}\left((x_i^+)^{-1} x_i^-\right)
$$
$$
= s_{+,*}^{-1}\left(l(\gamma_i)^{-1} (x_i^-)^{-1} l(\gamma_i) x_i^- \right)
= s_{+,*}^{-1}\left(\left[l(\gamma_i)^{-1},(x_i^-)^{-1}\right]\right)
$$
where $l(\gamma_i)$ is the oriented longitude of $\gamma_i$ defined by the framing.
One deduces from Claim \ref{claim:free_commutators} that, provided it has zero framing,
the string-link $(B,\gamma)$ belongs to $\mathcal{S}_l[n]$ if, and only if,
each of its longitudes is trivial modulo $\pi_1(S,*)_{n}$; one also deduces that
considering $x_i^{-1} \cdot s_{+,*}^{-1}\circ s_{-,*}(x_i)$ modulo $\varpi_{n+2}$ 
is equivalent to considering only $s_{+,*}^{-1} l(\gamma_i)$ modulo $\varpi_{n+1}$. 
So, one is led to the monoid homomorphism
$$
\mu_n: \mathcal{S}_l[n] \to \varpi/\varpi_2 \otimes \varpi_n/\varpi_{n+1}
$$
defined by the formula
\begin{equation}
\label{eq:Milnor}
\mu_n(B,\gamma) = \sum_{i=1}^l x_i \otimes s_{+,*}^{-1}(l(\gamma_i))
\end{equation}
and called the \emph{$n$-th Milnor invariant}.

\begin{innerremark} Usually, ``Milnor invariants''  refer to some integers that are defined as follows. 
Consider the Magnus expansion $\mu: \varpi \to \Z[[X_1,\dots,X_l]]$ 
defined by $\mu(x_i)=1+X_i$, where the indeterminates $X_i$'s do not commute.
It induces by truncation a map 
$$
\mu_{\leq n}: \varpi/\varpi_{n+1} \to ( \hbox{degree $\leq n$ part of }\Z[[X_1,\dots,X_l]]).
$$
Then, for all $i,j_1,\dots,j_n\in \{1,\dots,l\}$,
\emph{Milnor's invariant} $\mu(j_1,\dots,j_n;i)$ 
is the coefficient of $X_{j_1} \cdots X_{j_n}$ 
in $\mu_{\leq n}\left(s_{+,*}^{-1}(l(\gamma_i))\right)$. 
Thus, in the case when $(B,\gamma) \in \mathcal{S}_l[n]$, one has
$$
\mu_n(B,\gamma) = \sum_{i=1}^l x_i \otimes 
\left(\sum_{j_1,\dots,j_n} \mu(j_1,\dots,j_n;i)\cdot X_{j_1} \cdots X_{j_n}\right)
$$
where $\varpi_n/\varpi_{n+1}$ is identified by $\mu$ to a subgroup of the degree $n$ part of $\Z[[X_1,\dots,X_l]]$.
\end{innerremark}

\begin{claim}
\label{claim:MJ}
The Milnor--Johnson correspondence sends the Johnson filtration to the Milnor filtration (up to a shift of levels),
and the Johnson homomorphisms to the Milnor invariants. Precisely, we have 
$$
\xymatrix{
{\Cyl(F_{g})[n]} \ar[d]_-{\tau_n} \ar[r]^-{\hbox{MJ}}_\simeq &
{\mathcal{S}_{2g}[n+1] } \ar[d]^-{\mu_{n+1}} \\
{\pi_1/\pi_2 \otimes \pi_{n+1}/\pi_{n+2}} \ar[r]_-\simeq &
{\varpi_1/\varpi_2 \otimes \varpi_{n+1}/\varpi_{n+2}}
}
$$
where $n \geq 1$ and  the bottom isomorphism is induced by the identification
$\pi\simeq \varpi$ which sends $\alpha_i$ to $x_{2i-1}^{-1}$
and $\beta_i$ to $x_{2i}$.
\end{claim}

Now that the definitions are well set, the proof of the claim is straightforward,
and it is enough to prove its second statement.
Let $(M,m)$ be a homology cylinder of Johnson filtration $n$,
and let $(B,\gamma)$ be its image by $\hbox{MJ}$: We regard the latter as a cobordism $(S,s)$ from $D_{2g}$ to $D_{2g}$. 
We are asked to prove that
$$
\forall i \in \{1,\dots,g\} \quad
\left\{
\begin{array}{l}
m_{+,*}^{-1}\left((\beta_i^+)^{-1}\beta_{i}^-\right) = -s_{+,*}^{-1}(l(\gamma_{2i-1}))\\
m_{+,*}^{-1}\left((\alpha_i^+)^{-1}\alpha_{i}^-\right) = -s_{+,*}^{-1}(l(\gamma_{2i}))
\end{array}\right.
$$
in $\pi_{n+1}/\pi_{n+2} \simeq \varpi_{n+1}/\varpi_{n+2}$
(where $\pi$ is identified to $\varpi$ as in Claim \ref{claim:MJ}).
Since $\incl_* \circ m_{+,*}: \varpi \simeq \pi \to \pi_1(M,*) \to \pi_1(S,*)$
and $s_{+,*}: \varpi \to \pi_1(S,*)$ are identical at the $\varpi_1/\varpi_2$ level, 
they induce the same map at the $\varpi_{n+1}/\varpi_{n+2}$ level.
So, we wish to prove that
\begin{equation}
\label{eq:goal_MJ}
\forall i \in \{1,\dots,g\} \quad
\left\{
\begin{array}{l}
(\beta_{i}^-)^{-1} \beta_i^+ = l(\gamma_{2i-1}) 
\ \in \pi_1(S,*)_{n+1}/\pi_1(S,*)_{n+2} \\
(\alpha_{i}^-)^{-1} \alpha_i^+= l(\gamma_{2i}) \ \in \pi_1(S,*)_{n+1}/\pi_1(S,*)_{n+2}.
\end{array}\right.
\end{equation}
When $g=1$, the picture to have in mind is 
\begin{figure}[h!]
\centerline{\relabelbox \small
\epsfxsize 4truein \epsfbox{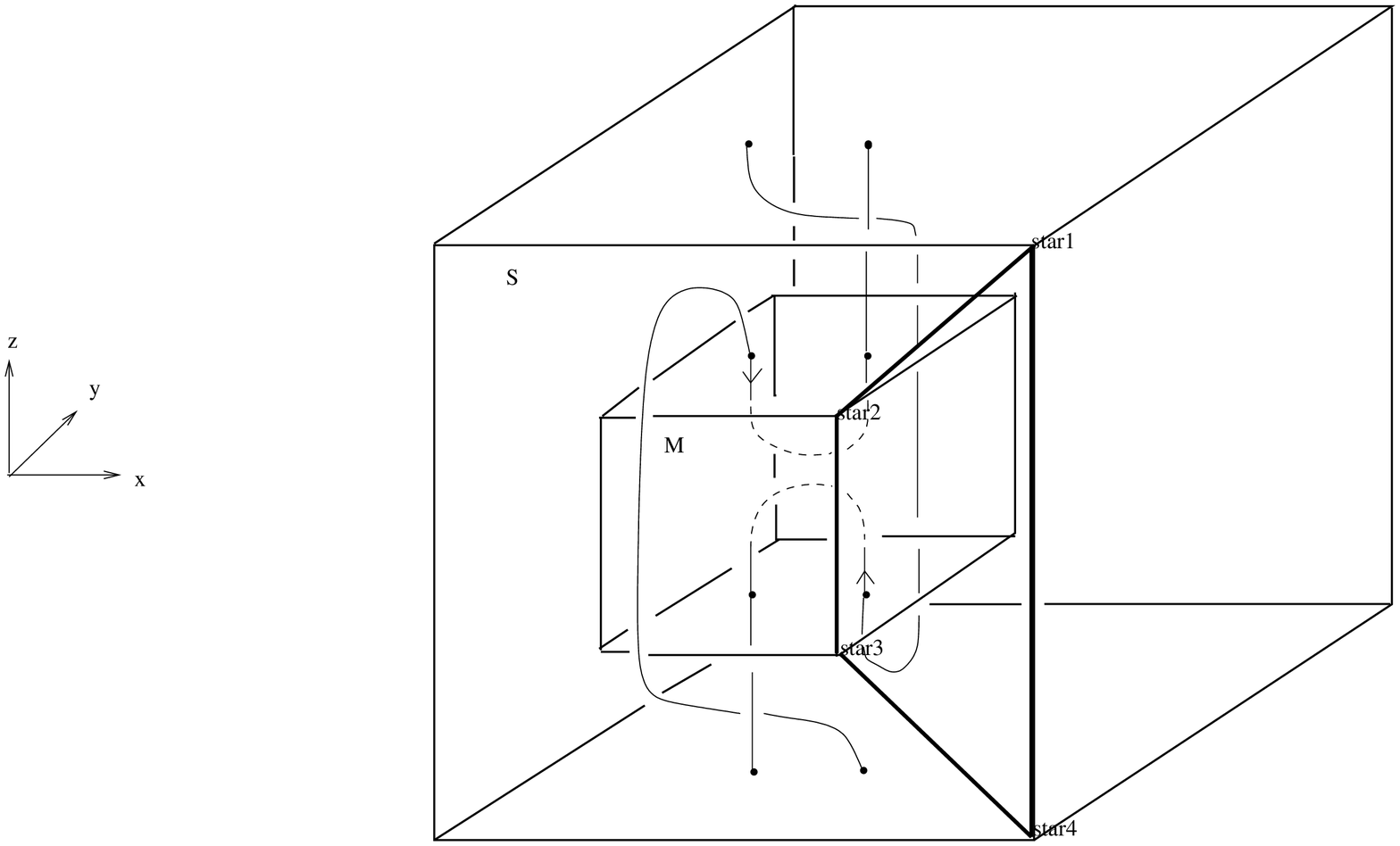}
\adjustrelabel <-0.1cm,-0.3cm> {star1}{\hbox{\LARGE *}}
\adjustrelabel <-0.15cm,-0.3cm> {star2}{\hbox{\LARGE *}}
\adjustrelabel <-0.15cm,-0.3cm> {star3}{\hbox{\LARGE *}}
\adjustrelabel <-0.15cm,-0.35cm> {star4}{\hbox{\LARGE *}}
\adjustrelabel <0cm,-0.2cm> {M}{$M$}
\adjustrelabel <0cm,-0.3cm> {S}{$S$}
\adjustrelabel <0cm,-0.cm> {x}{$x$}
\adjustrelabel <0cm,-0.cm> {y}{$y$}
\adjustrelabel <0cm,-0.cm> {z}{$z$}
\endrelabelbox}
\end{figure}

\noindent
where the $*$'s are base points for fundamental groups and are identified along
the bold lines. (One is allowed to do so since those bold lines bound a disk.)
When $g>1$, one should have in mind $g$ copies of that pair of cubes displayed
along the $x$ axis and glued along their $(y,z)$-faces and 
one should remember that the base points $*$'s are put on the rightmost pair of cubes. 

On that picture, one observes that
$$
l(\gamma_{2i-1})= \left( c_i (\beta_i^-)^{-1} c_i^{-1} \right) t_i  \left( c_i \beta_i^+ c_i^{-1} \right)
$$
where $t_i$ is a loop based on $\star$  which turns around the right $(y,z)$-face of the $i$-th cube,
and where $c_i$ is a loop based on $\star$ which grasps the $(g-i)$ rightmost copies of that strand of lowest $y$-coordinate
on the above picture. Thus, 
$t_g=1$ and $t_i= \left(x_{2i-1}^-\right)^{-1} \left(x_{2i}^-\right)^{-1} t_{i+1} x_{2i}^+ x_{2i-1}^+$.
Using the fact that $x_k^-= l(\gamma_k) x_k^+ l(\gamma_k)^{-1}$ 
and that $l(\gamma_k) \in \pi_1(S,*)_{n+1}$ commutes with everybody modulo $\pi_1(S,*)_{n+2}$,
one checks by decreasing induction on $i$ that $t_i \in \pi_1(S,*)_{n+2}$.
Since $(\beta_i^-)^{-1} \beta_i^+ \in \pi_1(S,*)_{n+1}$, this proves the first part of (\ref{eq:goal_MJ}).

On the same picture, one sees that
$$
l(\gamma_{2i}) = \left( c_i (\beta_i^-)^{-1} (\alpha_i^-)^{-1} \beta_i^- c_i^{-1} \right) 
\left( c_i (\beta_i^+)^{-1} \alpha_i^+ \beta_i^+ c_i^{-1} \right).
$$
Since $\beta_i^- (\beta_i^+)^{-1} \in \pi_1(S,*)_{n+1}$ and 
$(\alpha_i^-)^{-1} \alpha_i^+ \in \pi_1(S,*)_{n+1}$ commute with everybody
modulo $\pi_1(S,*)_{n+2}$, one obtains the  second part of (\ref{eq:goal_MJ}).
\end{proof}

Let us now recall  how the first non-vanishing Milnor invariant of a string-link, on the one hand, 
and the first non-vanishing Johnson  homomorphism of a homology cylinder, on the other hand, 
can be encoded diagrammatically. For this, we need to fix some notations.

For all set of variables $\set{l}^*=\{1^*,\dots,l^*\}$ and for all $n\geq 1$, define
$$
 D_n(\set{l}^*) := \Ker\left( (\Q \cdot \set{l}^*)  \otimes_\Q \hbox{Lie}_n(\set{l}^*) 
\stackrel{[\centerdot,\centerdot]}{\longrightarrow} \hbox{Lie}_{n+1}(\set{l}^*)\right)
$$
where $\hbox{Lie}(\set{l}^*)$ is the free Lie $\Q$-algebra generated by $\set{l}^*$ and
$\hbox{Lie}_n(\set{l}^*)$ denotes the subspace of length $n$ commutators. 
For all finite set $C$, denote
$$
\A^{t,c}(C) := \A(C) / (\hbox{ideal spanned by Jacobi diagrams that are disconnected or looped}).
$$
By writing Lie commutators as $r$-rooted binary trees whose leaves are colored with
$\set{l}^*$, we can regard $\hbox{Lie}_n(\set{l}^*)$ 
as a subspace of $\A^{t,c}_{n-1}(\set{l}^* \cup \{r\})$, the i-degree $(n-1)$ part of $\A^{t,c}(\set{l}^* \cup r\})$.
Thus, one defines an isomorphism
\begin{equation}
\label{eq:rational_iso}
\A^{t,c}_{n-1}(\set{l}^*) \stackrel{\simeq}{\longrightarrow} D_n(\set{l}^*)
\end{equation}
by sending a connected tree diagram $T$ to 
$(-1)^{n-1}\cdot \sum_{v} \hbox{color}(v) \otimes (T \hbox{ rooted at } v)$, 
where $v$ ranges over the set of external vertices of $T$.\\

\noindent
\textbf{Diagrammatic formulation of Milnor invariants \cite{HL,HM}.}
Let $(B,\gamma) \in \mathcal{S}_l[n]$. 
Because the $(n+1)$-st Artin representation of $(B,\gamma)$ 
sends $\prod_{i=1}^l x_i = \partial D_l$ to itself, the bracket of 
$$
\mu_n(B,\gamma) \in \frac{\varpi}{\varpi_2} \otimes \frac{\varpi_n}{\varpi_{n+1}}
\subset \left(\frac{\varpi}{\varpi_2} \otimes \frac{\varpi_n}{\varpi_{n+1}}\right) \otimes \Q
\simeq  (\Q \cdot \set{l}^*) \otimes_\Q \hbox{Lie}_n(\set{l}^*)
$$
vanishes. So, $\mu_n(B,\gamma)$ can be regarded in $\A^{t,c}_{n-1}(\set{l}^*)$.\\

\noindent
\textbf{Diagrammatic formulation of Johnson homomorphisms \cite{GL}.}
Let $(M,m)\in \Cyl(F_g)[n]$. Because the automorphism $m_{+,*}^{-1}\circ m_{-,*}$ of $\pi/\pi_{n+2}$
sends $\prod_{i=1}^g [\alpha_i,\beta_i] = \partial F_g$ to itself, the bracket of
$$
\tau_n(M,m) \in \frac{\pi}{\pi_2} \otimes \frac{\pi_{n+1}}{\pi_{n+2}} 
\subset \left(\frac{\pi}{\pi_2} \otimes \frac{\pi_{n+1}}{\pi_{n+2}}\right) \otimes \Q
\simeq  \left(\Q \cdot (\set{g}^+ \cup \set{g}^-) \right)
\otimes_\Q \hbox{Lie}_{n+1}(\set{g}^+ \cup \set{g}^-)
$$
vanishes. So, $\tau_n(M,m)$ can be regarded in $\A^{t,c}_{n}(\set{g}^+ \cup \set{g}^-)$.\\

Next theorem recalls how Milnor invariants of string-links can be extracted 
from their Kontsevich--LMO invariant:

\begin{theorem}[Habegger--Masbaum \cite{HM, Habegger, Moffatt}]
\label{th:HM}
Let $(B,\gamma)$ be a string-link in a homology cube with connected components labelled $1^*,\dots,l^*$.
If $\mu_{n+1}(B,\gamma) \in \A^{t,c}_{n}(\set{l}^*)$ is its first non-vanishing Milnor invariant,
then the \emph{tree-reduction} of the Kontsevich--LMO invariant of $(B,\gamma)$, namely
$$
\chi^{-1}Z(B,\gamma)\ \mod\ (\hbox{ideal spanned by looped Jacobi diagrams}) \quad \in \A^t(\set{l}^*),
$$
is equal to $\varnothing + \mu_{n+1}(B,\gamma) + (\deg > n+1)$.
\end{theorem}

\noindent
Here, the Kontsevich--LMO invariant $Z(B,\gamma)$ is as defined in \S \ref{subsec:Kontsevich-LMO},
with the assumption that $w_b(\gamma) = w_t(\gamma)$. 
We can now prove the following analogue of Theorem \ref{th:HM} for homology cylinders:

\begin{theorem}
\label{th:Johnson}
Let $M$ be a homology cylinder over $F_g$. 
If  $\tau_n(M) \in \A^{t,c}_n(\set{g}^+ \cup \set{g}^-)$ is its first non-vanishing Johnson homomorphism,
then the \emph{tree-reduction} of the LMO invariant of $M$, namely
$$
\ZtildeY(M)\ \mod\ (\hbox{ideal spanned by looped Jacobi diagrams}) \quad \in \A^{Y,t}(\set{g}^+ \cup \set{g}^-),
$$
is equal to $\varnothing +\tau_n(M) + (\ideg > n)$.
\end{theorem}

\noindent 
The proof of Theorem \ref{th:Johnson} needs the following technical

\begin{lemma}
\label{lem:chi}
Let $C$ be a finite set and let $D$ be a Jacobi diagram based on $\uparrow^C$. 
Let $D'$ be the Jacobi diagram that, as a vertex-oriented uni-trivalent graph, is $D$ 
but whose external vertices are now colored by $C$ in accordance with the way the corresponding
legs of $D$ were attached to $\uparrow^C$.
Then, one has that
$$
\chi^{-1}(D) = D' + (\hbox{i-filter} > \ideg(D)) \ \in \A(C).
$$
\end{lemma}

\begin{proof}
Let $x(D)$ be the non-negative integer $\sum_{c\in C} \max(0,|\uparrow^c \cap D|-1)$. 
If $x(D)=0$, then $\chi(D')=D$. So, we can proceed by induction on $x(D)$.
Using the STU relation, we see that $\chi(D')$ can be written as $D + \Sigma_i q_i \cdot E_i$ where
$q_i\in \Q$ and each $E_i$ is a Jacobi diagram on $\uparrow^C$ such that $x(E_i)=x(D)-1$ and $\ideg(E_i)=\ideg(D)+1$. So,
\begin{eqnarray*}
\chi^{-1}(D) & = &D' - \Sigma_i q_i \cdot \chi^{-1}(E_i) \\
&=& D' - \Sigma_i  q_i \cdot E'_i + (\hbox{i-filter}>\ideg(E_i)).
\end{eqnarray*}
Since $\ideg(E'_i)=\ideg(E_i)= \ideg(D)+1$, the conclusion follows.
\end{proof}

\begin{proof}[Proof of Theorem \ref{th:Johnson}]
Let $M$ be a homology cylinder over $F_g$ 
and let $\hbox{D}^{-1}(M)$ be the corresponding bottom-top $q$-tangle in a homology cube.
Recall that, by definition,  
\begin{equation}
\label{eq:recall}
\ZtildeY(M) = \left[- \sum_{i=1}^g \strutgraph{i^-}{i^+} \right] \sqcup
\left( \chi^{-1} Z(D^{-1}(M)) \circ \mathsf{T}_g \right).
\end{equation}
Let now $\hbox{MJ}(M)$ be the image of $M$ by the Milnor--Johnson correspondence. 
Using the functoriality and the tensor-preserving property 
of the Kontsevich--LMO invariant of $q$-tangles in homology cubes,
one easily defines from Figure \ref{fig:MJ} a linear map
$$
\Psi: \A\left(\cupright^{\set{g}^+} \capleft^{\set{g}^-}\right) 
\to \A\left(\uparrow^{\set{2g}^*}\right)
$$
such that $\Psi(Z(\hbox{D}^{-1}(N)))= Z(\hbox{MJ}(N))$ for all $N\in \Cyl(F_g)$.
Of course, this map is difficult to compute since its definition
involves several occurences of $\nu$ and of the  associator $\Phi$.
Nevertheless, one can prove the following two facts about $\Psi$:

\begin{claim}
\label{claim:about_i-degree}
For all Jacobi diagram $D \in \A\left(\set{g}^+ \cup \set{g}^-\right)$,
one has that 
$$
\chi^{-1} \Psi \chi (D) = \left[\sum_{i=1}^g \strutgraph{(2i-1)^*}{(2i)^*}\quad \quad \right] 
\sqcup \left(D\left/ \begin{array}{l} i^- \mapsto -(2i-1)^*\\ i^+ \mapsto (2i)^*\end{array}\right.\right) 
+ \left(\hbox{i-filter}> \ideg(D)\right).
$$
\end{claim}

\begin{claim}
\label{claim:about_loops}
The map $\chi^{-1} \Psi \chi $ sends a looped diagram to a linear combination of looped diagrams.
\end{claim}

\begin{proof}[Proof of Claim \ref{claim:about_i-degree}]
It is enough to prove that, 
for all Jacobi diagram $E$ based on the $1$-manifold $\cupright^{\set{g}^+} \capleft^{\set{g}^-}$,
$$
\chi^{-1} \Psi(E) = \left[\sum_{i=1}^g \strutgraph{(2i-1)^*}{(2i)^*}\quad \quad \right] \sqcup  
\left(E'\left/ \begin{array}{l} i^- \mapsto -(2i-1)^*\\ i^+ \mapsto (2i)^*\end{array}\right.\right)
+ \left(\ideg > \ideg(E)\right)
$$
where $E' \in \A(\set{g}^+ \cup \set{g}^-)$ is the diagram described in Lemma \ref{lem:chi}.
But, this follows from the facts that the Drinfeld associator is equal
to $\downarrow \downarrow \downarrow + (\hbox{i-filter} > 0)$ and that $\nu$ is equal to 
$\circlearrowleft + (\hbox{i-filter} >1)$.
\end{proof}

\begin{proof}[Proof of Claim \ref{claim:about_loops}] 
It is enough to check that the maps $\chi$, $\Psi$ and $\chi^{-1}$ 
send a looped Jacobi diagram to a linear combination of looped Jacobi diagrams. 
This is obvious for $\Psi$ and $\chi$. 
This is proved for $\chi^{-1}$ by an argument similar to the proof of Lemma \ref{lem:chi}.
\end{proof}

Let $T$ be the first non-vanishing term in the tree-reduction of $\ZtildeY(M)$
which comes, say, in i-degree $n$:
$$
\ZtildeY(M) = \varnothing + T + (\ideg> n) + (\hbox{looped diagrams}).
$$
We deduce from (\ref{eq:recall}) that $\chi^{-1} Z(\hbox{D}^{-1}(M))$ is equal to 
$$  
\left(\left[\sum_{i=1}^g \strutgraph{i^-}{i^+}\right]
\sqcup \left(\varnothing+T + (\ideg> n) + (\hbox{looped diagrams})\right) \right)\circ
\mathsf{T}_g^{-1}
$$
$$
= \mathsf{T}_g^{-1} + \left[\sum_{i=1}^g \strutgraph{i^-}{i^+}\right] \sqcup
\left(\left(T + (\ideg> n) + (\hbox{looped diagrams})\right) \star \left(\mathsf{T}_g^{-1}\right)^{Y}\right)
$$
where $\mathsf{T}_g^{-1}$ denotes the inverse of $\mathsf{T}_g$ in the algebra $\left(\tsA(g,g),\circ \right)$, 
namely $\mathsf{T}_g^{-1} = \chi^{-1} Z (\Id_w)$ 
if $w$ is an arbitrary non-associative word of length $g$.
Since $\mathsf{T}_g^{-1}$ is group-like by Lemma \ref{lem:linking_matrix}, 
$\left(\mathsf{T}_g^{-1}\right)^{Y}$ is equal to 
$\varnothing + (\ideg>0)$. We obtain that
$$
\chi^{-1} Z(\hbox{D}^{-1}(M)) = \mathsf{T}_g^{-1} +
\left[\sum_{i=1}^g \strutgraph{i^-}{i^+}\right] \sqcup T +
(\ideg> n) + (\hbox{looped diagrams}).
$$
By applying $\chi^{-1} \Psi \chi$ to this identity and by using the above two claims, we get 
$$
\chi^{-1}  Z(\hbox{MJ}(M)) = \chi^{-1} \Psi \chi\left(\mathsf{T}_g^{-1}\right)  +
\left(T\left/ \begin{array}{l} i^- \mapsto -(2i-1)^*\\ i^+ \mapsto (2i)^*\end{array}\right.\right) 
$$
$$
+ (\ideg> n) + (\hbox{looped diagrams}).
$$
Since $\chi^{-1} \Psi \chi\left(\mathsf{T}_g^{-1}\right)= \chi^{-1} \Psi Z(\Id_w) = 
\chi^{-1} Z(\hbox{trivial string-link})=\varnothing$,
we conclude thanks to Theorem \ref{th:HM} and Claim \ref{claim:MJ}.
\end{proof}

Since the Johnson filtration of the Torelli group has a trivial intersection,
Theorem \ref{th:Johnson} has the following 

\begin{corollary}
The LMO homomorphism of homology cylinders restricts to a group homomorphism
$$
\ZtildeY: \I(F_g) \to \left\{\hbox{units of } \left(\AY(\set{g}^+ \cup \set{g}^-), \star \right) \right\},
$$
which is injective.
\end{corollary}

\subsection{At the Casson invariant level}

As an illustration of the previous material, let us recover this formula
due to Morita that measures in terms of the first Johnson homomorphism how the Casson
invariant is far from defining a representation of the Torelli group \cite{Morita}.

Each homology cube $B$ has a \emph{Casson invariant} $\lambda(B)$, which
is the Casson invariant $\lambda(\hat{B})$ of the corresponding homology sphere  $\hat{B}$.

\begin{proposition}
\label{prop:Casson}
The reduction of $\Ztilde$ modulo $(\ideg >2)$ is a functorial extension 
of the Casson invariant to the category $\LCob$.
\end{proposition}

\begin{proof}
According to Proposition \ref{prop:low_degree}, the reduced LMO functor $\Ztilde$ factors through $\LCob$:
$$
\xymatrix{
{\LqCob} \ar[r]^-{\Ztilde} \ar@{->>}[d] & {\tsA}/_{ (\ideg>2)}\\
{\LCob} \ar@{-->}[ru]
}
$$
So, it is enough to recall from \cite{LMO} that, 
for all homology cube $B\in \LCob(0,0)$, one has
$$
\Ztilde(B) = \Omega(\hat{B}) = \varnothing + 
\frac{\lambda(\hat{B})}{2}\cdot \thetagraph + (\ideg > 2) \quad \in \A(\varnothing).
$$
\end{proof}

Any homology cylinder $M$ over $F_g$ can be ``filled-in'' to a homology cube:
$$
\overline{\underline{M}} := \epsilon^{\otimes g} \circ M \circ \eta^{\otimes g}
$$
where $\eta^{\otimes g}=C_g^0$ 
and $\epsilon^{\otimes g}= C_0^g$ are the genus $g$ handlebodies defined in \S \ref{subsec:cob}.

\begin{theorem}
For all $M,N\in \Cyl(F_g)$, one has that
$$
\lambda\left(\overline{\underline{M \circ N}}\right)
= \lambda\left(\overline{\underline{ M}}\right) + \lambda\left(\overline{\underline{ N}}\right)
+ 2 \cdot \left(\hbox{the $\thetagraph$-coordinate in }\tau_1(M) \star \tau_1(N)\right) \ \in \Z.
$$
\end{theorem}

\noindent
When $M,N \in \I(F_g)$, this is exactly Morita's formula \cite[Theorem 4.3]{Morita}.

\begin{proof}
By functoriality of $\Ztilde$, one has that
$$
\Ztilde(\overline{\underline{ M}})
=  \Ztilde(\epsilon^{\otimes g}) \circ \Ztilde(M) \circ \Ztilde(\eta^{\otimes g}).
$$
We deduce from Lemma \ref{lem:values} that
$$ 
\Ztilde(\overline{\underline{M}}) =
\varnothing \circ \Ztilde(M) \circ \varnothing 
= \left(\ZtildeY(M) / i^+ \mapsto 0, i^- \mapsto 0\right).
$$
In particular, the $\thetagraph$-coordinates of $\Ztilde(\overline{\underline{M}})$ and $\ZtildeY(M)$ are the same,
namely equal to $\lambda(\overline{\underline{M}})/2$.
Moreover, it follows from the identity 
$\ZtildeY(M\circ N) = \ZtildeY(M) \star \ZtildeY(N)$ that
$$
\ZtildeY_2(M\circ N) = \ZtildeY_2(M) + \ZtildeY_2(N) + \ZtildeY_1(M) \star \ZtildeY_1(N),
$$ 
where $\ZtildeY_i$ denotes the i-degree $i$ part of $\ZtildeY$.
Since $\ZtildeY_1(M) =\tau_1(M)$ and $\ZtildeY_1(N) =\tau_1(N)$, the conclusion follows.
\end{proof}

%
%
%
% Bibliography
%
%
%

\bibliographystyle{amsalpha}

\begin{thebibliography}{99}
\newcommand{\numero}{$\textrm{n}^{\circ}$}

\bibitem{BN}
D. Bar-Natan,
\emph{On the Vassiliev knot invariants},
Topology \textbf{34:2} (1995), 423--472.

\bibitem{BGRT1}
D. Bar-Natan, S. Garoufalidis, L. Rozansky, D. Thurston,
\emph{The Aarhus integral of rational homology $3$-spheres I},
Selecta Math. \textbf{8:3} (2002),  315--339.

\bibitem{BGRT2}
D. Bar-Natan, S. Garoufalidis, L. Rozansky, D. Thurston,
\emph{The Aarhus integral of rational homology $3$-spheres II},
Selecta Math. \textbf{8:3} (2002),  341--371.

\bibitem{BL}
D. Bar-Natan, R. Lawrence, \emph{A rational surgery formula for the LMO invariant},
Israel J. Math. \textbf{140} (2004), 29--60.


\bibitem{CL1}
D. Cheptea, T.T.Q. Le,
\emph{A TQFT associated to the LMO invariant of three-dimensional manifolds},
preprint (2005), \texttt{math.GT/0508220}, to appear in Comm. Math. Phys.

\bibitem{CL2}
D. Cheptea, T.T.Q. Le,
\emph{Three-cobordisms with their rational homology on the boundary},
preprint (2006), \texttt{math.GT/0602097}.

\bibitem{CY}
L. Crane, D. Yetter,
\emph{On algebraic structures implicit in topological quantum field theories},
J. Knot Th. Ramifications \textbf{8:2} (1999), 125--163.

\bibitem{Garoufalidis}
S. Garoufalidis,
\emph{The mystery of the brane relation},
J. Knot Theory Ramifications \textbf{11:5} (2002), 725--737.

\bibitem{GGP}
S. Garoufalidis, M. Goussarov, M. Polyak,
\emph{Calculus of clovers and FTI of $3$-manifolds},
Geom. Topol. $\mathbf{5}$ (2001), 75--108.

\bibitem{GL}
S. Garoufalidis, J. Levine,
\emph{Tree--level invariants of 3-manifolds, Massey products
and the Johnson homomorphism},
Proc. Sympos. Pure Math. \textbf{73} (2005), 173--203.

\bibitem{Goussarov_note}
M. Goussarov,
\emph{Finite type invariants and $n$-equivalence of 3-manifolds},
Compt. Rend. Acad. Sc. Paris $\mathbf{329}$ S\'erie I (1999), 517--522.

\bibitem{Goussarov_clovers}
M. Goussarov,
\emph{Knotted graphs and a geometrical technique of $n$-equivalence},
St. Petersburg Math. J. \textbf{12-4} (2001).

\bibitem{Habegger}
N. Habegger,
\emph{Milnor, Johnson and tree-level perturbative invariants},
preprint (2000).

\bibitem{HL}
N. Habegger, X.--S. Lin,
\emph{On link concordance and Milnor's $\overline{\mu}$ invariants},
Bull. London Math. Soc. \textbf{30:4} (1998), 419--428.

\bibitem{HM}
N. Habegger, G. Masbaum,
\emph{The Kontsevich integral and Milnor's invariant},
Topology \textbf{39} (2000), 1253--1289.

\bibitem{Habiro_claspers}
K. Habiro,
\emph{Claspers and finite type invariants of links},
Geom. Topol. \textbf{4} (2000), 1--83.

\bibitem{Habiro_BT}
K. Habiro,
\emph{Bottom tangles and universal invariants},
Algebr. Geom. Topol. \textbf{6} (2006), 1113--1214.

\bibitem{JMM}
D. Jackson, I. Moffatt, A. Morales, 
\emph{On the group-like behaviour of the Le--Murakami--Ohtsuki invariant}, 
preprint (2005), \texttt{math.QA/0511452}.

\bibitem{Kerler}
T. Kerler,
\emph{Bridged links and tangle presentations of cobordism categories},
Adv. Math. \textbf{141:2} (1999), 207--281. 

\bibitem{Kerler_towards}
T. Kerler,
\emph{Towards an algebraic characterization of $3$-dimensional cobordisms},
In: \emph{Diagrammatic morphisms and applications (San Francisco, CA, 2000)},  
141--173, Contemp. Math., 318, Amer. Math. Soc., Providence, RI, 2003.

\bibitem{Le_universality}
T.T.Q. Le,
\emph{An invariant of integral homology $3$-spheres which is universal for all finite type invariants},
In: \emph{Solitons, geometry, and topology: on the crossroad}, 
75--100, Amer. Math. Soc. Transl. Ser. 2, 179 (1997).

\bibitem{Le_Grenoble}
T.T.Q. Le,
\emph{The LMO invariant},
lecture notes given during the  \emph{Ecole d'Et\'e de Math\'ematiques : Invariants de noeuds
et de vari\'et\'es de dimension trois}, Institut Fourier (1999).

\bibitem{LM1}
T.T.Q. Le, J. Murakami,
\emph{Representation of the category of tangles by Kontsevich's iterated integral},
Comm. Math. Phys. \textbf{168:3} (1995), 535--562.

\bibitem{LM2}
T.T.Q. Le, J. Murakami,
\emph{The universal Vassiliev--Kontsevich invariant for framed oriented links},
Compositio Math. \textbf{102:1} (1996), 41--64.

\bibitem{LM3}
T.T.Q. Le, J. Murakami,
\emph{Parallel version of the universal Vassiliev--Kontsevich invariant},
J. Pure Appl. Algebra  \textbf{121:3}  (1997), 271--291.

\bibitem{LMO}
T.T.Q. Le, J. Murakami, T. Ohtsuki,
\emph{ On a universal perturbative invariant of $3$-manifolds},
Topology \textbf{37:3} (1998), 539--574.

\bibitem{Lescop}
C. Lescop,
\emph{Splitting formulae for the Kontsevich--Kuperberg--Thurston invariant
of rational homology $3$-spheres},
preprint (2004),  \texttt{math.GT/0411431}.

\bibitem{Massuyeau}
G. Massuyeau,
\emph{Finite-type invariants of three-manifolds and the dimension subgroup problem},
preprint (2006), \texttt{math.GT/0605497}, to appear in J. London Math. Soc.

\bibitem{Matveev}
S. Matveev,
\emph{Generalized surgery of three-dimensional manifolds 
and representations of homology spheres}, 
Mat. Zametki \textbf{42:2} (1987), 268--278.

\bibitem{MM}
J. Milnor, J. Moore,
\emph{On the structure of Hopf algebras},
Ann. Math. \textbf{81} (1965), 211--264.

\bibitem{Moffatt}
I. Moffatt,
\emph{The Aarhus integral and the $\mu$-invariants},
J. Knot Th. Ramifications \textbf{15} (2006), 361--377.

\bibitem{Morita}
S. Morita,
\emph{On the structure of the Torelli group and the Casson invariant},
Topology \textbf{30:4}, 603--621.

\bibitem{MN}
H. Murakami, Y. Nakanishi,
\emph{On a certain move generating link homology}, 
Math. Ann. \textbf{284} (1989), 75--89.

\bibitem{MO}
J. Murakami, T. Ohtsuki, 
\emph{Topological quantum field theory for the universal quantum invariant}, 
Comm. Math. Phys. \textbf{188:3} (1997), 501--520.

\bibitem{Ohtsuki} T. Ohtsuki, \emph{Quantum invariants. A study of knots, $3$-manifolds, and their sets},
Series on Knots and Everything $\mathbf{29}$, Word Scientific Publishing Co., 2002.

\bibitem{Stallings} 
J. Stallings, 
\emph{Homology and central series of groups},  
Journal of Algebra \textbf{2} (1965), 170-181. 

\end{thebibliography}

\end{document}